\documentclass[a4paper,twoside,12pt]{book}
\usepackage[english]{babel}
\usepackage[latin1]{inputenc}
\usepackage[T1]{fontenc}
\usepackage{amsmath}
\usepackage{amsfonts}
\usepackage{dsfont}
\usepackage{amssymb}
\title{Entropy of meromorphic maps and dynamics of birational maps}
\usepackage{dsfont}
\author{Henry de Thélin and Gabriel Vigny}
\begin{document}
\newtheorem{Theorem}{Theorem}
\newtheorem{theorem}{Theorem}[section]
\newtheorem{proposition}[theorem]{Proposition}
\newtheorem{defi}[theorem]{Definition}
\newtheorem{corollaire}[theorem]{Corollary}
\newtheorem{Hypothesis}[theorem]{Hypothesis}
\newtheorem{lemme}[theorem]{Lemma}
\newtheorem{Remark}[theorem]{Remark}
\newcommand{\U}{\mathcal{U}}
\newcommand{\C}{\mathcal{C}}
\renewcommand{\P}{\mathbb{P}}
\newcommand{\Cc}{\mathbb{C}}
\newcommand{\Nn}{\mathbb{N}}
\newcommand{\Rr}{\mathbb{R}}
\newcommand{\Qq}{\mathbb{Q}}
\newcommand{\Zz}{\mathbb{Z}}
\newcommand{\Acal}{\mathcal{A}}
\newcommand{\Bcal}{\mathcal{B}}
\newcommand{\Dcal}{\mathcal{D}}
\newcommand{\Ecal}{\mathcal{E}}
\newcommand{\Gcal}{\mathcal{G}}
\newcommand{\Hcal}{\mathcal{H}}
\newcommand{\Lcal}{\mathcal{L}}
\newcommand{\Ical}{\mathcal{I}}
\newcommand{\Pcal}{\mathcal{P}}
\newcommand{\Qcal}{\mathcal{Q}}
\newcommand{\Scal}{\mathcal{S}}
\newcommand{\Zcal}{\mathcal{Z}}
\date{}
\frontmatter

\maketitle

\vspace{0.5 cm}

\noindent \emph{Abstract.} \\

\noindent We study the dynamics of meromorphic maps for a compact Kähler 
manifold $X$. More precisely, we give a simple criterion that allows us 
to produce a measure of maximal entropy. We can apply this result to bound the Lyapunov exponents. 

 Then, we study the particular case of a family of generic birational
 maps of $\P^k$ for which we construct the Green currents and the equilibrium measure. We use for that the theory of
 super-potentials. We show that the measure is mixing and gives no
 mass to pluripolar sets. Using the criterion we get that the measure
 is of maximal entropy. It implies finally that the measure is hyperbolic. 
 
\vspace{1 cm} 

\noindent\textbf{MSC: 37Fxx, 32H04, 32Uxx, 37A35, 37Dxx}   \\
\noindent\textbf{Keywords: Complex dynamics, meromorphic maps,
  Super-potentials, currents, entropy, hyperbolic measure. }

\tableofcontents

\mainmatter

\chapter{Introduction}
Complex dynamics in several variables and more precisely the iteration 
of polynomial maps have received much attention in the last twenty years. 
This can be explained because of the links with real dynamics 
(especially for \emph{Hénon maps}) and also because of the possibility to 
use powerful methods from several complex variables. 

Let $P_1,\dots,P_k$ be polynomials in $k$ complex variables and let
$f=(P_1,\dots,P_k)$ be the associated polynomial map in $\mathbb{C}^k$. 
The issue is to study the behavior of the sequence of iterations 
$f^n$. 
As such, it is often easier to consider the dynamics 
in $\mathbb{P}^k$ instead of $\mathbb{C}^k$. So we study the dynamics of rational maps $f$ in $\P^k$ and 
more generally the dynamics of dominating meromorphic maps 
in a compact Kähler manifold $X$ (recall that a map is \emph{dominating} 
if its image contains an open set). The ``classical'' program is to 
construct an invariant measure that will describe the chaotic part 
of the dynamics. Then one tries to prove the basic properties of 
the measure: ergodicity or even mixing, computation of the entropy 
(with the question: is the measure of maximal entropy?) and estimation
of the Lyapunov exponents (or simply a bound). Essentially, 
one want to prove that the measure of maximal entropy is hyperbolic. \\

In dimension 1, a classical tool is Montel theorem: a family of maps from the unit disk $\mathbb{D}\subset \mathbb{C}$ 
to $\P^1$ minus three points is normal. 
There are no such simple results in higher dimension so one need
to use other techniques. For endomorphisms of $\P^k$, the measure was 
defined by Forn{\ae}ss and Sibony in \cite{FS, FS1, F.S} using pluripotential theory.  
They introduced a positive closed current of bidegree $(1,1)$ 
called the \emph{Green current} which carries informations on the 
dynamics of $f$, especially on its chaotic behavior 
(see also \cite{Sib}). 
Then the measure is
 defined as a Monge-Ampère of the Green current
and the authors show that it is mixing.
Briend and Duval gave a bound of the Lyapunov exponents and showed that the measure is hyperbolic
(see \cite{B, BDu1, BDu}). 

Pluripotential has also been used for Hénon maps. Sibony defined the 
Green current for Hénon maps and then the equilibrium measure as an intersection 
of the Green currents. Using pluripotential theory, Bedford, 
Lyubich and Smillie proved numerous properties of the currents and measures 
in a series of articles (see in particular \cite{B.S, BLS, BS}), 
see also the results of Forn{\ae}ss and Sibony in \cite{FS2}. \\
 
In order to define the currents and the measures, 
one has to deal with some \emph{dynamical degrees} $d_1,\dots, d_k$ (see for example \cite{RS}). 
Roughly speaking, the degree $d_l$ measures the asymptotic spectral 
radius of the action of $f$ on the cohomology group $H^{l,l}(X)$. 
The last degree $d_k$ is the \emph{topological degree}. It can be shown 
that the sequence of degrees is increasing up to a rank $s$ and then it is decreasing. 
When several dynamical degrees are equal, complications might happen 
and the program fails (see \cite{Gu1}). So the study takes place when 
there is a dynamical degree $d_s$ strictly larger than the others. 
When $s=k$, namely the topological degree is the largest dynamical degree, 
one can construct and study the measure directly (see \cite{Gu2}, \cite{DS3}). 
The other cases are harder and one has often to make additional hypotheses.\\
 
Another complication appears with indeterminacy sets in particular the \emph{second indeterminacy set}:
 $$ I':=\{z , \ \text{dim}(f^{-1}(z)) \geq 1 \}.$$ 
This set is of codimension $\geq 2$ thus is of mass zero for a positive closed currents of bidegree $(1,1)$. The presence of those indeterminacy sets implies difficulties to define pull-back, push-forward and intersections of currents and measures. So here again, one has to make some hypothesis on the indeterminacy sets to define those operations. Finally, when $s>1$, one has to deal with currents of bidegree $(s,s)$: the potentials of those currents are no longer quasi-plurisubharmonic functions but forms that can be singular. Consequently, very little has been done in the study of meromorphic maps in dimension $k>2$ for which the largest dynamical degree is not the topological degree. \\

On the other hand, the abstract theory of dynamical systems and
especially non uniformly hyperbolic dynamical systems is very
developed with the work of Yomdim, Pesin, Katok and others
(\cite{K.H}). Assuming the existence of an invariant measure, one can
define the (metric) entropy of the measure which describes how chaotic
the dynamics is. When the map is continuous, the variational principle
implies that the topological entropy of the map is given by the
supremum of the entropies for all the invariant measures. Moreover,
when one has a hyperbolic measure, the Oseledec-Pesin's theory permits to construct stable and unstable manifolds associated to non zero Lyapunov exponents and we have uniform estimates outside sets of small measure.

This describes fairly well the dynamics of the application. In a way,
such behavior is expected to be generic. A central and difficult
problem in dynamics is to construct examples of hyperbolic invariant measures. 
In the complex case, this can be done for holomorphic maps 
(see also polynomial-like maps \cite{DS3} and horizontal-like maps \cite{DNS}). 
So there is a need for dynamical models which admits a hyperbolic measure.  \\

The purpose of this study is to answer the two above questions:
getting results on the dynamics of meromorphic maps in general and
giving classes of examples where one can prove the hyperbolicity of an
invariant measure. More precisely, in a first part (Chapter
\ref{Henry}), we give a criterion that allows us to produce invariant
measure of maximal entropy for a meromorphic map of a compact Kähler
manifold $X$. This can then be applied to bound the Lyapunov
exponents. In a second part (Chapter \ref{Gabriel}), we study the more
precise case of a generic family of birational maps of $\P^k$ for
which we construct the equilibrium measure. We show that it is mixing
and using the results of the first part we show that it is of maximal
entropy. We deduce finally the hyperbolicity of the measure. Let us detail our results. \\

When $f: X \rightarrow X$ is a smooth map on a smooth Riemannian manifold
 $X$, it is known since the work of  Yomdin (see \cite{Y} and \cite{G}) and  Newhouse (see \cite{Ne}) that $f$
admits an ergodic measure of maximal entropy. If $f$ is a Hénon map of $\Cc^2$,
E. Bedford and J. Smillie have shown in \cite{BS} that 
the Green measure of $f$ is of maximal entropy. Their proof is based on  
Yomdin's theorem (see \cite{Y}) and also on the proof of the variational 
principle. This approach has been used several times since then in 
dynamics in order to bound from below the entropy of measures (e.g. 
\cite{Gu}, \cite{De1} and \cite{DS12}). In all these cases, one can use 
 Yomdin's theorem because the application $f$ is either holomorphic or 
 when it is meromorphic everything takes place in a stable open set where $f$ is holomorphic. \\

The purpose of the first part is to quantify Bedford
and Smillie's approach. We detail the setting first.

Let $(X, \omega)$ be a compact Kähler manifold of dimension $k$ and let $f$ be a 
dominating meromorphic map. We denote by $I$ the indeterminacy set of $f$ and for $l=0 \dots k$, we write:
$$\lambda_l(f):=\int_X f^{*}(\omega^l) \wedge \omega^{k-l}.$$
The $l$-th \emph{dynamical degree} of $f$ is defined by (see \cite{RS} and \cite{DS5}):
$$d_l:= \lim_{n \rightarrow + \infty} (\lambda_l(f^n))^{1/n}.$$
Now, we consider the sequence of measures:
$$\mu_n:= \frac{1}{n} \sum_{i=0}^{n-1} f^i_{*}\left( \frac{(f^n)^* \omega^{l}
  \wedge \omega^{k-l}}{\lambda_l(f^n)} \right).$$
 It is a well defined sequence of probability measures (see Section \ref{merom}).
  Remark that in the cases where we know how to construct a measure of maximal entropy $\mu$, the measure $\mu$ is the limit of
$\mu_n$ with $l=s$ where $d_s$ is the largest dynamical degree (see
 \cite{BS} for Hénon maps,  \cite{Sib} for regular automorphisms of $\Cc^k$, \ldots ).
  In fact, it is likely that in the case where $d_s$ is strictly 
  larger than the other dynamical degrees then $(\mu_n)$ will always 
  converge to the measure of maximal entropy (see also \cite{Gu1}).\\

In Chapter \ref{Henry}, we do not assume that one of the dynamical degree is larger than the others. 
We suppose that there exists a subsequence
$\mu_{\psi(n)}$ of $\mu_n$ which converges to a measure $\mu$ with:
$$(H) \mbox{  :   }  \lim_{n \rightarrow + \infty}
\int \log d(x,I)  d \mu_{\psi(n)} (x) = \int \log d(x,I)  d \mu(x) > - \infty.$$
Here $d$ is a distance in $X$ and $I$ is the indeterminacy set of $f$. When $I = \emptyset$, we define $d(x,I):=1$ for all $x \in X$.

The hypothesis allows us to measure the way that the orbits get near
$I$. Then, using a quantitative version of Yomdin's theorem
and of the variational principle, we obtain a bound for the entropy:
\begin{Theorem}\label{ENTROPY}
If the Hypothesis (H) is satisfied, then  $\mu$ is an invariant measure 
of metric entropy greater or equal to $\log d_l$.
\end{Theorem}
This result is interesting even in the holomorphic case. 
Indeed, in that situation  $I$ is empty so (H) is satisfied.
So we have measures of maximal entropy: we just take the sequence
 $\mu_n$ with $l=s$ where $d_s$ is the highest dynamical degree and we take a cluster value.
Indeed, we always have the bound from above of the entropy 
by $\log d_s$ (see \cite{DS9} for the projective case and \cite{DS5} for the Kähler case).
More generally, if we can prove the convergence of the sequence $\mu_n$ 
with the hypothesis (H) with $l=s$ where $d_s$ is the highest dynamical degree, 
we obtain for the same reason explicit measures of maximal entropy $\log d_s$.

Remark that the criterion can be extended to the case where $(X,
\omega)$ is a compact hermitian manifold. In that case, we do not know if the limit:
$$d_l:= \lim_{n \rightarrow + \infty} (\lambda_l(f^n))^{1/n}$$
exists, but it is sufficient to replace $d_l$ by $\limsup_n
(\lambda_l(f^{\psi(n)}))^{1/\psi(n)}$ in the theorem.

Under an additional hypothesis on the integrability of $\log d(x,\C)$
where $\C$ is the critical set and when $d_s$ is strictly larger than
the other dynamical degrees we can use a result of the first author to give a bound on the Lyapunov exponents \cite{DT1} which
implies the hyperbolicity of the measure. \\

In Chapter \ref{Gabriel}, we study the dynamics of birational maps of
$\P^k$, that is maps that are meromorphic and biholomorphic outside
some analytic set. The study of birational maps started with in $\P^2$
with the dynamics of Hénon maps. For such a map $f$ of algebraic
degree $d$, Sibony introduced  the Green current $T^+$ as
$T^+=\lim_{n\to \infty} d^{-n}(f^n)^*(\omega)$ (here $\omega$ is the
Fubini-Study form on $\P^2$). The limit exists and for the same
reasons we can consider the current $T^-$ associated to
$f^{-1}$. Sibony's strategy is then to consider the measure $\mu:=T^+
\wedge T^-$ (well defined) and to prove the ergodic properties of the
measure (mixing, entropy, \ldots). This has been done for polynomial
automorphisms of $\mathbb{C}^2$ by Bedford, Smillie and Lyubich
(\cite{BS, B.S, BLS}) and also Forn{\ae}ss and Sibony in
\cite{FS2}. This strategy has been used for different families of
birational maps of surfaces (see for example \cite{Dil} and
\cite{DF}). Each time, the properties of the potential of those
currents play a big role to prove the existence of measures.

 Sibony  worked out these properties in the case of regular
 automorphisms of $\mathbb{C}^k$ (\cite{Sib}, \cite{GS1} and also
 \cite{Gu}). Sibony and Dinh extended these results to the case of
 \emph{regular birational maps} in $\P^k$ in \cite{DS10}. One can also
 study the dynamics of automorphism of compact Kähler manifolds (see
 \cite{Can}, \cite{DS12}). 
 
 In all the above works, the indeterminacy sets of $f$ and
 $f^{-1}$ are either empty (for automorphisms of compact Kähler
 manifold) or are disjoint from the support of the equilibrium
 measure. Roughly speaking, the cases considered by these authors satisfy the condition:
$$ \overline{\bigcup_{n\geq 0} f^{-n}I^+} \cap \overline{\bigcup_{n\geq 0} f^{n}I(f^-)} =\varnothing,$$
where $I^+$ is the indeterminacy set of $f$ and $I^-$ is the indeterminacy set of $f^{-1}$. \\
 
 Another approach in the case of surfaces, initiated by Bedford and Diller in \cite{BD1}, is to take a weaker, quantitative version of the above, namely:
 $$\sum_{n \geq 0} \left(\frac{1}{d}\right)^n\log \text{dist}(I^+,f^n(I^-)) >-\infty. $$
Using that hypothesis, the authors define the equilibrium measure and show that the potential of the Green current is integrable for the measure. They proved that the measure is mixing and hyperbolic. Using laminar currents, Dujardin computed the entropy and showed that the measure is of maximal entropy \cite{Duj}.  Diller and Guedj extended those results to a more general case in \cite{Dilgu}. Note also the extension to the case of meromorphic maps of a surface in the recent articles \cite{DDG1}, \cite{DDG2}, \cite{DDG3}.  \\

Here we explore both directions. We consider birational maps of $\P^k$ ($k\geq2$) and we authorize the indeterminacy sets to get close to each others. Namely, let $f:\P^k\rightarrow\P^k$ be a birational map of algebraic degree $d$ and let $\delta$ be the algebraic degree of $f^{-1}$. We assume that $\text{dim}(I^+)=k-s-1$ and $\text{dim}(I^-)=s-1$; in this case, we have that the largest dynamical degree of $f$ is $d^+_s=d^s$ and  the largest dynamical degree of $f^{-1}$ is $d^-_{k-s}=\delta^{k-s}=d^s$. We assume that: 
 \begin{align*}
&\sum_{n \geq 0} \left(\frac{1}{d}\right)^n\log \text{dist}(I^+,f^n(I^-)) >-\infty \\
& \qquad \qquad \qquad \text{and} \\
&\sum_{n \geq 0} \left(\frac{1}{\delta}\right)^n\log \text{dist}(I^-,f^{-n}(I^+)) >-\infty.
\end{align*}
In fact, we will assume a weaker hypothesis (which is equivalent to the previous one only in dimension $2$). 
The interest of the family of maps that we consider is that they are generic (see Theorem \ref{generic}).

Under that condition, we construct the Green current $T^+_s$ of order $s$ of $f$. 
Similarly, we define the Green current $T^-_{k-s}$ of order $k-s$. More precisely, we have (see Theorems \ref{convergencecurrent}, \ref{invariance} and \ref{extremality}):
\begin{Theorem}
Let $f$ be a birational map as above, then the sequence
$(d^{-s}f^*)^n(\omega^s)$ is well defined and converges in the sense of currents to a positive closed current $T^+_s$ of bidegree $(s,s)$ and of mass $1$. 

The current $T^+_s$ satisfies $f^*(T_s^+)=d^s T^+_s$ and is extremal in the set of positive closed currents.   
\end{Theorem}
We prove some equidistribution results
on the currents.
Then we consider the intersection $T^+_s \wedge T^-_{k-s}$ and we prove (Theorem \ref{constructionmeasure}, Proposition \ref{invariancemeasure} and Theorem \ref{mixing}):
\begin{Theorem}
The wedge-product $\mu:=T^+_s \wedge T^-_{k-s} $
is a well-defined invariant probability measure for which  the potential of the Green current of order $1$ is integrable. The measure $\mu$ is mixing for $f$.
\end{Theorem}
Using a space of test functions introduced by Dinh and Sibony in \cite{DS4} and studied by the second author \cite{moi2}, we prove that the measure gives no mass to pluripolar sets. In particular, the measure gives no mass to analytic subsets.   

Then we use the results of Chapter \ref{Henry} to prove that (Theorem \ref{entropy} and Theorem \ref{hyperbolic}): 
\begin{Theorem}
The measure $\mu$ is of maximal entropy $\log d^s$ and is hyperbolic.
\end{Theorem}
In order to prove the convergences, we deal directly with positive closed currents of bidegree $(s,s)$. The potentials $U$ of a positive closed current $S$ of bidegree $(s,s)$ are no longer quasi-plurisubharmonic (qpsh for short) functions but currents satisfying $dd^cU+\omega^s=S$. Two such potentials $U$ and $U'$ differ by a $dd^c$ closed current. Such object can be singular. So we use the new theory of super-potentials introduced by Dinh and Sibony \cite{DS6} (and also \cite{DS12} for the Kähler case). It provides a calculus on $(s,s)$ positive closed currents.

The idea is to consider super-potentials $\U$ of $S$ not as a form of bidegree $(s-1,s-1)$ 
but as a function on positive closed currents of bidegree $(k-s+1, k-s+1)$. 
Super-potentials can be seen as qpsh functions on the set of positive closed 
currents of bidegree $(k-s+1, k-s+1)$ and they inherit the properties of qpsh functions. 
 
We sum up the properties of super-potentials that we used in an appendix. \\

The two parts are fairly independent as we only use the results 
of Chapter \ref{Henry} at the end of Chapter \ref{Gabriel}. 
So they can be read in any order.

\chapter{Entropy of meromorphic maps}\label{Henry}
\section{Push-forward of measures by meromorphic maps}\label{merom}
Let $(X,\omega)$ be a compact Kähler manifold of dimension $k$. We assume that the diameter of $X$ is less than $1$. 
Let $f$ be a dominating meromorphic map and let $I$ be the indeterminacy set of $f$.
Recall that for $l=0 \dots k$, we write:
$$\lambda_l(f):=\int_X f^{*}(\omega^l) \wedge \omega^{k-l}.$$

We start by recalling how to define the push-forward by $f$
of a measure that gives no mass to $I$. In all this text, a measure
will be a finite positive Radon measure.

Let $\nu$ be such a measure. On $X \setminus I$, $f$ is a measurable map. So we can define
 $f_{*} \nu$ by the formula:
$$(f_{*} \nu)(A):=\nu(\{x \in X \setminus I \mbox{ with } f(x)
\in A \})= \nu( f^{-1}(A) \cap (X \setminus I)).$$

When a measure $\nu$ gives no mass to the indeterminacy set, we have:
$$ \int \varphi \circ f d \nu= \int \varphi d (f_{*} \nu)$$
for all $\varphi \in L^1(f_{*} \nu)$. It is implicitly assumed that the integral is on $X
\setminus I$. The equality follows from the approximation of function in $L^1$ by
characteristic functions.

The operator $f_{*}$ has the good property of continuity. Indeed, we have:
\begin{lemme}{\label{continuite}}
Let $\nu_n$ be a sequence of measures that give no mass to $I$. Then if $(\nu_n)$ converges to $\nu$ and
$\nu(I)=0$ then $(f_{*}(\nu_n))$ converges to $f_{*}\nu$.
\end{lemme}
\emph{Proof.} Of course, the mass of $\nu_n$ converges to the mass of $\nu$. 
Now, let $\varphi$ be a continuous function and let 
$0\leq \chi_{\varepsilon} \leq 1$ be a smooth function equal to $0$ in
an $\varepsilon$-neighborhood $I_{\varepsilon}$ of $I$ and equal to $1$ outside a
$2\varepsilon$-neighborhood $I_{2 \varepsilon}$ of $I$. Then, we have:
$$\int \varphi  d (f_{*} \nu_n)= \int \varphi \circ f d \nu_n= \int (1-
\chi_{\varepsilon}) \varphi \circ f d \nu_n +  \int  \chi_{\varepsilon} \varphi
\circ f d \nu_n.$$
The first term is bounded in absolute value by  $\| \varphi \|_{\infty} \nu_n( I_{2\varepsilon})$
which can be taken arbitrarily small by taking $\varepsilon$
small then $n$ large (because $\nu$ gives no mass to $I$). The second term
  converges to $ \int  \chi_{\varepsilon} \varphi
\circ f d \nu$ since $\chi_{\varepsilon} \varphi
\circ f$ is a continuous function. Finally, if $\varepsilon$ is small enough,
 $ \int  \chi_{\varepsilon} \varphi
\circ f d \nu$ is as close as we want from $ \int  \varphi \circ f d \nu$ since
 $\nu$ gives no mass to $I$.
  \hfill $\Box$ \hfill \\

In this section, we consider in particular the push-forward of the measures 
$$\nu_n:=\frac{(f^n)^* \omega^{l}
  \wedge \omega^{k-l}}{\lambda_l(f^n)}.$$
The $\nu_n$ are well defined probability measures. Indeed, $(f^n)^*
\omega^{l}$ is a form with coefficients in $L^1$ so it gives no mass to analytic sets of dimension $<k$.
This implies that 
$$\frac{(f^n)^*\omega^{l} \wedge \omega^{k-l}}{\lambda_l(f^n)}$$
is a probability that gives no mass to
 $\cup_{i \in \Nn} f^{-i}(I)$ (because $f$ is dominating). So we can 
push-forward this probability by $f^i$ and we get again a probability. 
We also make the observation:
$$(f^i)_{*}(f^j)_{*} \frac{(f^n)^*
  \omega^{l} \wedge \omega^{k-l}}{\lambda_l(f^n)}=(f^{i+j})_{*} \frac{(f^n)^*
  \omega^{l} \wedge \omega^{k-l}}{\lambda_l(f^n)},$$
since $(f^j)_{*} \frac{(f^n)^*
  \omega^{l} \wedge \omega^{k-l}}{\lambda_l(f^n)}$ puts no mass on
  analytic sets of dimension $< k$. In particular, we can write
  $f_{*}^i \nu_n$ or $(f^i)_{*} \nu_n$, it is the same.

We also have the notion of invariance. Namely, a measure $\mu$ 
that gives no mass to $I$ is \emph{invariant} (or  $f_{*}$-invariant) if  $f_{*}(\mu)= \mu$.
One has the following easy lemma:
\begin{lemme}{\label{inv}}
Let $\mu$ be a measure that gives no mass to $I$. Then the following properties are equivalent:
\begin{itemize}
\item $\mu$ is invariant.
\item For any continuous function $\varphi$, we have:
$$\int \varphi \circ f d \mu= \int \varphi d \mu$$
where the left-hand side integral is taken over $X \setminus I$.
\end{itemize}
\end{lemme}
When these properties are true, we even have:
$$\int \varphi \circ f d \mu= \int \varphi d \mu$$
for any $\varphi$ in $L^1(\mu)$ (with the same abuse of notation for the left-hand side integral that we will do in the whole section).\\

We give now some properties of meromorphic maps that will be useful in the proof of 
Theorem \ref{ENTROPY}. First recall that we denote:
$$\mu_n:= \frac{1}{n} \sum_{i=0}^{n-1} f^i_{*}\left( \frac{(f^n)^* \omega^{l}
  \wedge \omega^{k-l}}{\lambda_l(f^n)} \right).$$
We have seen that it is a well defined sequence of probabilities. Since $f$ est dominating,
these measures give no mass to analytic sets of dimension $<k$. \\

We need an invariant measure to consider the metric entropy. So we will need the following lemma:
\begin{lemme}{\label{invariante}}
If  $(\mu_{\psi(n)})$ converges to a measure $\mu$ that gives no mass to $I$, then $\mu$ is $f_*$-invariant.
\end{lemme}
\emph{Proof.} To simplify the notations, assume that $(\mu_n)$ converges to $\mu$.

We can write $f_{*} (\mu_n)= \mu_n + \alpha_n$ with $\alpha_n$ going to zero.
Using Lemma \ref{continuite},
$f_{*}(\mu_n)$ converges to $f_{*} \mu$ and the lemma follows.  \hfill $\Box$ \hfill \\

Now, since we have an invariant probability measure that gives no mass to $I$, its mass is $1$ on $\Omega= X \setminus  \cup_{i \in \Nn}
f^{-i}(I)$. Since $f(\Omega) \subset \Omega$, we can define the metric entropy of $\mu$
using partitions (see \cite{Gu1} and \cite{K.H}).

We recall the following estimate that we use later:
\begin{lemme}{\label{estimee}}(see \cite{DiDu} Lemma 2.1)

\noindent There exist constants $K$ and $p$ such that:
$$\| D f(x) \|  \leq K d(x, I)^{-p}.$$
\end{lemme}

\section{Yomdin's theorem}{\label{Yomdin}}
In this paragraph, we recall some facts 
on Yomdin's theorem (see \cite{Y}) using Gromov's version (see \cite{G}).

Let $l$ be an integer between $1$ and $2k$. If $Y$ is a subset of $\Cc^k$
 (for example a submanifold of real dimension $l$), we call $C^r$-\emph{size}  (with $r \in \mathbb{N}^*$) of $Y$, the lower bound of the numbers $t \geq 0$ for which there exists a $C^r$-map of the unit $l$-cube into $\Cc^k$,
$h: [0,1]^{l} \mapsto \Cc^k$,  with $Y \subset
h([0,1]^{l})$ and $\| D_r h \| \leq t$. Here $ D_r h $ is the vector of the partial derivatives of $h$ of order $1 ,
\dots , r$. The norm refers to supremum over $x\in [0,1]^{l}$:
$$ \| D_r h \|= \sup_{x}  \|D_r h(x)\| $$ 
We make some comments on $C^r$-size first.

First, the $C^1$-size bounds the (real) $l$-dimensional volume of
$Y$ and its diameter. More precisely
$$ C^1-\text{size of} \ Y \geq \max((
\text{l-dimensional volume }(Y))^{1/l}, l^{-1/2} \mbox{Diameter}(Y)).$$
A process that we will use in what follows is the division of a set of
$C^r$-size. If $Y$ is a set of $C^r$-size smaller than $t$, we can divide
$Y$ in $j^l$ pieces of $C^r$-size smaller than $t/j$. For that it is 
sufficient to divide the $l$-cube $[0,1]^{l}$  in $j^l$ equal pieces and
then to scale: for example $R: [0,1]^{l} \mapsto [0 , j^{-1}]^{l}$ and 
similarly for the  $j^{l} -1$ other cubes. The composition of $h: [0,1]^{l}
\mapsto \Cc^k$ which covers $Y$ with the scaling $R$ satisfies 
$\| D_r (h \circ R) \|  \leq t/j$ and the union of the images of these  
$j^{l}$ maps covers $Y$.

Here is now the principal result of Gromov-Yomdin that we will need
(see Lemma 3.4 in \cite{G}).
\begin{theorem}[\cite{G}]
Let $Y$ be an arbitrary subset in the graph $\Gamma_g \subset
[0,1]^{l} \times \Cc^k$ of a $C^r$-map $g: [0,1]^l \mapsto \Cc^k$
and take some positive number $ \epsilon \leq 1$. Then $Y$ can be divided into
$N \leq C(k,l,r) \epsilon^{-l} (1 + \| \partial_r g \|
)^{l/r}$ sets of $C^r$-size  $\leq C(k,l,r) \epsilon
\mbox{Diameter}(Y)$, where $\partial_r g$ denotes the vector 
of the partial derivatives of $g$  of order exactly $r$ and $C(k,l,r)$ is a universal constant.
\end{theorem}
Here is the application of the above theorem that we will use:
it is a small variation of Corollary 3.5 in \cite{G}.
\begin{proposition}{\label{Gromov}}
Let $V$ be an open set of $\Cc^k$ and $f: V \rightarrow \Cc^k$ a map of class $C^r$. 
Let $Y_0 \subset V$ be a set of $C^r$-size smaller than $1$
such that $d(Y_0 , \partial V) \geq \sqrt{l}$. Then the intersection
of $f(Y_0)$ with a ball of  $\Cc^k$ of radius $\beta$ can be divided into 
$N \leq C(k,l,r) \left(1+  \frac{\| D_r f \|}{\beta} \right)^{l/r}$
pieces of $C^r$-size less than $\beta$.
\end{proposition}
\emph{Proof.} We want to divide $f(Y_0) \cap B(a, \beta)$ into pieces of $C^r$-size $\leq \beta$. 
If $H(a, 1 / \beta)$ denotes the homothety  of center $a$ and ratio $1 /
\beta$ in $\Cc^k$, it is equivalent to divide 
$$H(a, 1 / \beta)(f(Y_0) \cap B(a, \beta)) = H(a, 1 / \beta)(f(Y_0)) \cap
B(a,1)$$
into sets of $C^r$-size less than $1$.

By hypothesis, there exists a map $h : [0,1]^l \rightarrow \Cc^k$ of class $C^r$ with
$\| D_r h \| \leq1$ and $Y_0 \subset h( [0,1]^l )$. Define $g:=H(a , 1 /
\beta ) \circ f \circ h$. By the chain rule, we have 
$$\| D_r g \|
\leq \frac{C'(k,l,r)}{\beta} \| D_r f \|.$$
We apply now the previous theorem to $Y$ the graph of $g$ intersected with 
$[0,1]^l \times B(a,1)$. So we have that we can cover $Y$ by a number:
$$N \leq C(k,l,r)\left(1+ \frac{C'(k,l,r)}{\beta} \| D_r f \|\right)^{l/r} \leq
C(k,l,r) \left(1+ \frac{\| D_r f \|}{\beta} \right)^{l/r}$$
sets of $C^r$-size  $\leq 1$ (changing the constant $C(k,l,r)$ if necessary). Since the image of $Y$ by the projection $[0,1]^l
\times \Cc^k \mapsto \Cc^k$ covers  $H(a, 1 / \beta)(f(Y_0)) \cap
B(a,1)$, the proposition follows.  \hfill $\Box$ \hfill

\section{Proof of Theorem \ref{ENTROPY}}
The hypothesis we made assure us that there exists a subsequence
$(\mu_{\psi(n)})$ which converges to a measure $\mu$ with:
$$(H) \mbox{  :   }  \lim_{n \rightarrow + \infty}
\int \log d(x,I)  d \mu_{\psi(n)} (x) = \int \log d(x,I)  d \mu(x) > -
\infty.$$

In order to clarify the exposition, we shall write $\psi(n)=n$ .

When $s(x)$ is a function defined on $X$ with values in $\Rr^{+}$, we define (see \cite{M}):
$$B(x, s , n , f):=\left\{y \mbox{, } d(f^i(x),f^i(y)) \leq
s(f^i(x)) \mbox{ for } i \in [0,n-1] \right\}.$$
We shall use these dynamical balls with for $s(x)$ the function $\rho(x)$ or $\eta(x)$ where :
$$\rho(x)= \left( \frac{ d(x,I) \times \dots \times
  d(f^{m-1}(x),I)}{K^m} \right)^{p}$$
(here $K$ and $p$ are defined at the end of Section \ref{merom} and $m
\in \Nn$ will be chosen later) and:
$$\eta(x)= \left( \frac{ d(x,I)}{K} \right)^{p}.$$

When $f$ is holomorphic (i.e. $I = \varnothing$), take $d(x,I):=1$ and
$p=2$ in these expressions.

If $n \in \Nn$ is fixed, by the Euclidean algorithm, we write $n= \phi(n)m + r(n)$ with $0
\leq r(n) < m$. In what follows, we will consider the following dynamical balls:
$$B_n(x):= B(x, \rho ,  \phi(n), f^m) \cap f^{-\phi(n)m +
  m}(B(f^{\phi(n)m -m}(x) , \eta, r(n)+m,f)).$$

Here is the plan of the proof. As in the article of Bedford and Smillie (see \cite{BS}), 
it is based in one hand on Yomdin's theorem and in the other hand on the proof 
of the variational principle. Of course, we  have to quantify precisely those two 
parts because of the presence of the indeterminacy set.
More precisely, in a first section we show that 
there are ``a lot'' of dynamical balls. Indeed, thanks to a quantification of Yomdim's
theorem, we bound from above: 
$$\nu_{n}(B_n(x)) : =\frac{(f^{n})^{*} \omega^{l}  \wedge \omega^{k-l}}{\lambda_l(f^{n})}(B_n(x))$$
by $d_l^{-n}$ for generic points of $\nu_n$ (that we call \emph{good points}). In a second section,
we show that the presence of these $d_l^{n}$ dynamical balls allows us to bound from below 
the entropy of the measure $\mu$. We use for that ideas that lie in the 
proof of the variational principle.

\subsection{Upper bound of $\nu_n(B_n(x))$}
We give some notations first. First of all we can put on $X$ 
a family of chart $(\tau_x)_{x \in X}$ such that $\tau_x(0)=x$, $\tau_x$ is defined
on $B(0,\epsilon_0) \subset
\Cc^k$ with $\epsilon_0 > 0$ independent of $x$ and such that the norm
of the derivatives of order $1$ of the $\tau_x$ is bounded from above by a constant independent of $x$.
These charts are obtained from a finite family $(U_i, \psi_i)$ of charts of $X$ by composing them 
with translations. In $\Cc^k$, we also consider $\pi_1 \mbox{, } \dots \mbox{, }
\pi_i$ the projections from $\Cc^k$ onto 
the vectorial subspaces of dimension $k-l$. In what follows, the choice of these
coordinates is supposed to be generic and $\beta_j$ denotes the standard volume form on $\pi_j(\Cc^k)$. \\

Fix $x \in X$ and:
$$\Omega:= (\tau_x)_{*}( \pi_1^{*} \beta_1 + \dots + \pi_i^{*}
\beta_i).$$
We want to compute :
$$\nu_n(B_n(x)) =\frac{(f^n)^*  \omega^{l}  \wedge \omega^{k-l}}{\lambda_l(f^n)}(B_n(x)).$$
Taking $K$ large enough, we can assume that $B_n(x)\subset \tau_x(B(0,\epsilon_0))$ so the previous quantity is less than:
$$C(X) \frac{(f^n)^*  \omega^{l}  \wedge \Omega}{\lambda_l(f^n)}(B_n(x))=C(X)  \sum_{j=1}^{i} \int \int_{B_n(x) \cap \tau_x(Y_j(t))} \frac{(f^n)^*
  \omega^{l}}{\lambda_l(f^n)} dt$$
where $Y_j(t)$ is equal to $\pi_j^{-1}(t)$ for $t$ in the $j$-th
subspace of dimension $k-l$ and $dt$ stands for the Lebesgue measure
  on that space (we used Fubini theorem: see \cite{C} p. 334). Remark
  that $t$ lives in a ball $B(0, \epsilon_0)$. \\

\noindent So we have a upper bound of $\nu_n(B_n(x))$ by:
$$\frac{C(X)}{\lambda_l(f^n)} \sum_{j=1}^{i} \int \int_{f^n(B_n(x) \cap \tau_x(Y_j(t)))} \omega^{l} dt.$$
To control this integral, we have to bound from above the $2l$-dimensional volume of 
$f^n(B_n(x) \cap \tau_x(Y_j(t)))$ for some good points $x$ of $\nu_n$. In order to do that, 
we explain first what are these good points for $\nu_n$ then
we will bound the volume using Yomdin's approach and finally we will finish the bound of $\nu_n(B_n(x))$.

\subsubsection{{\bf Good points for the measure $\nu_n$ \\}}

\noindent In what follows, we consider a constant $L>0$ 
and an integer $n_0$ such that: 
$$\int \log d(x,I)  d \mu_{n} (x) \geq -L,$$
for $n \geq n_0$. The existence of these constants follows 
easily from Hypothesis (H).\\

Let $\delta > 0$. Our goal is to show that the entropy of  $\mu$ is greater than $\log d_l - \delta$.
We choose  a constant $C_0$ large enough ($1/C_0 \ll \delta$).

We are going to show that Hypothesis (H) implies that the orbits of  generic points
of the measure $\nu_{n}= \frac{(f^n)^* \omega^{l}
\wedge \omega^{k-l}}{\lambda_l(f^n)}$ are not close to the indeterminacy set $I$. They are
going to be the good points. 

\begin{lemme}{\label{visite}}
For $n \geq n_0$, there exists a set $A_n$ of $\nu_n$-measure greater or equal to $1 - C_0^{-1}$
whose points $x \in A_n$ satisfy:
$$  \prod_{i \in [0,n-1]} d(f^i(x),I)  \geq  e^{- C_0 L  n}.$$
\end{lemme}
\emph{Proof.}
We have
$$ \frac{1}{n} \int \log \left( \prod_{i \in [0,n-1] } d(f^i(x),I) \right) d \nu_n(x) =\frac{1}{n}  \int \sum_{i=0}^{n-1}
  \log d(f^i(x),I)d \nu_n(x).$$
Since $\mu_n = \frac{1}{n} \sum_{i=0}^{n-1} (f^i)_{*} \nu_n$:
 $$ \frac{1}{n} \int \log \left( \prod_{i \in [0,n-1] } d(f^i(x),I)
\right) d \nu_n(x) = \int \log d(x,I) d \mu_n(x).$$
Thanks to our hypothesis, this last integral is $ \geq -  L$.

Now, if we denote $h(x)=\frac{1}{n} \log \left( \prod_{i \in [0,n-1]} d(f^i(x),I) \right)$ and $A_n:= \{x \mbox{, } h(x) \geq -
   C_0 L  \}$, we have:
$$ \int_{A_n} h(x) d \nu_n(x) + \int_{X \setminus A_n} h(x) d \nu_n(x)
\geq -  L.$$
But $\int_{A_n} h(x) d \nu_n(x) \leq 0$ and $\int_{X \setminus A_n}
h(x) d \nu_n(x) \leq - C_0 L  \nu_n(X \setminus A_n)$.

This implies that $\nu_n(X \setminus A_n) \leq 1/C_0$. 

The set $A_n$ is indeed of measure  $\geq 1 -
 C_0^{-1}$ and if $x \in A_n$ then:
$$  \prod_{i \in [0,n-1]} d(f^i(x),I)  \geq  e^{- C_0 L n}, $$ 
which is what we wanted. \hfill $\Box$ \hfill \\ 

The orbit of points in  $A_n$ are not too close to $I$. These are the \emph{good points} for the measure $\nu_n$.

We now prove the upper bound of the volume.

\subsubsection{Upper bound for the volume of $f^n(B_n(x)
  \cap \tau_x(Y_j(t)))$ for $x\in A_n$}
Let $Y_0$ denote one of the $\tau_x (Y_j(t))$ (where $Y_j(t)$ is the fiber of $\pi_j$ with $t$ in the $j$-th
subspace of dimension $k-l$). Our aim is to prove:
\begin{proposition}{\label{recouvrement}}
The $2l$-dimensional volume of $f^{n}(Y_0 \cap B_n(x))$
is less or equal to:
$$ C(X,l,r)^{n/m+2m} \times K^{
  \frac{2npl}{r} + \frac{4mpl}{r}} \times \prod_{0 \leq i \leq n-1}
  d(f^i(x),I)^{\frac{-4pl}{r}}.$$ 
Here $C(X,l,r)$ is a constant that depends only on $X$, of the complex dimension 
$l$ of $Y_0$ and the regularity $r$ that we chose. The constants $K=K(f)$ and $p=p(f)$ are those of paragraph \ref{merom}.
\end{proposition}
Observe that the upper bound does not depend on the fiber $Y_j(t)$ that we consider.

Before proving the proposition, we give the upper bound of the $2l$-dimensional
 volume of $f^n(B_n(x) \cap \tau_x(Y_j(t)))$ that follows from the proposition.

Recall that we fixed $\delta$ and $C_0$. Now, let $r$ be such that $\frac{1}{r} \log K < \delta$ and $\frac{C_0 L}{r} <
\delta$. Then,  we choose $m$ so that
$\frac{1}{m} \log(C(X,l,r)) < \delta$ where $C(X,l,r)$ is the constant from the previous
proposition. Reformulating the previous proposition we have that the $2l$-dimensional
 volume of $f^{n}(Y_0 \cap B_n(x))$ is bounded by:
$$ C(X,l,r,m,p,K) e^{\delta n} \times e^{2\delta npl} \times
  \prod_{0 \leq i \leq n-1}
  d(f^i(x),I)^{\frac{-4pl}{r}} .$$ 
Finally, if $x$ is in $A_{n}$ (i.e. if $x$ is a good point for the measure $\nu_n$), the $2l$-dimensional
 volume of $f^{n}(Y_0 \cap B_n(x))$ is bounded from above by (see Lemma \ref{visite}):
$$e^{ 4 \delta npl} e^{\frac{4pl C_0 L n}{r}} \leq e^{ 8 \delta
  npl},$$
if $n$ is large (independently of $x \in A_n$).

It is this upper bound that we use now to finish the upper bound of
 $\nu_n(B_n(x))$ for $x \in A_n$. \\

\noindent{\bf End of the proof of the upper bound of $\nu_n(B_n(x))$ for $x \in A_n$}

Recall that we have bounded $\nu_n(B_n(x))$ by:
$$\frac{C(X)}{\lambda_l(f^n)} \sum_{j=1}^{i} \int \int_{f^n(B_n(x) \cap \tau_x(Y_j(t)))} \omega^{l} dt.$$
Now, if $x\in A_{n}$, we get:
$$\nu_n(B_n(x)) \leq \frac{e^{10 \delta npl}}{\lambda_l(f^n)},$$
for $n$ large enough. The ``$n$ large enough'' does not depend on $x \in A_n$. This quantity is approximately $d_l^{-n}$ and it stands for $x \in A_{n}$ which is a set of measure
 $\geq 1 - \frac{1}{C_0}$ for $\nu_{n}$. This is the upper bound that we wanted and
it will allow us to bound the entropy of $\mu$.
 
It remains to prove Proposition \ref{recouvrement}, which is the purpose of 
rest of this section. \\

\noindent{\bf Proof of Proposition \ref{recouvrement} \\}

 Consider $g=f^a$ an iterate of $f$ and let $x \in X$. We define $g_x= \tau_{g(x)}^{-1} \circ g \circ \tau_x$. We also define $g_{x, s(x)}= h(0,\frac{1}{s(x)})  \circ g_x \circ h(0,s(x))$  where $h(0,t)$ is the homothety of center $0$ and ratio $t$ in $\Cc^k$. Here, $s(x)$ is defined by:
$$s(x)=s_a(x)= \left( \frac{ d(x,I) \times \dots \times
  d(f^{a-1}(x),I)}{K^a} \right)^{p}.$$
We will consider later the particular cases $a=1$
(i.e. $s(x)=\eta(x)$) and $a=m$ (i.e. $s(x)=\rho(x)$).

In what follows, we are going to consider $C^r$-sizes associated to $2l$
 (i.e. sets that will be cover by some $h([0,1]^{2l})$ with $h\in C^r$).

First, we prove the following lemma by induction:
\begin{lemme}
Let $Z_0$ be a set of complex dimension $l$ such that the 
$C^r$-size of $\tau_x^{-1}(Z_0
\cap B(x,s(x)))$
is $\leq s(x)$.
 
Then, for $j \geq 1$, we can cover $g^{j-1}(Z_0 \cap B(x, s , j , g))$
by a union of $N_j$ sets $Z$ for which the $C^r$-size of
$\tau_{g^{j-1}(x)}^{-1}(Z)$ is $\leq s(g^{j-1}(x))$ and $N_j$ is bounded from above by
:
$$C(X,l,r)^{j-1} \prod_{0 \leq i \leq j-1} s(g^i(x))^{-2l/r} .$$
Here $B(x,s , j , g)$ is the dynamical ball:
$$B(x, s , j , g)= \{y \mbox{, } d(g^i(x),g^i(y)) \leq
s(g^i(x)) \mbox{ pour } i \in [0,j-1] \}.$$
\end{lemme}
\emph{Proof.} For $j=1$, the lemma stands by hypothesis.

Assume now that the induction assumption stands for $j-1$.

Observe that:
$$g^{j-1} (Z_0 \cap B(x, s, j ,g))= g( g^{j-2}(Z_0 \cap B(x, s ,j-1,g)))
\cap B(g^{j-1}(x), s(g^{j-1}(x))).$$
Let $Z$ be one of the $N_{j-1}$ sets whose union covers $g^{j-2}(Z_0 \cap
B(x, s ,j-1,g))$. The $C^r$-size of $\tau_{g^{j-2}(x)}^{-1}(Z)$ is 
$\leq s(g^{j-2}(x))$ by the induction assumption. To prove the lemma,
we bound from above the numbers of sets $Y$ which cover 
$g(Z) \cap B(g^{j-1}(x), s(g^{j-1}(x)))$  for which the $C^r$-size of
  $\tau_{g^{j-1}(x)}^{-1}(Y)$ is $\leq s(g^{j-1}(x))$.

We consider $\widetilde{Z}= h(0, 1/ s(g^{j-2}(x))) \circ
  \tau_{g^{j-2}(x)}^{-1} (Z)$. The $C^r$-size of $\widetilde{Z}$ is $ \leq s(g^{j-2}(x)) \times
  \frac{1}{s(g^{j-2}(x))}=1$. Furthermore, since $Z$ is in the ball
  $B(g^{j-2}(x), s(g^{j-2}(x)))$ (else we only consider the part of $Z$ 
that is in the ball and we still denote it $Z$), $\widetilde{Z}$ is in the ball
  $B(0, C(X))$ (where $C(X)$ is a constant that depends only on $X$). Using Proposition \ref{Gromov} of Section \ref{Yomdin} with $f=g_{g^{j-2}(x),
  s(g^{j-2}(x))}$ and $Y_0= \widetilde{Z}$ we get that we can cover $g_{g^{j-2}(x),
  s(g^{j-2}(x))} ( \widetilde{Z}) \cap B ( 0, \beta)$ (we take
  $\beta= C(X) \frac{s(g^{j-1}(x))}{s(g^{j-2}(x))} $) by
 $$C(X,l,r) \left(1+ \frac{\| D_r g_{g^{j-2}(x), s(g^{j-2}(x))} \| }{
  \beta} \right)^{2l/r}$$
sets $\widetilde{Y}$ of $C^r$-size $\leq C(X) \frac{s(g^{j-1}(x))}{s(g^{j-2}(x))}$.
Here the norm $\|.\|$ is taken over the ball $B ( 0, C(X) +
  \sqrt{2l})$. The images $Y$ of the $\widetilde{Y}$ by $\tau_{g^{j-1}(x)} \circ
  h(0, s(g^{j-2}(x)))$ cover 
\begin{equation*}
\begin{split}
&\tau_{g^{j-1}(x)} \circ  h(0,  s(g^{j-2}(x))) (g_{g^{j-2}(x),
  s(g^{j-2}(x))} ( \widetilde{Z}) \cap B(0 , \beta))\\
& =  g(Z)  \cap  \tau_{g^{j-1}(x)} \circ  h(0,  s(g^{j-2}(x)))(B(0,
  \beta))
\end{split}
\end{equation*}
which contains
$$ g(Z) \cap B(g^{j-1}(x), s(g^{j-1}(x))).$$
This is the set we wanted to cover and
  $\tau_{g^{j-1}(x)}^{-1}(Y)= h(0, s(g^{j-2}(x)))( \widetilde{Y})$
 is of $C^r$-size $\leq s(g^{j-1}(x))$ up to dividing it into $C(X)^{2l}$ pieces as in Section
  \ref{Yomdin} (this multiplies $N_j$ by a universal constant). 

To finish the proof, we have to count the number of pieces
$Y$ that we constructed for which the $C^r$-size of $\tau_{g^{j-1}(x)}^{-1}(Y)$ is bounded form above by
$s(g^{j-1}(x))$. Indeed, the union of those sets covers $g^{j-1} (Z_0 \cap B(x, s, j
,g))$. 

To control $N_j$, we need a control of the norm  $\| D_r g_{g^{j-2}(x),
  s(g^{j-2}(x))}  \|$ on the ball
$B(0, C(X) + \sqrt{2l})$. 

We admit temporarily that this norm is $ \leq C(X,l,r) s(g^{j-2}(x))^{-1}$.

Then:
$$N_j \leq N_{j-1} C(X,l,r) \left(1 + \frac{\| D_r g_{g^{j-2}(x), s(g^{j-2}(x))}
  \| s(g^{j-2}(x)) }{ s(g^{j-1}(x))} \right)^{2l/r} ,$$ 
which is smaller than:
$$N_{j-1}  C(X,l,r) \left( \frac{2 C(X,l,r) }{s(g^{j-1}(x))} \right)^{2l/r}
\leq N_{j-1} C(X,l,r)   s(g^{j-1}(x))^{-2l/r} $$    
up to changing $C(X,l,r)$. This concludes the proof of the lemma up 
to the upper bound of the norm of $\| D_r g_{g^{j-2}(x),
  s(g^{j-2}(x))}  \|$ on the ball $B ( 0, C(X) + \sqrt{2l})$. \\

\noindent {\bf Upper bound of the norm $\| D_r g_{g^{j-2}(x),
    s(g^{j-2}(x))}  \|$ on $B ( 0, C(X) + \sqrt{2l})$ }

\noindent Since 
$$g_{g^{j-2}(x),  s(g^{j-2}(x))}= h(0,\frac{1}{s(g^{j-2}(x))})  \circ g_{g^{j-2}(x)} \circ h(0,
s(g^{j-2}(x))),$$
 $\| \partial_r g_{g^{j-2}(x), s(g^{j-2}(x))}  \|$ is equal to $ s(g^{j-2}(x))^{r-1} \| \partial_r g_{g^{j-2}(x)} \|$ where that last norm is taken over the ball
 $$B(0,s(g^{j-2}(x))
  (C(X) + \sqrt{2l}))$$
(see Section \ref{Yomdin} for notations).

To prove the upper bound of the norm, we are going to prove that:
 $$g_{g^{j-2}(x)}(B(0,2 s(g^{j-2}(x))   (C(X) + \sqrt{2l})))$$
 is contained in the ball $B(0 , 1)$. We will then deduce the upper bound of
$\| \partial_r g_{g^{j-2}(x)} \|$  on $B(0,s(g^{j-2}(x))
  (C(X) + \sqrt{2l}))$ by $C(X,r) (s(g^{j-2}(x)) (C(X) +
  \sqrt{2l}))^{-r}$ thanks to Cauchy inequalities. This gives exactly 
the upper bound that we want.

So, we show that:
$$g_{g^{j-2}(x)}(B(0,2 s(g^{j-2}(x))   (C(X) + \sqrt{2l})))$$
 is contained in $B(0 ,1)$.

If we let $y=g^{j-2}(x)$, we have: 
$$g_{g^{j-2}(x)}(B(0,2 s(g^{j-2}(x))   (C(X) + \sqrt{2l})))=
g_y(B(0, 2 s(y) (C(X) + \sqrt{2l} )))$$
which is equal to:
$$ \tau_{f^a(y)}^{-1} \circ f^a \circ \tau_y (B(0, 2 s(y) (C(X) + \sqrt{2l} )))$$
because $g=f^a$. Furthermore:
$$\tau_{f^a(y)}^{-1} \circ f^a \circ \tau_y= f_{f^{a-1}(y)} \circ
\dots \circ f_y,$$
with $f_x := \tau_{f(x)}^{-1} \circ f \circ \tau_x$.

Now we use Lemma \ref{estimee} of Section
\ref{merom} to control the differential of
$f_y$ on $B(0, 2 s(y) (C(X) + \sqrt{2l} ))$.

If $z$ is a point of the ball $B(0, 2 s(y) (C(X) + \sqrt{2l} ))$
then the distance between $\tau_y(z)$ and $I$ is $\geq d(y, I) -
2 s(y) C(X) (C(X) + \sqrt{2l} )$. But that last quantity is $\geq \frac{d(y,I)}{2}$ since by definition of $s(y)$, we have
 $s(y) \leq \frac{d(y,I)}{K}$ and we can assume that $K$ is large compared to the constants
that depend only on $X$ and $l$ (recall that
 $l$ is the complex dimension  of $Z_0$: it is between  $0$ and $k$, 
so in particular they are only a finite number of such quantities).
Using Lemma \ref{estimee}, we get an upper bound of $\| D f_y
\|$ on the ball $B(0, 2 s(y) (C(X) + \sqrt{2l} ))$ by $K C(X) 2^p
d(y,I)^{-p}$. Using the control over the differential, this implies
that the image of  $B(0, 2 s(y) (C(X) + \sqrt{2l} ))$ by $f_y$
is contained in $B(0, K C(X) 2^p d(y,I)^{-p} 2 s(y) (C(X) +
\sqrt{2l}))$. But since:
$$s(y)=  \left( \frac{ d(y,I) \times \dots \times
  d(f^{a-1}(y),I)}{K^a} \right)^{p},$$
we have:
$$ K C(X) 2^p d(y,I)^{-p} 2 s(y) (C(X) + \sqrt{2l}) \leq \left( \frac{
  d(f(y),I) \times \dots \times d(f^{a-1}(y),I)}{K^{a-1}} \right)^{p},$$
since we can assume that $K$ is large compared to the $C(X)$.

So we have proved that the image of $B(0, 2 s(y) (C(X) + \sqrt{2l} ))$ by $f_y$
is contained in 
$$B \left( 0, \left( \frac{ d(f(y),I) \times \dots \times
  d(f^{a-1}(y),I)}{K^{a-1}} \right)^{p} \right).$$
Now, if we do again what we just did for  $f(y)$ instead of
$y$, we get that the image by $f_{f(y)} \circ f_y$ 
of the ball $B(0, 2 s(y) (C(X) + \sqrt{2l} ))$ is contained in the ball:
$$B \left( 0, \left( \frac{ d(f^2(y),I) \times \dots \times
  d(f^{a-1}(y),I)}{K^{a-2}} \right)^{p} \right),$$
and so on. At the end, we have that the image of the ball $B(0, 2
s(y) (C(X) +\sqrt{2l} ))$ by $f_{f^{a-1}(y)} \circ
\dots \circ f_y=\tau_{f^a(y)}^{-1} \circ f^a \circ \tau_y$ is contained in the ball:
$$B \left( 0,  K C(X) 2^p d(f^{a-1}(y),I)^{-p}  \left( \frac{
  d(f^{a-1}(y),I)}{K} \right)^{p} \right),$$
which is contained in $B(0, 1)$ for $K$ large enough.

This concludes the proof of the upper bound of the norm $\| D_r g_{g^{j-2}(x),
  s(g^{j-2}(x))}  \|$ on the ball $B ( 0, C(X) + \sqrt{2l})$ and that concludes the proof of the lemma.
 \hfill $\Box$ \hfill \\ 

Now we will use that lemma to prove Proposition \ref{recouvrement}. 
Recall some notations first. The set
$Y_0$ is one the fiber 
$\tau_x (Y_j(t))$, $n=m \phi(n) + r(n)$ with
$0 \leq r(n) <m$,
$$\rho(x)= \left( \frac{ d(x,I) \times \dots \times
  d(f^{m-1}(x),I)}{K^m} \right)^{p}$$
and

$$\eta(x)= \left( \frac{ d(x,I)}{K} \right)^{p}.$$
Finally, we denote:
$$B_n(x):= B(x, \rho ,  \phi(n), f^m) \cap f^{-\phi(n)m +
  m}(B(f^{\phi(n)m -m}(x) , \eta, r(n)+m,f)).$$
Applying the previous lemma for $g=f^m$ (and thus $s(x)=\rho(x)$),
$j=\phi(n)$ and $Z_0= Y_0 \cap B(x, \rho(x))$ (whose image by $\tau_{x}^{-1}$ is of $C^r$-size$\leq \rho(x)$ up to dividing into $C(X)^{2l}$
pieces because $Y_j(t)$ is a linear subspace), we get that we can cover
$f^{m(\phi(n)-1)}(Y_0 \cap B(x, \rho ,  \phi(n) , g))$ by a number $N_{\phi(n)}$  of sets  $Z$ 
for which the $C^r$-size of
$\tau_{g^{\phi(n)-1}(x)}^{-1}(Z)=\tau_{f^{m(\phi(n)-1)}(x)}^{-1}(Z)$ is 
$\leq \rho(g^{\phi(n)-1}(x))$ and $N_{\phi(n)}$ bounded from above by:
$$C(X,l,r)^{\phi(n)} \prod_{0 \leq i \leq \phi(n)-1} \rho(g^i(x))^{-2l/r}.$$
So we went up to $f^{m(\phi(n)-1)}(x)$ and we still have to go to $f^n(x)$. 

For that, we use the above lemma again with for $Z_0$ one of the $N_{\phi(n)}$ pieces $Z$,
$g=f$ (so now $s(x) = \eta(x)$), $j= n - m(\phi(n)-1)  =
r(n) + m$ and $x=f^{m(\phi(n)-1)}(x)$. We can do that because the $C^r$-size of
$\tau_{f^{m(\phi(n)-1)}(x)}^{-1}(Z_0)$ is $ \leq 
\rho(f^{m(\phi(n)-1)}(x)) \leq \eta(f^{m(\phi(n)-1)}(x))$. 
So we get that we can cover $f^{r(n)+m-1}(Z \cap
B(f^{m(\phi(n)-1)}(x),\eta,r(n)+m,f))$ by a union of $M$ sets $Y$ for which the $C^r$-size
of $\tau_{f^{n-1}(x)}^{-1}(Y)$ is $\leq \eta(f^{n-1}(x))$ and $M$ is less than:
$$C(X,l,r)^{m+r(n)-1} \prod_{1 \leq i \leq m+r(n)}
\eta(f^{n-i}(x))^{-2l/r}.$$

The sets $Y$ that we constructed belong to (up to keeping only the part that lies in it):
$$B(f^{r(n)+m-1+m(\phi(n)-1)}(x), \eta(f^{r(n)+m-1+m(\phi(n)-1)}(x)))=B(
f^{n-1}(x), \eta(f^{n-1}(x))).$$
The $C^1$-size of these $Y$ is smaller than $C(X) \eta(f^{n-1}(x))$
which implies that the diameter of $h([0,1]^{2l})$ (where $h$ is the map in $C^r$ associated to $Y$) is 
smaller than $C(X,l)
\eta(f^{n-1}(x))$. So, the set $h([0,1]^{2l})$ is contained in
$$B \left( f^{n-1}(x),\frac{d(f^{n-1}(x),I)}{2} \right).$$
Since the differential of $f$ in this last ball is bounded by
$$K2^p d(f^{n-1}(x),I)^{-p}$$ 
using Lemma \ref{estimee},
one gets that the images by $f$ of those $Y$ are of $C^1$-size
bounded by $C(X) \eta(f^{n-1}(x)) K2^p d(f^{n-1}(x),I)^{-p}$. So their 
$2l$-dimensional volume is $\leq 1$.

Summing up, we have covered
$$f^{r(n)+m}(f^{m(\phi(n)-1)}(Y_0 \cap B(x,\rho,\phi(n),g)) \cap
B(f^{m(\phi(n)-1)}(x),\eta,r(n)+m,f))$$
which contains $f^n(B_n(x) \cap Y_0)$
by a number $N$ of sets $Y$ of volume $\leq 1$ with:
$$N \leq C(X,l,r)^{\phi(n)+2m} \prod_{0 \leq i \leq \phi(n)-1}
\rho(g^i(x))^{-2l/r} \prod_{1 \leq i \leq m+r(n)} \eta(f^{n-i}(x))^{-2l/r}.$$
Using now the fact that:
$$\rho(y)= \left(\frac{ d(y,I) \times \dots \times
  d(f^{m-1}(y),I)}{K^m} \right)^{p},$$
and
$$\eta(y)=\left(\frac{ d(y,I)}{K} \right)^{p},$$
we have:
$$\prod_{0 \leq i \leq \phi(n)-1} \rho(g^i(x))^{-2l/r} \leq
  K^{\frac{2m \phi(n)pl}{r}}   \prod_{0 \leq i \leq \phi(n)m-1}
  d(f^i(x),I)^{\frac{-2pl}{r}},$$
and
$$ \prod_{1 \leq i \leq m+r(n)} \eta(f^{n-i}(x))^{-2l/r} \leq  K^{\frac{4mpl}{r}}   \prod_{1 \leq i \leq m+r(n)}
  d(f^{n-i}(x),I)^{\frac{-2pl}{r}}.$$
Finally, we have covered $f^n(B_n(x) \cap Y_0)$ by a number $N$
of sets $Y$ of volume $\leq 1$ with:
$$N \leq C(X,l,r)^{n/m + 2m} K^{\frac{2npl}{r} + \frac{4mpl}{r}}
\prod_{0 \leq i \leq n-1} d(f^i(x),I)^{\frac{-4pl}{r}}.$$
That concludes the proof of Proposition  \ref{recouvrement}.

\subsection{Lower bound for the entropy of $\mu$}
Recall that we consider a cluster value $\mu$ of the sequence
$$\mu_n= \frac{1}{n} \sum_{i=0}^{n-1} f^i_{*}\left( \frac{(f^n)^* \omega^{l}
  \wedge \omega^{k-l}}{\lambda_l(f^n)} \right)$$
and that in order to simplify the notations we assume that $(\mu_n)$
converges to $\mu$. By assumption,  $\mu$ gives no mass to the indeterminacy set $I$
and it is invariant by Lemma \ref{invariante}. The aim of this section is to prove
that the metric entropy $h_{\mu}(f)$ is  $\geq \log d_l -  \delta$. 
This will implies Theorem \ref{ENTROPY} by letting $\delta \rightarrow
  0$.

So we have to bound $h_{\mu}(f)$. Here is the plan of this section: 
first we will construct partitions of finite entropy for $\mu$ that will be used 
latter with the proof of the variational principle to get the lower bound of the entropy 
that we want.

\subsubsection{{\bf Construction of the partitions}\\}
The proof is the same than the one of Ma{\~n}é (see Lemma $2$ in \cite{M}). 
We give it for the sake of the reader since we will use it in what follows.
We consider a function $s(x)$
comprised between $0$ and $1$. Later, we will take the values
 $\rho(x)$ or $\eta(x)$ for $s(x)$.
\begin{proposition}{\label{partition}}
We can construct a countable partition $\Pcal$ of $X \setminus
\{ s = 0 \}$ such that:
\begin{enumerate}
\item If  $x \in X \setminus \{ s=0 \}$, then $\mbox{diam} \Pcal(x) <
   s(y)$ for all $y \in \Pcal(x)$ (here $\Pcal(x)$ denotes the atom of 
the partition that contains $x$).

\item For any probability measure $\nu$ such that  $\int \log s(x) d \nu(x) > -
   \infty$, we have $H_{\nu}(\Pcal) < + \infty$. Here $H_{\nu}(\Pcal)$
denotes the entropy of the partition $\Pcal$ for the measure $\nu$.
\end{enumerate}
\end{proposition}
Before proving the proposition, recall the following Ma{\~n}é's lemma (see Lemma $1$ in \cite{M}):
\begin{lemme}{\label{Mane}}

If $\sum_{n=0}^{+ \infty} x_n$ is a series with $0 \leq x_n \leq 1$
for all $n$ and if $ \sum_{n=0}^{+ \infty} n x_n < + \infty$ then
$$\sum_{n=0}^{+ \infty} x_n \log(1/x_n) < + \infty$$
with the convention that $x_n \log(1/x_n)=0$ when $x_n=0$.
\end{lemme}
\emph{Proof of Lemma \ref{Mane}}
Here is the Ma{\~n}é's proof:

Let $S$ be the set of integers $n \geq 0$ for which $x_n \neq 0$
and $\log(1/x_n) \leq  n$. If $n \notin S$ then $x_n \leq e^{-n}$. Furthermore:
$$ \sum_{n=0}^{+ \infty} x_n \log(1/x_n) \leq \sum_{n \in S} n x_n +
\sum_{n \notin S} (\sqrt{x_n}) (\sqrt{x_n}) \log(1/x_n).$$
But since $(\sqrt{t}) \log(1/t) \leq 2 e^{-1}$ for all $t \geq 0$,
we have:
$$ \sum_{n=0}^{+ \infty} x_n \log(1/x_n) \leq  \sum_{n=0}^{+ \infty} n
x_n + 2e^{-1} \sum_{n \notin S} \sqrt{x_n}$$
which is less than:
$$ \sum_{n=0}^{+ \infty} n x_n  + 2e^{-1} \sum_{n=0}^{+ \infty}
e^{-n/2}.$$
And that gives the lemma. \hfill $\Box$ \hfill \\

\noindent \emph{Proof of Proposition \ref{partition}}
First of all, there are constants $C>0$ and $r_0>0$ such that for
 $0 < r \leq r_0$, there exists a partition $\Pcal_r$ of $X$ whose 
elements have a diameter $\leq r$ and such that the number of elements of the
partition $| \Pcal_r |$ is $\leq C(1/r)^{2k}$.

Now, we define $V_n:=\{x \mbox{, } e^{-(n+1)} < s(x) \leq
  e^{-n} \}$ for $n \geq 0$.

Since the function $s$ is less than $1$, we have that  $X
\setminus \{ s = 0 \} = \cup_{n \geq 0} V_n$.

Let $\Pcal$ be the partition defined as follows: for $n$ fixed,
we consider the sets $Q \cap V_n$ for $Q \in \Pcal_{r_n}$ with
$r_n = e^{-(n+1)}$. This defines a partition of $V_n$. Now,
we get the partition $\Pcal$ of $X \setminus \{ s=0 \}$ by taking all
the $n$ between $0$ and $+ \infty$.

If $x\notin\{ s = 0 \}$, then $x \in V_n$ for some $n
\geq 0$ and then the atome of $\Pcal$ containing $x$,
$\Pcal(x)$, is contained in an atom of $\Pcal_{r_n}$, we have:
$$\mbox{diam} \Pcal(x) \leq e^{-(n+1)} < s(y)$$
for all $y \in \Pcal(x) \subset V_n$. This proves the first point of
Proposition \ref{partition}. \\

We now consider a measure $\nu$ such that $\int \log s(x) d \nu(x) > -
   \infty$. We want to show that $H_{\nu}(\Pcal) < + \infty$.

We have:
$$H_{\nu}(\Pcal) = \sum_{n=0}^{+ \infty} \left( - \sum_{P \in \Pcal
  \mbox{, } P \subset V_n} \nu(P) \log \nu(P) \right).$$
Using the inequality:
$$- \sum_{i=1}^{m_0} x_i \log x_i \leq \left( \sum_{i=1}^{m_0} x_i \right)
\left( \log m_0 - \log \sum_{i=1}^{m_0} x_i \right)$$
which comes from the convexity of the function $\phi(x)=x \log(x)$ for
$x \geq 0$, we get:
$$H_{\nu}(\Pcal) \leq  \sum_{n=0}^{+ \infty} \nu(V_n)( \log |
\Pcal_{r_n} | - \log \nu(V_n)).$$
Since the number $|\Pcal_{r_n} |$ of elements of $\Pcal_{r_n}$ is less than
 $C e^{2k(n+1)}$, we have:
$$H_{\nu}(\Pcal) \leq \log C + 2k \sum_{n=0}^{+ \infty} (n+1)
\nu(V_n) + \sum_{n=0}^{+ \infty} \nu(V_n) \log \left(
\frac{1}{\nu(V_n)} \right).$$
By assumption:
$$\int \log s(x) d \nu(x)= \int_{\cup_{n \geq 0} V_n}
\log s(x) d \nu(x) > -  \infty.$$
This implies that:
$$ \sum_{n=0}^{+ \infty} n \nu(V_n) < + \infty,$$
And the proposition is then deduced from Lemma \ref{Mane}.  \hfill $\Box$ \hfill

\subsubsection{ Lower bound for the entropy of $\mu$}
 In what follows, we denote $\Pcal$ (resp. $\Qcal$) the partition previously constructed for
$s(x)=\rho(x)$ (resp. $s(x) = \eta(x)$). We consider the 
restriction of $\Pcal$ and $\Qcal$ to $\Omega= X \setminus \cup_{i \geq 0} f^{-i} (I)$ 
(that we still denote
$\Pcal$ and $\Qcal$). They are partitions of $\Omega$. The advantage of those 
partitions over $\Omega$ is that the $f^i$ are well-defined on them.
In particular, we can define for example the partition 
$f^{-i}(\Pcal)$: its atoms are the $f^{-i}(P):=\{ x \in \Omega
\mbox{ with } f^{i}(x) \in P \}$ where the $P$ are the atoms of
$\Pcal$. Since $f(\Omega) \subset \Omega$, we get a partition of
$\Omega$. The measures that we consider ($\nu_n$, $\mu_n$ or $\mu$) have a mass $1$ on 
 $\Omega$. The parts of $X$ that we drop are of mass 0 for them. 
We remark that with our convention, we have:
$f^{-a}(f^{-b}(P))=f^{-a-b}(P)=\{x \in \Omega \mbox{ with } f^{a+b}(x)
\in P \}$. 
Recall that we denote:
$$\nu_n= \frac{(f^{n})^*
  \omega^{l}  \wedge \omega^{k-l}}{\lambda_l(f^{n})} $$
and that $\nu_n(A_{n}) \geq 1 - \frac{1}{C_0}$
  (see Lemma \ref{visite}).

In what follows, we denote $\nu_n':= \frac{{\nu_n}_{|
    A_{n}}}{\nu_n(A_{n})}$ (i.e. $\nu_n'(B)= \frac{\nu_n(B
    \cap A_{n})}{\nu_n(A_{n})}$).

Define the joint partition $\Pcal_{-n}$ of the partitions $\Pcal$ and
$\Qcal$ by (recall that $n=\phi(n)m +r(n)$ with $0 \leq r(n) < m$):
$$\Pcal_{-n}:= \Pcal \vee f^{-1}( \Pcal) \vee \dots \vee
f^{-\phi(n)m+m}(\Pcal) \vee f^{- \phi(n)m +m-1} ( \Qcal) \vee \dots
\vee f^{-n+1}(\Qcal).$$
First, we have the lemma:
\begin{lemme}
If $n$ is large enough, then
$$\nu_n' ( \Pcal_{-n}(x)) \leq  \frac{e^{10  \delta n l p}}{ \lambda_l(f^n)} \frac{1}{1 - \frac{1}{C_0}}.$$
for every atom $\Pcal_{-n}(x)$ of $\Pcal_{-n}$. 
\end{lemme}
\emph{Proof.} We have shown in the previous paragraph that if $n$ is large enough then for every
 $x \in A_n$:
$$\nu_n(B_n(x)) \leq \frac{e^{10 \delta npl}}{\lambda_l(f^n)}.$$
Consider now a $n$ large enough so that the previous property
is satisfied. Here:
$$B_n(x):= B(x, \rho ,  \phi(n), f^m) \cap f^{-\phi(n)m +
  m}(B(f^{\phi(n)m -m}(x) , \eta, r(n)+m,f)).$$
If $\Pcal_{-n}(x)$ does not contain any points of $A_{n}$ then
 $\nu_n'( \Pcal_{-n}(x)) =0$ and the lemma is true.
So we can assume that there exists $y \in
\Pcal_{-n}(x) \cap A_{n}$.

By definition of the joint partition, we have $\Pcal_{-n}(x)$ which is equal to:
\begin{equation*}
\begin{split}
&\Pcal(x) \cap  \dots \cap
f^{-\phi(n)m+m}( \Pcal(f^{\phi(n)m-m}(x)))\\
&\cap f^{-\phi(n)m+m-1}
(\Qcal(f^{\phi(n)m-m+1}(x))) \cap \dots \cap f^{-n+1}(
\Qcal(f^{n-1}(x))).
\end{split}
\end{equation*}

In particular, $f^i(y) \in \Pcal(f^i(x))$ for $i=0 \dots \phi(n)m-m$
and then $f^i(y) \in \Qcal(f^i(x))$ for $i= \phi(n)m-m+1 \dots n-1$. By
Proposition \ref{partition}, the diameter of $
\Pcal(f^i(x))$ is $\leq \rho(f^i(y))$ for $i=0 \dots
\phi(n)m-m$ and the diameter of $\Qcal(f^i(x))$ is $\leq \eta(f^i(y))$
for $i= \phi(n)m-m+1 \dots n-1$ which means:
$$\Pcal_{-n}(x) \subset B_n(y).$$
The lemma follows then first from the estimation of the previous paragraph
since $y \in A_{n}$ and secondly from the fact that $\nu_n(A_{n})$ is
 $ \geq 1 - \frac{1}{C_0}$.  \hfill $\Box$ \hfill \\

Thanks to that estimation on $\nu_n' ( \Pcal_{-n}(x)) $, we can bound the
entropy of $\mu$ using a variation of the proof of the variational principle.
We refer the reader to \cite{W} p.188-190 for the proof of the principle
and to \cite{BS}, \cite{De1} or \cite{Gu} for its use to bound from below 
the entropies of measure in holomorphic or meromorphic dynamics.

Let $q$ be an integer $2m<q<n$ (with $m$ from the above paragraph). For $0 \leq j \leq q-1$,
we let $a(j)=\left[ \frac{n-j}{q} \right]$ and then
$$\{0, 1 , \dots , n-1 \}= \{ j +rq+i \mbox{, } 0 \leq r \leq a(j)-2
  \mbox{, } 0 \leq i \leq q-1 \} \cup S(j)$$

where $S(j)= \{0,1, \dots , j-1 , j +(a(j)-1)q,  j +(a(j)-1)q +1, \dots, n-1 \}$ 
is of cardinality less than $3q$ since $j+(a(j)-1)q \geq j + \left( \frac{n-j}{q}-2 \right) q  =
n-2q$. We took the indexes  $r$ up to $a(j)-2$ so that 
$S(j)$ contains $n-q \dots n-1$ and so in particular $\phi(n)m -m +1
\dots n-1$ (we take $q$ large with respect to $m$). We denote $S_1(j)$
the elements of $S(j)$ other than $\phi(n)m -m +1 \dots n-1$ and
$S_2(j)$ the elements $\phi(n)m -m +1 \dots n-1$. 

Now, we have (see for example Proposition 4.3.3 of \cite{K.H}):
$$H_{\nu_n'}(  \Pcal_{-n} ) \geq - \log (
\sup_{P \in \Pcal_{-n}} \nu_n'(P)) \geq  -10 \delta n l p + \log \lambda_l(f^n) + \log \left( 1 - \frac{1}{C_0} \right),$$
by the previous lemma.

On the other hand, by the proof of the variational principle for $0 \leq j \leq q-1$, we have:
$$ \Pcal_{-n} = \bigvee_{r=0}^{a(j)-2} \left( f^{-(rq+j)}
\bigvee_{i=0}^{q-1} f^{-i} \Pcal \right) \vee \bigvee_{t \in S_1(j)} f^{-t}
\Pcal \vee \bigvee_{t \in S_2(j)} f^{-t} \Qcal.$$
So, (again by Proposition 4.3.3 in \cite{K.H}):
$$H_{\nu_n'}( \Pcal_{-n}  ) \leq
\sum_{r=0}^{a(j)-2} H_{\nu_n'} (f^{-(rq+j)} \bigvee_{i=0}^{q-1}
f^{-i} \Pcal)  + \sum_{t \in S_1(j)} H_{\nu_n'}( f^{-t}
\Pcal) + \sum_{t \in S_2(j)} H_{\nu_n'}( f^{-t}
\Qcal)$$
which is equal to: 
$$\sum_{r=0}^{a(j)-2} H_{f^{rq+j}_{*} \nu_n'} (\bigvee_{i=0}^{q-1}
f^{-i} \Pcal) + \sum_{t \in S_1(j)} H_{\nu_n'}( f^{-t} \Pcal) + 
\sum_{t \in S_2(j)} H_{\nu_n'}( f^{-t} \Qcal)$$
Summing this relation for $j=0 \dots q-1$, we get:
\begin{equation*}
\begin{split}
& q \left( -10 \delta n l p + \log \lambda_l(f^n)+\log \left( 1 -
 \frac{1}{C_0} \right) \right) \\
& \leq
\sum_{j=0}^{q-1} \sum_{r=0}^{a(j)-2} H_{f^{rq+j}_{*} \nu_n'  }
(\bigvee_{i=0}^{q-1} f^{-i} \Pcal) + \sum_{j=0}^{q-1}   \left( \sum_{t \in S_1(j)} H_{\nu_n'
}( f^{-t} \Pcal) + \sum_{t \in S_2(j)} H_{\nu_n' }( f^{-t} \Qcal)
 \right).
\end{split}
\end{equation*}
The integers $j +rq$ for $0 \leq j \leq q-1$ and $0 \leq r \leq a(j)-2$
are all distinct and $ \leq n-2q$. So we have that (using the convexity
of the function $\Phi(x)=x \log(x)$ for $x >0$):
$$\mbox{  (I):  } \frac{q}{n-2q+1} \left( -10 \delta n l p + \log
\lambda_l(f^n)+ \log \left( 1 - \frac{1}{C_0} \right)    \right)$$
which is less than:
$$H_{\frac{1}{n-2q+1} \sum_{p=0}^{n-2q} f^{p}_{*} \nu_n' }
(\bigvee_{i=0}^{q-1} f^{-i} \Pcal) + \sum_{j=0}^{q-1} \left(  \sum_{t \in S_1(j)} \frac{H_{\nu_n'}(
f^{-t} \Pcal)}{n-2q+1}  +  \sum_{t \in S_2(j)} \frac{H_{\nu_n'}(
f^{-t} \Qcal)}{n-2q+1}   \right)      .$$ 
Here is the plan of the rest of the proof. In a first time, we deduce 
from that inequality a lower bound of $\frac{1}{q}
H_{\frac{1}{n-2q+1} \sum_{p=0}^{n-2q} f^{p}_{*} \nu_n}(\bigvee_{i=0}^{q-1} f^{-i} \Pcal)$. 
Then we will pass to the limit in that inequality. \\

\noindent {{\bf  1) Lower bound of $\frac{1}{q}
    H_{\frac{1}{n-2q+1} \sum_{p=0}^{n-2q} f^{p}_{*} \nu_n}(\bigvee_{i=0}^{q-1} f^{-i} \Pcal)$} \\}

By definition, $\nu_n' := \frac{{\nu_{n}}_{|
    A_{n}}}{\nu_n(A_{n})}$. In particular, $ \nu_n' \leq
    \frac{\nu_n}{1 - \frac{1}{C_0}}$ and
$$\frac{1}{n-2q+1} \sum_{p=0}^{n-2q} f^{p}_{*} \nu_n'  \leq
    \frac{1}{(1 - \frac{1}{C_0})(n-2q+1)}  \sum_{p=0}^{n-2q} f^{p}_{*}
    \nu_n.$$
In order to control  $\frac{1}{q} H_{\frac{1}{n-2q+1} \sum_{p=0}^{n-2q}
 f^{p}_{*} \nu_n}(\bigvee_{i=0}^{q-1} f^{-i} \Pcal)$ with the inequality (I), we are 
going to use the following lemma:
\begin{lemme}
Let $\nu$ and $\nu'$ be two  probabilities such that $\nu' \leq \beta
\nu$ for some $\beta >1 $. Then for any partition $\Qcal$, we have:
$$H_{\nu'}(\Qcal) \leq \beta( H_{\nu}(\Qcal)+ 1).$$
\end{lemme}
\emph{Proof.} The function $\Phi(x)= - x \log(x)$ is increasing on 
$[0 , e^{-1}]$ and decreasing on $[e^{-1},1]$.

So we have:
\begin{equation*}
\begin{split}
H_{\nu'}(\Qcal) &= \sum_{Q \in \Qcal} - \nu'(Q) \log \nu'(Q)\\
&=\sum_{Q \in \Qcal \mbox{, } \nu(Q) \leq \frac{e^{-1}}{\beta}} - \nu'(Q) \log
\nu'(Q) + \sum_{Q \in \Qcal \mbox{, } \nu(Q) > \frac{e^{-1}}{\beta}} - \nu'(Q) \log
\nu'(Q)
\end{split}
\end{equation*}
which is less than:
$$\sum_{Q \in \Qcal \mbox{, } \nu(Q) \leq \frac{e^{-1}}{\beta}} - \beta \nu(Q) \log
(\beta \nu(Q))  + \sum_{Q \in \Qcal \mbox{, } \nu(Q) > \frac{e^{-1}}{\beta}} - \nu'(Q) \log
\nu'(Q).$$
Since they are at most $  \frac{\beta}{e^{-1}}$ of $Q \in \Qcal$ with
$\nu(Q) > \frac{e^{-1}}{\beta}$ and because on the interval $[0,1]$, the function $\Phi(x)$ is 
non negative and bounded by $e^{-1}$, we have:
\begin{align*}
& \quad & \quad  & \quad & \quad  &H_{\nu'}(\Qcal) \leq \beta H_{\nu}(\Qcal)  +  \beta.    & \quad & \quad &\quad &\quad &\quad &\quad &\quad &\Box  
\end{align*}

Now, since
$$ \frac{1}{n-2q+1} \sum_{p=0}^{n-2q} f^{p}_{*} \nu_n'  \leq
    \beta \frac{1}{n-2q+1}  \sum_{p=0}^{n-2q} f^{p}_{*}
    \nu_n,$$
with $\beta= \frac{1}{1- \frac{1}{C_0}}$ and
$$ \nu_n'  \leq  \beta \nu_n $$
for that same $\beta>1$, we have:
$$\frac{q}{n-2q+1} \left(-10 \delta n l p + \log
 \lambda_l(f^n)+\log \left( 1 - \frac{1}{C_0} \right) \right)$$
which is less than:
\begin{equation*}
\begin{split}
&\frac{1}{1 - \frac{1}{C_0}} [  H_{\frac{1}{n-2q+1} \sum_{p=0}^{n-2q} f^{p}_{*} \nu_n }
(\bigvee_{i=0}^{q-1} f^{-i} \Pcal) + 1\\
&+ \sum_{j=0}^{q-1} \left( \sum_{t \in S_1(j)} \frac{H_{\nu_n}(
f^{-t} \Pcal)}{n-2q+1} + \sum_{t \in S_2(j)} \frac{H_{\nu_n}(
f^{-t} \Qcal)}{n-2q+1} \right) + \frac{3q^2}{n-2q+1}  ]
\end{split}
\end{equation*}

(since the cardinality of $S(j)$ is $\leq 3q$).

This implies a lower bound of $\frac{1}{q}  H_{\frac{1}{n-2q+1} \sum_{p=0}^{n-2q} f^{p}_{*} \nu_n }
(\bigvee_{i=0}^{q-1} f^{-i} \Pcal)$ by
\begin{equation*}
\begin{split}
&\left( 1 - \frac{1}{C_0} \right) \left( \frac{1}{n-2q+1} \left(-10 \delta n l p + \log
 \lambda_l(f^n)+\log \left( 1 - \frac{1}{C_0} \right)\right) \right)\\
& - \frac{1}{q} - \frac{1}{q} \left(  \sum_{j=0}^{q-1}  \sum_{t \in S_1(j)} \frac{H_{\nu_n}(
f^{-t} \Pcal)}{n-2q+1} + \sum_{j=0}^{q-1}  \sum_{t \in S_2(j)} \frac{H_{\nu_n}(
f^{-t} \Qcal)}{n-2q+1} \right) - \frac{3q}{n-2q+1}.
\end{split}
\end{equation*}

It remains now to take the limit of that inequality when $n$ goes to $\infty$. \\

{{\bf  2) Pass to the limit $n \rightarrow + \infty$}\\}

First:
$$\frac{1}{n-2q+1} \left( -10 \delta n l p + \log
 \lambda_l(f^n)+\log \left( 1 - \frac{1}{C_0} \right) \right)$$
goes to  $- 10 \delta  lp +  \log d_l$ when $n\to \infty$.
We need the following proposition.
\begin{proposition}\label{passtothelimit}
We have:
\begin{enumerate}
\item $\int \log \rho d \mu > - \infty$.
\item For all $q > 2m$,
$$H_{\frac{1}{n-2q+1} \sum_{p=0}^{n-2q} f^{p}_{*} \nu_n  }
(\bigvee_{i=0}^{q-1} f^{-i} \Pcal)$$
converges to $H_{\mu}(\bigvee_{i=0}^{q-1} f^{-i} \Pcal)$ when  $n\to \infty$.
\item For $q > 1$:
$$\frac{1}{q} \sum_{j=0}^{q-1}  \sum_{t \in S_1(j)} \frac{H_{\nu_n }(
f^{-t} \Pcal)}{n-2q+1}$$
converges to $0$ when $n \to \infty$.
\item For $q > 1$:
$$\frac{1}{q} \sum_{j=0}^{q-1}  \sum_{t \in S_2(j)} \frac{H_{\nu_n }(
f^{-t} \Qcal)}{n-2q+1}$$
converges to $0$ when $n \to \infty $.
\end{enumerate}
\end{proposition}
We assume temporarily that the proposition is true. We finish the lower bound of the entropy
of $\mu$.

If we pass to the limit in the inequality of the previous paragraph, we get: 
$$\frac{1}{q} H_{\mu}(\bigvee_{i=0}^{q-1} f^{-i} \Pcal) \geq  \left( 1 - \frac{1}{C_0}
\right) (- 10 \delta lp + \log d_l ) - \frac{1}{q}.$$
If we let $q$ go to $\infty$, we have:
$$ h_{\mu}(f) \geq \left( 1 - \frac{1}{C_0} \right) (- 10 \delta  lp +
 \log d_l ) $$ 
because the entropy of $\Pcal$ is finite for $\mu$ from Proposition
\ref{partition} and the first point of the above Proposition.

This proves the theorem by letting $C_0$ go to $\infty$ then by 
letting $\delta$ go to $0$.

Up to the proof of Proposition \ref{passtothelimit}, we have 
proved Theorem \ref{ENTROPY}. \\

In order to simplify the notations, we denote:
$$\mu_n'= \frac{1}{n-2q+1} \sum_{p=0}^{n-2q} f^{p}_{*} \nu_n.$$
For the proof of the four points of the proposition, we will use the following lemma:
\begin{lemme}\label{estimation_integrale}
For $i=0 \dots q-1$, we have
$$ 0 \geq \int_{ \{ \rho \leq \epsilon \} } \log \rho \mbox{ } d((f^i)_{*} \mu_n')
\geq - \delta(\epsilon)$$
if $n$ is large enough. Here $\delta( \epsilon)$ goes to $0$
when $\epsilon$ goes to $0$.
\end{lemme}
\emph{Proof.} \\

{\bf First step.}
In a first time, we are going to bound from above  $(f^i)_{*} \mu_n'( \{x \mbox{, }  \rho(x)
\leq \epsilon \} )$ by $\delta'( \epsilon)$ for $n$ large (with $\delta'( \epsilon)$ going to  $0$
when $\epsilon$ goes to $0$).

Recall that:
$$ \rho(x)= \left( \frac{d(x,I) \times \dots \times
  d(f^{m-1}(x),I)}{K^m} \right)^p.$$
In particular, we have:
$$ \{x \mbox{, }  \rho(x) \leq \epsilon \} \subset \{x \mbox{, } d(x,I) \leq \epsilon^{
  \frac{1}{mp}} K \} \cup \dots \cup  \{ x \mbox{, } d(f^{m-1}(x),I) \leq \epsilon^{
  \frac{1}{mp}} K \}.$$
Now,
$$(f^i)_{*} \mu_n'( \{x \mbox{, }  \rho(x) \leq \epsilon \} ) \leq
\sum_{l=0}^{m-1}  (f^i)_{*} \mu_n'( \{x \mbox{, } d(f^l(x),I) \leq
\epsilon^{ \frac{1}{mp}} K  \} )$$
which is equal to:
$$\sum_{l=0}^{m-1}  (f^l)_{*} (f^i)_{*} \mu_n'( \{x \mbox{, } d(x,I) \leq
\epsilon^{ \frac{1}{mp}} K  \} ).$$
The measure  $\sum_{l=0}^{m-1}  (f^l)_{*} (f^i)_{*} \mu_n'$ is
lower than $\frac{mn}{n-2q+1} \frac{1}{n} \sum_{l=0}^{n-1}
(f^l)_{*} \nu_{n}=\frac{mn}{n-2q+1} \mu_n $ which converges to $m \mu$. 
Using Hypothesis (H), we know: 
$$\mu( \{x \mbox{, } d(x,I) \leq
\epsilon^{ \frac{1}{mp}} K  \}) \leq \delta'(\epsilon)$$
with $\delta'(\epsilon)$ converging to $0$ when $\epsilon$ goes to
$0$ since $\log d(x,I)$ is integrable for the measure $\mu$ and so $\mu$ puts no mass on $I$.

We have then:
$$(f^i)_{*} \mu_n'( \{x \mbox{, }  \rho(x) \leq \epsilon \} ) \leq
m(1 + \epsilon) ( \delta'( \epsilon) + \epsilon)= \delta'(\epsilon)$$
if $n$ is large and up to changing $\delta'(\epsilon)$ (of course it
depends on $m$). This gives the first step. \\

\noindent {\bf Second step}

We now prove the lower bound of  $\int_{ \{ \rho \leq \epsilon \}
} \log \rho \mbox{ } d((f^i)_{*} \mu_n')$.

By the definition of  $\rho$, we have:
\begin{equation*}
\begin{split}
&\int_{ \{ \rho \leq \epsilon \} } \log \rho \mbox{ } d((f^i)_{*}
\mu_n')\\
&= p \sum_{l=0}^{m-1} \int_{ \{ \rho \leq \epsilon \} } \log
d(f^l(x),I)  \mbox{ } d((f^i)_{*} \mu_n') -mp \log K  (f^i)_{*}
\mu_n'(\{ \rho \leq \epsilon \}) .
\end{split}
\end{equation*}
By the first step , we get:
$$\int_{ \{ \rho \leq \epsilon \} } \log \rho \mbox{ } d((f^i)_{*}
\mu_n') \geq p \sum_{l=0}^{m-1} \int_{ \{ \rho \leq \epsilon \} } \log
d(f^l(x),I)  \mbox{ } d((f^i)_{*} \mu_n') -mp \log K  \delta'(
\epsilon).$$
It remains to control $\sum_{l=0}^{m-1} \int_{ \{ \rho \leq \epsilon \} } \log
d(f^l(x),I)  \mbox{ } d((f^i)_{*} \mu_n')$.

For that, we split these integrals into two parts:
\begin{equation*}
\begin{split}
&\int_{ \{ \rho \leq
  \epsilon \} \cap \{x \mbox{, }  d(f^l(x),I) \leq \delta'( \epsilon) \}  } \log
d(f^l(x),I)  \mbox{ } d((f^i)_{*} \mu_n')\\
&+ \int_{ \{ \rho \leq \epsilon \} \cap \{x \mbox{, }  d(f^l(x),I) > \delta'( \epsilon) \} } \log
d(f^l(x),I)  \mbox{ } d((f^i)_{*} \mu_n').
\end{split}
\end{equation*}
The second part is greater than:
$$ \delta'(\epsilon) \log \delta'( \epsilon)$$
if $n$ is large enough using the first step. That quantity goes to $0$ when
$\epsilon$ goes to $0$.
 
For the first part, we have:

$$\sum_{l=0}^{m-1} \int_{ \{ \rho \leq
  \epsilon \} \cap \{x \mbox{, }  d(f^l(x),I) \leq \delta'( \epsilon) \}  } \log
d(f^l(x),I)  \mbox{ } d((f^i)_{*} \mu_n') $$
which is greater than:
$$ \sum_{l=0}^{m-1} \int_{\{x \mbox{, }  d(f^l(x),I) \leq \delta'( \epsilon) \}  } \log
d(f^l(x),I)  \mbox{ } d((f^i)_{*} \mu_n') $$
which is equal to:
$$\int_{\{x \mbox{, }  d(x,I) \leq \delta'( \epsilon) \}  } \log
d(x,I)  \mbox{ } d(  \sum_{l=0}^{m-1} (f^l)_{*}    (f^i)_{*} \mu_n'). $$
As in the first step, $\sum_{l=0}^{m-1} (f^l)_{*}    (f^i)_{*}
\mu_n'$ is less than the measure:
$$\frac{mn}{n-2q+1} \frac{1}{n} \sum_{l=0}^{n-1}(f^l)_{*} \nu_{n}=   \frac{mn}{n-2q+1} \mu_n $$
and the above integral is bounded from below by:
$$m(1 + \epsilon) \left( \int_{\{x \mbox{, }  d(x,I) \leq \delta'( \epsilon) \}  } \log
d(x,I)  \mbox{ } d \mu(x) - \epsilon \right)$$
when $n$ is large. Indeed, on one hand we have that:
$$\int \log d(x,I)
\mbox{ } d \mu_n(x)$$
converges to  $\int \log d(x,I)  \mbox{ } d
\mu(x)$ by Hypothesis (H). On the other hand: 
$$\int_{\{x \mbox{, }  d(x,I) > \delta'(
  \epsilon) \}  } \log d(x,I)  \mbox{ } d \mu_n (x)$$
converges to $\int_{\{x \mbox{, }  d(x,I) >
  \delta'( \epsilon) \}  } \log d(x,I)  \mbox{ } d \mu(x)$ up to 
choosing $\epsilon$ generic so that $\mu$ gives no mass to $\{x \mbox{, }  d(x,I)= \delta'( \epsilon) \} $.

Finally, since $\int_{\{x \mbox{, }  d(x,I) \leq \delta'( \epsilon) \}  } \log
d(x,I)  \mbox{ } d \mu(x)$ goes to $0$ when  $\epsilon$ converges
to $0$ by dominated convergence, the lemma follows. \hfill $\Box$ \hfill \\

\noindent \emph{End of the proof of Proposition \ref{passtothelimit}.\\}

{\bf First point of the proposition}

By Hypothesis (H), we have:
$$\int \log d(x, I) d \mu > - \infty$$
the integrability of $\log \rho$ follows from the invariance of the measure $\mu$.

$$ $$

{\bf Second point of the proposition}

We are going to prove by induction on $j=1 \dots q$ that
$H_{\mu_n'}(\bigvee_{i=0}^{j-1} f^{-i} \Pcal)$ converges to
$H_{\mu}(\bigvee_{i=0}^{j-1} f^{-i} \Pcal )$. The sequence
 $\mu_n'$
converges to $\mu$. The difficulty lies in the fact that $\Pcal$ is a countable
 partition and not a finite partition.

{\bf For $j=1$}

Here, we show that $H_{\mu_n'}(\Pcal)$ converges to $H_{\mu}(\Pcal)$.

We have:
$$H_{\mu_n'} (\Pcal)= \sum_{s=0}^{+ \infty} \sum_{P \in \Pcal \mbox{, } P \subset
  V_s} - \mu_n'(P) \log \mu_n'(P) $$
that we divide as:
$$ \sum_{s=0}^{s_0 - 1} \sum_{P \in \Pcal \mbox{, } P \subset
  V_s} - \mu_n'(P) \log \mu_n'(P) +  \sum_{s=s_0}^{+ \infty} \sum_{P \in \Pcal \mbox{, } P \subset
  V_s} - \mu_n'(P) \log \mu_n'(P).$$
Up to moving slightly the boundaries of the partition $\Pcal$, we can assume that 
$\mu$ gives no mass to the boundary of its elements. In particular, the first above term converges to:
$$ \sum_{s=0}^{s_0 - 1} \sum_{P \in \Pcal \mbox{, } P \subset
  V_s} - \mu(P) \log \mu(P)$$
when $n$ goes to infinity since we only consider a finite number of elements.

We now show that the second term is small if we take 
$s_0$ large then $n$ large.

We follow here the notations and the ideas of the proof of Proposition \ref{partition}.
$$\sum_{s=s_0}^{+ \infty} \sum_{P \in \Pcal \mbox{, } P \subset
  V_s} - \mu_n'(P) \log \mu_n'(P)$$
is less than:
$$\sum_{s=s_0}^{+ \infty} \mu_n'(V_s)( \log | \Pcal_{r_s} | - \log
\mu_n'(V_s))$$
which is in turn less than:
$$ \log C \sum_{s=s_0}^{+ \infty} \mu_n'(V_s) + 2k \sum_{s=s_0}^{+
  \infty} (s+1) \mu_n'(V_s) + \sum_{s=s_0}^{+ \infty} \mu_n'(V_s) \log \left(
  \frac{1}{\mu_n'(V_s)} \right).$$
But, first:
$$\sum_{s=s_0}^{+ \infty} \mu_n'(V_s) \leq  \mu_n'( \{ \rho \leq e^{-
  s_0} \} )$$
is as small as we want if we take $s_0$ large enough then $n$
large enough (this is exactly what we proved in the first step of the previous lemma with $i=0$).

Then:
$$\sum_{s=s_0}^{+  \infty} s \mu_n'(V_s) \leq - \int_{ \{ \rho \leq
  e^{- s_0} \} } \log \rho  d \mu_n'$$
is also as small as we want if we take $s_0$ large enough then 
 $n$ large enough thanks to the previous lemma with $i=0$.

Finally, following the proof of Lemma \ref{Mane}, we have:
$$\sum_{s=s_0}^{+ \infty} \mu_n'(V_s) \log \left(
  \frac{1}{\mu_n'(V_s)} \right) \leq \sum_{s=s_0}^{+ \infty} s
  \mu_n'(V_s) + 2 e^{-1} \sum_{s=s_0}^{+ \infty} e^{-s/2}$$
is also as small as we want if we take $s_0$ large enough then 
 $n$ large enough.

We have indeed shown that:
$$\sum_{s=s_0}^{+ \infty} \sum_{P \in \Pcal \mbox{, } P \subset
  V_s} - \mu_n'(P) \log \mu_n'(P)$$
is small and since all that we did remains true if we replace
$\mu_n'$ by $\mu$:
$$\sum_{s=s_0}^{+ \infty} \sum_{P \in \Pcal \mbox{, } P \subset
  V_s} - \mu(P) \log \mu(P)$$
is as small as we want for $s_0$ large enough. 

In particular, this implies that $H_{\mu_n'}(\Pcal)$ converges
to $H_{\mu}(\Pcal)$.

We continue the induction: we assume that $H_{\mu_n'}(\bigvee_{i=0}^{j-1} f^{-i} \Pcal )$ 
converges to
$$H_{\mu}(\bigvee_{i=0}^{j-1} f^{-i} \Pcal)$$
for some $j$ less than $q-1$ and we are going to show 
that the property holds for the rank $j+1$. \\

{\bf For $j+1$}

First, we have:
$$H_{\mu_n'}(\bigvee_{i=0}^{j} f^{-i} \Pcal)= H_{\mu_n'}( \bigvee_{i=0}^{j-1} f^{-i} \Pcal \vee f^{-j} (\Pcal)
)= H_{\mu_n'} ( \bigvee_{i=0}^{j-1} f^{-i} \Pcal) + H_{\mu_n'} (f^{-j} (\Pcal) |
\bigvee_{i=0}^{j-1} f^{-i} \Pcal)$$
by Proposition 4.3.3 in \cite{K.H}.

The first term converges to $H_{\mu} ( \bigvee_{i=0}^{j-1} f^{-i} \Pcal)$ 
by the induction assumption. We now show that the second term converges
to $ H_{\mu} (f^{-j} (\Pcal) |  \bigvee_{i=0}^{j-1} f^{-i} \Pcal)$. This will finish the 
induction and thus gives the second point of the proposition.

By definition, $H_{\mu_n'} (f^{-j} (\Pcal) |  \bigvee_{i=0}^{j-1} f^{-i} \Pcal)$ is equal to:
$$- \sum_{P_1 \in \bigvee_{i=0}^{j-1} f^{-i} \Pcal} \mu_n'(P_1) \sum_{s=0}^{+ \infty}
  \sum_{P_2 \in \Pcal \mbox{, } P_2 \subset V_s} \frac{\mu_n'(
  f^{-j}(P_2) \cap P_1)}{ \mu_n'(P_1)} \log \left(   \frac{\mu_n'(
  f^{-j}(P_2) \cap P_1)}{ \mu_n'(P_1)} \right).$$
We divide that term into two parts:
$$A= - \sum_{P_1 \in \bigvee_{i=0}^{j-1} f^{-i} \Pcal} \mu_n'(P_1) \sum_{s=0}^{s_0 -1}
  \sum_{P_2 \in \Pcal \mbox{, } P_2 \subset V_s} \frac{\mu_n'(
  f^{-j}(P_2) \cap P_1)}{ \mu_n'(P_1)} \log \left(   \frac{\mu_n'(
  f^{-j}(P_2) \cap P_1)}{ \mu_n'(P_1)} \right),$$
and:
$$B= - \sum_{P_1 \in \bigvee_{i=0}^{j-1} f^{-i} \Pcal} \mu_n'(P_1) \sum_{s=s_0}^{+ \infty}
  \sum_{P_2 \in \Pcal \mbox{, } P_2 \subset V_s} \frac{\mu_n'(
  f^{-j}(P_2) \cap P_1)}{ \mu_n'(P_1)} \log \left(   \frac{\mu_n'(
  f^{-j}(P_2) \cap P_1)}{ \mu_n'(P_1)} \right).$$
First, we show that the second term is as small as we want if we take   
$s_0$ large enough then $n$ large enough. We will deal with $A$ after that.

We have:
$$B=- \sum_{s=s_0}^{+ \infty}  \sum_{P_2 \in \Pcal \mbox{,
  } P_2 \subset V_s}  \left( \sum_{P_1 \in \bigvee_{i=0}^{j-1} f^{-i} \Pcal}
 \mu_n'(P_1) \phi \left( \frac{\mu_n'( f^{-j}(P_2) \cap
  P_1)}{ \mu_n'(P_1)} \right) \right)$$
where $\phi(x)= x \log(x)$. Since that function is convex on $[0 , +
  \infty[$, we deduce:
$$B \leq  -  \sum_{s=s_0}^{+ \infty}  \sum_{P_2 \in \Pcal \mbox{,
  } P_2 \subset V_s} \phi \left(  \sum_{P_1 \in \bigvee_{i=0}^{j-1} f^{-i} \Pcal} \mu_n'(
  f^{-j}(P_2) \cap P_1) \right).$$
That means:
$$B \leq  - \sum_{s=s_0}^{+ \infty}  \sum_{P_2 \in \Pcal \mbox{,
  } P_2 \subset V_s}  \mu_n'( f^{-j}(P_2)) \log ( \mu_n'(
  f^{-j}(P_2)) ).$$
This term is controlled as in the case $j=1$. Indeed,

$$B \leq \sum_{s=s_0}^{+ \infty}   \mu_n'( f^{-j}(V_s))( \log | \Pcal_{r_s} | - \log
  \mu_n'( f^{-j}(V_s))),$$
is smaller than:
\begin{equation*}
\begin{split}
& \log C \mu_n'( f^{-j}( \{ \rho \leq e^{-s_0} \}) )\\
& + 2k \sum_{s=s_0}^{+ \infty} (s+1) \mu_n'( f^{-j}(V_s)) + \sum_{s=s_0}^{+
  \infty}  \mu_n'( f^{-j}(V_s))  \log \left( \frac{1}{ \mu_n'(
  f^{-j}(V_s)) } \right).
\end{split}
\end{equation*}
But since:
$$\sum_{s=s_0}^{+ \infty} s \mu_n'( f^{-j}(V_s)) \leq -  \int_{  \{ \rho \leq e^{-s_0} \} } \log \rho \mbox{ } d( f^j_{*}
  \mu_n')$$
is as small as we want if we take   
$s_0$ large enough then $n$ large enough thanks to the previous lemma, we have that $B$ 
is as small as we want using as for $j=1$ the proof of Lemma \ref{Mane}.

To conclude, it remains to deal with

$$A= - \sum_{P_1 \in \bigvee_{i=0}^{j-1} f^{-i} \Pcal} \mu_n'(P_1) \sum_{s=0}^{s_0 -1}
  \sum_{P_2 \in \Pcal \mbox{, } P_2 \subset V_s} \frac{\mu_n'(
  f^{-j}(P_2) \cap P_1)}{ \mu_n'(P_1)} \log \left(   \frac{\mu_n'(
  f^{-j}(P_2) \cap P_1)}{ \mu_n'(P_1)} \right).$$

We divide the sum $\sum_{P_1 \in \bigvee_{i=0}^{j-1} f^{-i} \Pcal} $ into two parts: 
$$\sum_{P_1 \in
  \bigvee_{i=0}^{j-1} f^{-i} \Pcal \mbox{, } P_1 , \dots , f^{j-1}(P_1) \subset
  \cup_{s=0}^{s_1 -1} V_s }$$
and
$$\sum_{P_1 \in
  \bigvee_{i=0}^{j-1} f^{-i} \Pcal \mbox{, } \exists l \in [0 , j-1] \mbox{, } f^{l}(P_1) \nsubseteq \cup_{s=0}^{s_1 -1} V_s }.$$

The first sum is finite, so:
$$ \sum_{ 
 \begin{array}{c}
 \scriptstyle{P_1 \in \bigvee_{i=0}^{j-1} f^{-i} \Pcal} \\
 \scriptstyle{P_1,  \dots , f^{j-1}(P_1) \subset
  \cup_{s=0}^{s_1 -1} V_s} 
\end{array}
}  \mu_n'(P_1) \sum_{s=0}^{s_0 -1} \sum_{P_2 \in \Pcal \mbox{, } P_2 \subset V_s} \frac{\mu_n'(
  f^{-j}(P_2) \cap P_1)}{ \mu_n'(P_1)} \log \left(   \frac{\mu_n'(
  f^{-j}(P_2) \cap P_1)}{ \mu_n'(P_1)} \right)$$
converges to:
$$ \sum_{
 \begin{array}{c}
 \scriptstyle{P_1 \in \bigvee_{i=0}^{j-1} f^{-i} \Pcal}\\
 \scriptstyle{ P_1 , \dots , f^{j-1}(P_1) \subset
  \cup_{s=0}^{s_1 -1} V_s }
\end{array}
}  \mu(P_1) \sum_{s=0}^{s_0 -1}
  \sum_{P_2 \in \Pcal \mbox{, } P_2 \subset V_s} \frac{\mu(
  f^{-j}(P_2) \cap P_1)}{ \mu(P_1)} \log \left(   \frac{\mu(
  f^{-j}(P_2) \cap P_1)}{ \mu(P_1)} \right)$$
when $n$ goes to $\infty$.

Now, the second sum is less than:
$$- \sum_{l=0}^{j-1} \sum_{
 \begin{array}{c}
 \scriptstyle{P_1 \in  \bigvee_{i=0}^{j-1} f^{-i} \Pcal}\\
\scriptstyle{ f^{l} (P_1) \subset \cup_{s=s_1}^{+ \infty} V_s
  }
\end{array}
}   \mu_n'(P_1) \sum_{s=0}^{s_0 -1}
  \sum_{P_2 \in \Pcal \mbox{, } P_2 \subset V_s} \frac{\mu_n'(
  f^{-j}(P_2) \cap P_1)}{ \mu_n'(P_1)} \log \left(   \frac{\mu_n'(
  f^{-j}(P_2) \cap P_1)}{ \mu_n'(P_1)} \right)$$ 
(we might have add some $\geq0$ terms since $-x
  \log(x) \geq 0$ on $[0,1]$).

Furthermore, since the function $- x \log(x)$ is smaller 
than $e^{-1}$, we deduce that this term is less than:
$$C(s_0) \sum_{l=0}^{j-1} \mu_n'( f^{-l}( \{ \rho \leq e^{- s_1} \}
)),$$
which is also as small as we want if we take $s_1$ 
large enough with respect to $s_0$ then $n$ large enough.

Finally, up to replacing  $\mu_n'$ by $\mu$ in what we just did,
we also have that $ H_{\mu} (f^{-j} (\Pcal) |   \bigvee_{i=0}^{j-1} f^{-i} \Pcal)$
is as close as we want to:
$$- \sum_{
 \begin{array}{c}
 \scriptstyle{P_1 \in \bigvee_{i=0}^{j-1} f^{-i} \Pcal}\\
\scriptstyle{ P_1 , \dots , f^{j-1}(P_1) \subset
  \cup_{s=0}^{s_1 -1} V_s } 
\end{array}
} \mu(P_1)  \sum_{s=0}^{s_0 -1}
  \sum_{P_2 \in \Pcal \mbox{, } P_2 \subset V_s} \frac{\mu(
  f^{-j}(P_2) \cap P_1)}{ \mu(P_1)} \log \left(   \frac{\mu(
  f^{-j}(P_2) \cap P_1)}{ \mu(P_1)} \right)$$
if we take $s_0$ large enough then $s_1$ large with respect to $s_0$.

So we have proved that $H_{\mu_n'} (f^{-j} (\Pcal) | \bigvee_{i=0}^{j-1} f^{-i} \Pcal )$
converges to:
 $$H_{\mu} (f^{-j} (\Pcal) |  \bigvee_{i=0}^{j-1} f^{-i} \Pcal)$$
and that concludes the induction. So the second point of Proposition \ref{passtothelimit} is proved. \\

\noindent {\bf Third point of Proposition \ref{passtothelimit}}

We show that:
$$\frac{1}{q} \sum_{j=0}^{q-1}  \sum_{t \in S_1(j)} \frac{H_{\nu_n}(
f^{-t} \Pcal)}{n-2q+1}$$
converges to $0$ when $n$ goes to $\infty$.

We start by dividing:
\begin{equation*}
\begin{split}
&\frac{1}{q} \sum_{j=0}^{q-1}  \sum_{t \in S_1(j)} \frac{H_{\nu_n}(
f^{-t} \Pcal)}{n-2q+1}\\
&= \frac{1}{q(n-2q+1)} \sum_{j=0}^{q-1}\sum_{t \in S_1(j)} \sum_{s=0}^{+ \infty} - \sum_{P \in
  \Pcal \mbox{, } P \subset V_s} \nu_n(f^{-t}(P)) \log
  \nu_n(f^{-t}(P))
\end{split}
\end{equation*}
into two parts:
$$\frac{1}{q(n-2q+1)} \sum_{j=0}^{q-1} \sum_{t \in S_1(j)} \sum_{s=0}^{s_0 -1 } - \sum_{P \in
  \Pcal \mbox{, } P \subset V_s} \nu_n(f^{-t}(P)) \log  \nu_n(f^{-t}(P)) $$
and:
$$\frac{1}{q(n-2q+1)} \sum_{j=0}^{q-1}\sum_{t \in S_1(j)} \sum_{s=s_0}^{+ \infty} - \sum_{P \in
  \Pcal \mbox{, } P \subset V_s} \nu_n(f^{-t}(P)) \log  \nu_n(f^{-t}(P)). $$
For $s_0$ fixed, the first term goes to $0$ when $n$ goes to $\infty$ since the function $-x
\log x$ is bounded by $e^{-1}$ et since there are only a finite number
  of terms.

For the second term, we remark that it is less than (see again the proof of Proposition
\ref{partition}):
$$\frac{1}{q(n-2q+1)} \sum_{j=0}^{q-1} \sum_{t \in S_1(j)} \sum_{s=s_0}^{+ \infty} \nu_n(f^{-t}( V_s))( \log | \Pcal_{r_s}
| - \log \nu_n (f^{-t}( V_s)) )$$
which is in turn less than:
\begin{equation*}
\begin{split}
&\frac{1}{q(n-2q+1)} \sum_{j=0}^{q-1} \sum_{t \in S_1(j)} ( \log C  + 2k
\sum_{s=s_0}^{+ \infty} (s+1)  \nu_n(f^{-t}( V_s))\\
& +\sum_{s=s_0}^{+ \infty} \nu_n(f^{-t}( V_s)) \log
\left(\frac{1}{\nu_n(f^{-t}( V_s))} \right) ).
\end{split}
\end{equation*}
As previously, this term is smaller than:
$$ \frac{1}{q(n-2q+1)} \sum_{j=0}^{q-1} \sum_{t \in S_1(j)} ( \log C  + 2k + (2k+1) \sum_{s=s_0}^{+ \infty} s  \nu_n(f^{-t}( V_s)) +
2 e^{-1} \sum_{s=s_0}^{+ \infty} e^{-s/2} ).$$
But the terms:
$$ \frac{1}{q(n-2q+1)} \sum_{j=0}^{q-1} \sum_{t \in S_1(j)} ( \log C  + 2k)$$
and:
$$\frac{1}{q(n-2q+1)} \sum_{j=0}^{q-1} \sum_{t \in S_1(j)} 2 e^{-1}
\sum_{s=s_0}^{+ \infty} e^{-s/2} $$
converge to  $0$ when $n$ goes to $\infty$ (because the cardinality of 
$S_1(j)$ is smaller than $3q$). It remains to control:
$$\frac{2k+1}{q(n-2q+1)} \sum_{j=0}^{q-1} \sum_{t \in S_1(j)} \sum_{s=s_0}^{+ \infty} s  \nu_n(f^{-t}( V_s)).$$
This term is equal to:
\begin{equation*}
\begin{split}
&\frac{2k+1}{q(n-2q+1)} \sum_{j=0}^{q-1} \sum_{t \in S_1(j)} \sum_{s=s_0}^{+ \infty} s  ((f^t)_{*} \nu_n)(
  V_s) \\
&\leq - \frac{2k+1}{q(n-2q+1)} \sum_{j=0}^{q-1} \sum_{t \in S_1(j)} \int_{
  \{ \rho \leq e^{-s_0} \} } \log
\rho \mbox{ } d (f^t)_{*} \nu_n.
\end{split}
\end{equation*}
But:
$$\sum_{t \in S_1(j)}  (f^t)_{*}  \nu_n 
\leq  (\phi(n)m-m+1)\mu_n''$$
with $\mu_n''= \frac{1}{\phi(n)m-m+1} \sum_{p=0}^{\phi(n)m-m} f^p_{*} \nu_n$.  

Following Lemma \ref{estimation_integrale} with $\mu_n''$ instead of $\mu_n'$ and
$i=0$ (this is indeed possible because the indexes $p$ in $\mu_n''$ goes to 
$\phi(n)m-m$ which is $\leq n-1-(m-1)$), we deduce that:
$$\frac{2k+1}{q(n-2q+1)} \sum_{j=0}^{q-1} \sum_{t \in S_1(j)} \sum_{s=s_0}^{+ \infty} s  ((f^t)_{*} \nu_n)(
  V_s)$$
is as small as we want by taking $s_0$ large enough then $n$ large enough.

This gives the third point of Proposition \ref{passtothelimit}.

$$ $$

{\bf Fourth point of Proposition \ref{passtothelimit}}

The proof is the same than for the third point replacing $S_1(j)$ by $S_2(j)$, $\Pcal$ by $\Qcal$, $\rho$ by
$\eta$.

At the end, we have to bound from above:
$$\frac{2k+1}{q(n-2q+1)} \sum_{j=0}^{q-1} \sum_{t \in S_2(j)}
\sum_{s=s_0}^{+ \infty} s  \nu_n(f^{-t}( V_s)),$$

(here the $V_s$ correspond to the partition $\Qcal$ and to the function
$\eta$).

That term is equal to:
\begin{equation*}
\begin{split}
&\frac{2k+1}{q(n-2q+1)} \sum_{j=0}^{q-1} \sum_{t \in S_2(j)} \sum_{s=s_0}^{+ \infty} s  ((f^t)_{*} \nu_n)(
  V_s) \\
&\leq - \frac{2k+1}{q(n-2q+1)} \sum_{j=0}^{q-1} \sum_{t \in S_2(j)} \int_{
  \{ \eta \leq e^{-s_0} \} } \log
\eta \mbox{ } d (f^t)_{*} \nu_n.
\end{split}
\end{equation*}
Finally:
$$\sum_{t \in S_2(j)}  (f^t)_{*}  \nu_n 
\leq n \mu_n$$
and since $\int _{  \{ \eta \leq e^{-s_0} \} } \log \eta d \mu_n$
converges to $\int _{  \{ \eta \leq e^{-s_0} \} } \log \eta d \mu$,
it is also as small as we want if $s_0$ is large enough then $n$ is large enough.
This gives the fourth point of  Proposition \ref{passtothelimit}, and the proposition follows.    \hfill $\Box$ \hfill

\chapter{Dynamics of birational maps of $\P^k$}\label{Gabriel}

\section{A family of birational maps}
Recall that a meromorphic map $f:\P^k\rightarrow\P^k$ is holomorphic
outside an analytic subset $I(f)$ of codimension $\geq 2$ in $\P^k$. Let
$\Gamma$ denote the
closure of the graph of the restriction of $f$ to $\P^k\setminus
I(f)$. This is an irreducible analytic set of dimension $k$ in
$\P^k\times\P^k$. 

Let $\pi_1$ and $\pi_2$ denote the canonical
projections of $\P^k\times \P^k$ on its factors. 
The {\it indeterminacy locus} $I(f)$ of $f$ is the set of points $z\in\P^k$ such that $\dim
\pi_1^{-1}(z)\cap\Gamma\geq 1$. We assume that $f$ is {\it dominant},
that is, $\pi_2(\Gamma)=\P^k$. The {\it second indeterminacy set} of
$f$ is the set $I'(f)$  of points $z\in\P^k$ such that $\dim
\pi_2^{-1}(z)\cap\Gamma\geq 1$. 
If $A$ is a subset of $\P^k$, define
$$f(A):=\pi_2(\pi_1^{-1}(A)\cap\Gamma)\quad \mbox{and}\quad
f^{-1}(A):=\pi_1(\pi_2^{-1}(A)\cap\Gamma).$$
Define formally for a current $S$ on $\P^k$, not necessarily positive
or closed, the pull-back $f^*(S)$ by
\begin{equation} \label{eq_pullback_def}
f^*(S):=(\pi_1)_*\big(\pi_2^*(S)\wedge [\Gamma]\big)
\end{equation}
where $[\Gamma]$ is the current of integration of $\Gamma$.

Similarly, the operator $f_*$ is formally defined by 
\begin{equation}  \label{eq_pushforward_def}
f_*(R):=(\pi_2)_*\big(\pi_1^*(R)\wedge [\Gamma]\big).
\end{equation}
For $0\leq q\leq k$ and $n>0$, define
$\lambda_q(f^n)$ which gives a size for the action of $f^n$ on the cohomology group $H^{q,q}(\P^k)$ as:
\begin{align}
\lambda_q(f^n) & := \|(f^n)^*(\omega^q)\| = \int_{\P^k} (f^n)^*(\omega^q)\wedge\omega^{k-q} \nonumber \\
               & = \|(f^n)_*(\omega^{k-q})\|=\int_{\P^k} (f^n)_*(\omega^{k-q})\wedge \omega^q.
\end{align} 
We have in particular that $\lambda_1(f)=d$ is the algebraic degree. We define the {\it dynamical degree of order $q$} of
$f$ by:
\begin{align}
d_q & :=  \lim_{n\rightarrow\infty} (\lambda_q(f^n))^{\frac{1}{n}}
\end{align}
These limits always exist and $d_q\leq d_1^q$ \cite{DS9}.
The last degree $\lambda_k(f)=d_k$ is {\it the topological degree} of $f$.
It is equal to $\#f^{-1}(z)$ for $z$ generic. A result by Gromov \cite[Theorem 1.6]{Gromov}
implies that  $q\mapsto \log d_q$ is concave in $q$. In particular, there exists a $q_0$ such that:
$$ 1= d_0 \leq d_1 \leq \dots \leq d_{q_0} \geq  \dots \geq  d_k.$$ 
Of course, $q_0$ can be equal to $k$ which is the case for holomorphic endomorphisms of $\P^k$. \\

Here, we consider a birational map $f$ of algebraic degree $d\geq 2$. That is a map such that $\#f^{-1}(z)=1$ for $z$ generic ($d_k=1$). Let $\delta$ be the algebraic degree of $f^{-1}$ and denote by $\lambda^-_q(f^n)$ and $d^-_q$ the quantities previously defined for $f$. We denote $I^+:=I(f)$ and $I^-=I'(f)=I(f^{-1})$ the indeterminacy sets of $f$ and $f^{-1}$. \\

We also consider the critical sets $\mathcal{C}^+$ (or $\C(f)$) and $\mathcal{C}^-$ (or $\C(f^{-1})$) defined by:
\begin{align*}
 \mathcal{C}^+&:= f^{-1}(I^-) \\
 \mathcal{C}^-&:=(f^{-1})^{-1}(I^+).
 \end{align*}
Write $f=[P_1:\dots:P_{k+1}]$ where the $P_i$ are homogeneous polynomials of degree $d$. Let $F=(P_1,\dots,P_{k+1})$
be the induced map on $\mathbb{C}^{k+1}$. 
Similarly, write $f^{-1}=[Q_1:\dots: Q_{k+1}]$ where the $Q_i$ are homogeneous polynomials of degree $\delta$ and let $F^{-1}=(Q_1,\dots,Q_{k+1})$. There is of course an abuse of notation since $F \circ F^{-1} \neq Id $ instead, we have that:
$$F \circ F^{-1} = P(z_1,\dots,z_{k+1}) \times \left(z_1,\dots,z_{k+1}\right),$$
where $P$ is a homogeneous polynomial of degree $d \delta-1$ equal to $0$ in $\pi^{-1}(\C^-\cup I^-) $ where $\pi : \mathbb{C}^{k+1} \to \P^k$ is the canonical projection.  That implies that the critical set $\C^-$ is an analytic set of codimension 1 and that $I^- \subset \C^-$. Similarly, we have that $\C^+$ is of codimension $1$ and $I^+ \subset \C^+$ (see also Proposition 3.3 in \cite{Dil} and \cite{Sib}).
So, $f:\P^k \backslash  \mathcal{C}^+ \rightarrow \P^k \backslash  \mathcal{C}^-$ is a biholomorphism. \\

Let $s$ be such that $\text{dim} (I^+)=k-s-1$, then we have the proposition (see also Proposition 2.3.2 in \cite{Sib}): 
\begin{proposition}
For any birational map $f$, we have that:
$$\text{dim}(I^+)+\text{dim}(I^-)\geq k-2 .$$ 
Furthermore, we have $\lambda_q(f)=d^q$ for $q\leq s$ and $\lambda^-_q(f)=\delta^q$ for $q \leq k-\text{dim}(I^-)-1$. In particular, if $I^+\cap I^-=\varnothing$, then  $s-1\leq \text{dim}(I^-)\leq s$.
\end{proposition}
\emph{Proof.} For $q \leq s$, the $2k-2q+1$-Hausdorff dimension of $I^+$ is equal to $0$ since the real dimension of $I^+$ is $2k-2s-2$. Hence, $f^*(\omega)^{q}$ is well defined and its mass is equal to the mass of $f^*(\omega)$ at the power $q$, that is $d^q$  (see Corollary 4.11 in \cite[Chapter III]{dem2} or \cite{Sib}). On the other hand, the currents $(f^*(\omega))^q$ and $f^*(\omega^q)$ are equal outside $I^+$ and they cannot give mass to algebraic sets of dimension $\leq k-1$ thus they cannot give mass to $I^+$. In particular, these currents are equal and the mass of $f^*(\omega^q)$ is equal to $d^q$. 

 So $\lambda_q(f)=d^q$ for $q\leq s$.  This implies that $\lambda_q^-(f)=d^{k-q}$ for $q \geq k-s$. We have proved that $\lambda_q(f)$ is increasing up to the rank $k-1-\text{dim}(I^+)$. Applying the same argument 
to $f^{-1}$ gives that $\lambda_q(f^{-1})$ is increasing up to 
the rank $k-1-\text{dim}(I^-)$. Hence we have $k-s\geq k-1-\text{dim}(I^-)$. 
So  $\text{dim}(I^+)+\text{dim}(I^-)\geq k-2 $.
 
 For $f^{-1}$, we have as for $f$ that  $\lambda_q(f^{-1})=\delta^q$ for $q \leq k-\text{dim}(I^-)-1$. \hfill $\Box$ \hfill \\

In all the cases studied, one has $\text{dim}( I^-)=s-1$. In the case of regular automorphisms of $\mathbb{C}^k$, this is because the indeterminacy sets lie on the line at infinity which is isomorphic to $\mathbb{P}^{k-1}$ (hence $\text{dim}( I^+)+\text{dim}( I^-)=k-2$). For $k= 2$, $\text{dim}( I^-)=s$ is impossible since the indeterminacy sets are of codimension $\geq 2$, which means that they are points. Finally, in \cite{DS10}, the hypothesis of $s$-pseudoconvexity of some neighborhood of $I^+$ implies  $\text{dim}(I^-)=s-1$. 

Still, this is not always the case. Take for example $f$ in $\P^3$ given by $[yz:xz:zt+y^2:z^2]$ then $f$ is birational with inverse  $f^{-1}=[yt:xt:t^2:zt-x^2]$. Then $I(f)=\{y=0\}\cap \{z=0\}$ and $I(f^{-1})=\{x=0 \}\cap \{t=0 \}$. So they are both of dimension $1$ and $I(f)\cap I(f^{-1})=\varnothing$. \\

So we need to formulate a hypothesis: from now on, we are going to assume that $I^+$ and $I^-$ are of pure dimension and satisfies 
\begin{align}
\text{dim}( I^+)=k-s-1 \quad \text{and} \quad \text{dim}(I^-) =s-1.
\end{align}
for $s$ with $1\leq s\leq k-1$. \\

 In particular, the previous proposition becomes:
\begin{proposition}\label{312}
Let $f$ be as above, then we have $\lambda_q(f)=d^q$ for $q\leq s$ and $\lambda_q(f)=\delta^{k-q}$ for $q \geq s$. In particular, $d^s=\delta^{k-s}$.
\end{proposition}
We introduce some notations. Let $\C_s$ denote the convex compact set of (strongly) positive closed currents
$S$ of bidegree $(s,s)$ on $\P^k$ and of mass 1, i.e. $\|S\|:=\langle
S,\omega^{k-s}\rangle =1$.
For a positive closed current $T$ of mass $m(T)>0$, 
we denote by $T_{\text{nor}}$ the renormalization of $T$ (that is $T_{\text{nor}}=m(T)^{-1}T$). 
Denote for simplicity $L:=\lambda_q(f)^{-1}f^*$ and $\Lambda:=(\lambda_{k-q})^{-1}f_*=(\lambda^-_q(f))^{-1} f_*$ 
which are well defined operators on the elements of $\C_q$ which 
are smooth near $I^-$ (resp. $I^+$). We make an abuse of notations and write $L$ 
instead of $L_q$, this is not a problem since in what follows $L(S)$ will always be the current $f^*(S)_{\text{nor}}$. The theory of super-potential (see the appendix) allows us to extend the operator $L$ (resp. $\Lambda$) to the currents in $\C_q$ such that their super-potentials are finite at one point of the form $\Lambda(S)$ for $S\in \C_{k-q+1}$ smooth near $I^+$ (resp. at one point of the form $L(S)$ for $S\in \C_{k-q+1}$ smooth near $I^-$).\\

In order to work with the currents in cohomology, we need a hypothesis on the indeterminacy sets so that $(f^n)^*=(f^*)^n$ on the cohomology group $H^{q,q}(\P^k)$. If so, we say that the map is \emph{algebraically $q$-stable} (see \cite{Sib} and \cite{DS6}). 

 We introduce the following condition on $f$:  
\begin{align}\label{weak}
\bigcup_{n\geq 0} f^{-n}I(f) \cap \bigcup_{n\geq 0} f^{n}I(f^{-1}) =\varnothing.
\end{align} 
In the case where $q=1$ and $k=2$, this condition is equivalent to the algebraic stability.

No we show  that a map which satisfies (\ref{weak}) is in fact algebraically $q$-stable for all $q$. That is to say no mass is lost on the indeterminacy set by pull-back. More precisely, we have the proposition that uses the theory of super-potential (see the appendix):
\begin{proposition}\label{s-stability}
Let $f$ be a birational map satisfying (\ref{weak}), then $(f^*)^n=(f^n)^*$ on $\C_q$ for all $q$, $0\leq q \leq k$.
More precisely, $\lambda_q(f^n)= (\lambda_q(f))^n$ so $d_q=\lambda_q(f)$ for all $q$.
\end{proposition}
\emph{Proof.} We have to compute the integral:
\begin{align*}
\lambda_q(f^n) = \|(f^{n})^*(\omega^q)\| = \int_{\P^k} (f^{n})^*(\omega^q)\wedge\omega^{k-q}. 
\end{align*}          
The proof is by induction on $n$:  $(f^{n-1})^*(\omega^q)$ is a form in $L^1$ smooth  near $I^-$ by  (\ref{weak}). 
So we can define its pull-back by $f$ which is of mass $\lambda_q(f)\lambda_{q}(f^{n-1})$. On the other hand, $\omega^q$ is smooth near $I(f^{-n})$ so it is $(f^n)^*$-admissible and the mass is of $(f^n)^*(\omega^q)$ is $\lambda_{q}(f^{n})$. 

We will now prove that $f^*((f^{n-1})^*(\omega^q))=(f^n)^*(\omega^q)$. 

Let ${\pi_1}_{|\Gamma}$ and ${\pi_2}_{|\Gamma}$ be the restriction of $\pi_1$ and $\pi_2$ to the graph $\Gamma$ of $f$. That way, $f^*(S)=(\pi_1)_* ({\pi_2}_{|\Gamma})^*(S)$ where $S\in \C_q$ is smooth near $I^-$. We will take $S=(f^{n-1})^*(\omega^q)_{\text{nor}}$.

 Let $V$ be a small neighborhood of $I^-$ such that $S$ is smooth here. Outside ${\pi_2}^{-1}(V)\cap \Gamma$, ${\pi_2}_{|\Gamma}$ is a finite map, hence ${\pi_2}^*(S)\wedge [\Gamma]$ is well defined and depends continuously of $S$ here by \cite{DS13} (Theorem 1.1). Furthermore, if $S_{|\P^k\backslash V}$ does not give mass to a Borel set $A$ then ${({\pi_2}_{|\Gamma})}^*(S)$ does not give mass to $({\pi_2}_{|\Gamma})^{-1}(A)$ outside $\pi_2^{-1}(V)\cap \Gamma$. Since $\pi_1$ is holomorphic, the same is true for $f^*(S_{|\P^k\backslash V})$.  And on $V$, $S$ is smooth, hence $f^*(S_{|V})$ is a form in $L^1$ (see e.g. \cite{DS9}). \\

We consider $S=(f^{n-1})^*(\omega^q)$: $(f^{n-1})^*(\omega^q)$ is a form in $L^1$ hence it does not give mass to  algebraic sets of dimension $\leq k-1$; so $f^*((f^{n-1})^*(\omega^q))$ is a current that does not give mass to algebraic sets of dimension $\leq k-1$. We obtain that $f^*((f^{n-1})^*(\omega^q))$
and $(f^n)^*(\omega^q)$ are equal wherever they are smooth that is outside analytic sets of dimension $\leq k-1$.  We deduce that they are equal hence they have the same mass.\hfill $\Box$ \hfill \\

The following corollary of the previous proof will be used later:
\begin{corollaire}\label{composition_smooth}
Let $R\in \C_q$ be a smooth form, then for all $i,j\geq 0$, we have that $R$ is $(f^j)^*$-admissible and $(f^j)^*((f^i)^*(R))=(f^{i+j})^*(R)$. 
\end{corollaire}
Let $j \geq 0$ and $q \leq k$. For a current $S \in \C_q$ which is $(f^j)^*$-admissible, we can define $L_j(S)$ as $(\lambda_q(f^j))^{-1}(f^j)^*(S)$ (similarly we define $\Lambda_j$).  By Proposition \ref{s-stability}, we have that $\lambda_q(f^j)=\lambda_q(f)^j$ so we can also write $L_j(S)=\lambda_q(f)^{-j}(f^j)^*(S)$. From Corollary \ref{composition_smooth}, we have that $L^j(S)=L_j(S)$ on smooth forms, the question is: does it also stand for $(f^j)^*$-admissible currents ? The following lemma answers positively:
\begin{lemme}\label{stability}
Let $S\in \C_q$ for $q\leq k$. Let $n>0$ such that $S$ is $(f^n)^*$-admissible then for all $j$ with  $0\leq j\leq n-1$, $L^j(S)$ is well defined, $f^*$-admissible and  $L^{j+1}(S)=L_{j+1}(S)$. In particular, $L^n(S)=L_n(S)$. 
\end{lemme} 
\emph{Proof.} Let $S$ be as above, then a super-potential of $S$ is finite at $\Lambda_{n}(\omega^{k-q+1})$ by hypothesis. Since $f$ satisfies (\ref{weak}), we have that $\Lambda(\omega^{k-q+1})$ is smooth near $I(f^{n-1})$, the previous corollary implies that  $$\Lambda_{n}(\omega^{k-q+1})=\Lambda_{n-1}(\Lambda(\omega^{k-q+1})).$$
  So the super-potentials of $S$ are finite at the image by $\Lambda_{n-1}$ of a current smooth near $I(f^{n-1})$: it is $(f^{n-1})^*$-admissible (see the appendix). An immediate induction gives that $S$ is $(f^j)^*$-admissible for $j\leq n$. \\

Now we prove by induction on $j$ that $L^j(S)$ is $f^*$-admissible and that $L^{j+1}(S)=L_{j+1}(S)$. For $j=0$, it is just the fact that $S$ is $f^*$-admissible.  Now, assume the property holds up to the rank $j$. A super-potential of $L^j(S)=L_j(S)$ is by Proposition \ref{prop_pull_general}:
$$\U_{L_j(S)}= \U_{L_j(\omega^q)}+\frac{\lambda_{q-1}(f^j)}{\lambda_{q}(f^j)} \U_{S}\circ \Lambda_j$$ 
on forms smooth near $I(f^{j})$. Taking the value at $\Lambda(\omega^{k-q+1})$ (which is smooth near $I(f^{j})$) gives:
$$\U_{L_j(S)}(\Lambda(\omega^{k-q+1}))= \U_{L_j(\omega^q)}(\Lambda(\omega^{k-q+1}))+\frac{\lambda_{q-1}(f^j)}{\lambda_{q}(f^j)} \U_{S}(\Lambda_j(\Lambda(\omega^{k-q+1}))).$$ 
The current $L_j(\omega^q)=L^j(\omega^q)$ is $f^*$-admissible since it is smooth near $I^-$, that means that $\U_{L_j(\omega^q)}(\Lambda(\omega^{k-q+1}))$ is finite. Similarly, applying the previous corollary to $f^{-1}$ gives that $\Lambda_j(\Lambda(\omega^{k-q+1}))=\Lambda^{j+1}(\omega^{k-q+1})$ and since $S$ is $(f^{j+1})^*$-admissible then $\U_{S}(\Lambda_j(\Lambda(\omega^{k-q+1})))$ is also finite. 

 That gives that $\U_{L_j(S)}(\Lambda(\omega^{k-q+1}))$ is finite so $L_j(S)$ is $f^*$-admissible. We can now apply Proposition \ref{prop_pull_general} to $L^{j}(S)$:
 \begin{align*}
 \U_{L^{j+1}(S)} &=  \U_{L(\omega^q)}+ \frac{\lambda_{q-1}(f)}{\lambda_{q}(f)}\U_{L^{j}(S)}\circ \Lambda \\ 
                 &=  \U_{L(\omega^q)}+ \frac{\lambda_{q-1}(f)}{\lambda_{q}(f)}( \U_{L_j(\omega^q)}+\frac{\lambda_{q-1}(f^j)}{\lambda_{q}(f^j)} \U_{S}\circ \Lambda_j)\circ \Lambda                  
 \end{align*}
on smooth forms. Since $\U_{L(\omega^q)}+ \frac{\lambda_{q-1}(f)}{\lambda_{q}(f)} \U_{L_j(\omega^q)} \circ \Lambda = \U_{L^{j+1}(\omega^q)}$ on smooth forms, and since $\Lambda_j\circ \Lambda =\Lambda_{j+1} $ on smooth forms, we deduce from Proposition \ref{s-stability} that:
$$\U_{L^{j+1}(S)}=\U_{L_{j+1}(S)}$$
on smooth forms, hence  $L^{j+1}(S)=L_{j+1}(S)$ by Proposition \ref{prop_unique_sp}. That gives the lemma. \hfill $\Box$ \hfill \\

As usual, for two sets $E$ and $F$, we denote $\inf_{x\in E,y\in F} \text{dist}(x,y)$ by $\text{dist}(E,F)$. In \cite{BD1}, the authors asked for a quantitative and stronger version of (\ref{weak}) similar to:
\begin{Hypothesis}\label{distancefort}
The birational mapping $f$ satisfies:
\begin{align*}
&\sum^{\infty}_{n=0}      \left(\frac{1}{d}\right)^n\log \text{dist}(I^+,f^n(I^-)) >-\infty \\
& \qquad \qquad \qquad \mbox{and} \\
&\sum^{\infty}_{n=0}  \left(\frac{1}{\delta}\right)^n\log \text{dist}(I^-,f^{-n}(I^+)) >-\infty
\end{align*}
\end{Hypothesis}
In $\P^2$, $\sum^{\infty}_{n=0} d^{-n}\log \text{dist}(I^+,f^n(I^-)) >-\infty $ and $\sum^{\infty}_{n=0} \delta^{-n}\log \text{dist}(I^-,f^{-n}(I^+)) >-\infty$ are equivalent (see \cite{Dil}), it has no reason to be true in higher dimension. \\

Let $f$ be a birational map satisfying $\text{dim}( I^+)=k-s-1$ and $\text{dim}(I^-) =s-1$. Recall that a quasi-potential of a current $T\in \C_q$ is a current $U$ of bidegree $(q-1,q-1)$ such that $T=\omega^q+dd^c U$. We know from the appendix that it is always possible to take $U$ negative. Here, we will use a hypothesis that differs a bit from \ref{distancefort}. In what follows, for an irreducible analytic set $A$, we define $[f(A)]$ as the current of integration over $f(A)$ counting the multiplicity of $f$ at $A$ and if $A$ is not irreducible, we decompose it into irreducible components $(A_i)$ and we define $[f(A)]$ as $\sum_i [f(A_i)]$.  

Assume that $I^+ \cap f^j(I^-)=\varnothing$ for $j\leq n$, then $f^n(I^-)$ is well defined and the form $U_{L(\omega)}L(\omega^{s-1})$ is smooth in $f^n(I^-)$ so the following integral is well defined:
$$\int_{[f^n(I^-)]} U_{L(\omega)}L(\omega^{s-1}).$$
The terms $(\text{deg}(I^-))^{-1} $ and $(\text{deg}(I^+))^{-1}$ in the following hypothesis are just here to normalize the integrals. 
\begin{Hypothesis}\label{distance}
Let $f$ be a map satisfying (\ref{weak}). Let $U_{L(\omega)}$ be a negative quasi-potential of $L(\omega)$ and let $U_{\Lambda(\omega)}$ be a negative quasi-potential of $\Lambda(\omega)$. The birational mapping $f$ satisfies:
\begin{align*}
&\sum^{\infty}_{n=0} \left(\frac{1}{d^s}\right)^n (\text{deg}(I^-))^{-1}\int_{f^n(I^-)} U_{L(\omega)}L(\omega^{s-1}) >-\infty \\
& \qquad \qquad \qquad \text{and} \\
&\sum^{\infty}_{n=0} \left(\frac{1}{\delta^{k-s}}\right)^n (\text{deg}(I^+))^{-1} \int_{f^{-n}(I^+)} U_{\Lambda(\omega)}\Lambda(\omega^{k-s-1}) >-\infty
\end{align*}
\end{Hypothesis}
 In the case of $\P^2$, it is equivalent to Hypothesis \ref{distancefort} (see \cite[Theorem 4.3]{BD1} and  Theorem \ref{convergencecurrent} below). That is because the distance between the supports of the currents is a good distance for Dirac masses but not for currents of higher bidimension.
We will see in Theorem \ref{convergencecurrent} that Hypothesis \ref{distancefort} implies Hypothesis \ref{distance}. 

We can apply Proposition \ref{s-stability} to a map satisfying Hypothesis \ref{distance}. 
We will see in Theorem \ref{convergencecurrent} that Hypothesis \ref{distance} has a clear interpretation in term of super-potentials (it means that the super-potential of the Green current of order $s$ is finite at $[I^-]_{\text{nor}}$). Its interest is that it is generic (see Theorem \ref{generic}) so that we can construct many examples. \\

In what follows, by positive closed currents we mean strongly positive closed currents. So inequalities between positive closed currents have to be understood in the strong sense namely $S_1\leq  S_2$ means that $S_2-S_1$ is a strongly positive closed current.  \\

We sum up the setting we are in. From now on we consider a birational map $f$ of $\P^k$ with:
\begin{itemize}
\item $\text{dim}( I^+)=k-s-1$ and $\text{dim}(I^-) =s-1$ for some $s$ with $1\leq s\leq k-1$.
\item The map $f$ satisfied  Hypothesis \ref{distance}.
\end{itemize} 
Observe that the set of maps which satisfy those conditions is stable by iteration.

 \section{Construction of the Green currents}
 Recall that we assume that Hypothesis (\ref{distance}) holds for $f$. Using Propositions \ref{312} and \ref{s-stability}, we have that $f$ is \emph{algebraically $q$-stable} for all $q$ and for $q\leq s$, we have $\lambda_q(f^n)=(d^q)^n$ for all $n$ so $d_q=d^q$.
 
 Let $q\leq s$. Recall that for $S\in \C_q$ which is $f^*$-admissible, $L(S)$ is the element of $\C_q$ defined as $d^{-q}f^*(S)$. Furthermore, any current smooth in a neighborhood of $I^-$ is $f^*$-admissible. By Proposition \ref{s-stability}, $L^{n-1}(\omega^q)$ is $f^*$-admissible since $f$ satisfies (\ref{weak}) so we can define $L^n(\omega^q)$. 
  
   Now, let $\U_{L(\omega^q)}$ denote a negative super-potential of $L(\omega^q)$ (it is always possible by Proposition \ref{resolutionquasipotential}). 
  
 So, we have that by Proposition \ref{prop_pull_general} for $m>0$, a super-potential of $L^m(\omega^q)=L(L^{m-1}(\omega^q))$ is given on currents in $\C_{k-q+1}$ smooth in a neighborhood of $I^+$ by:
$$ \U_{L(\omega^q)}+\frac{1}{d} \U_{L^{m-1}(\omega^q)}\circ\Lambda.$$
So, by induction, for an element $R\in \C_{k-q+1}$ such that $\Lambda^n(R)$ is smooth near $I^+$ for all $0\leq n\leq m-1$, we have that a super-potential of $L^m(\omega^q)$ is given on $R$ by:
\begin{equation}\label{sum_green}
\sum_{n= 0}^{ m-1} \left(\frac{1}{d}\right)^n \U_{L(\omega^q)} \circ \Lambda^{n}(R).
\end{equation}
In particular, by (\ref{weak}), we have that for smooth forms in $\C_{k-q+1}$, a super-potential of  $L^m(\omega^q)$ is given by (\ref{sum_green}).  \\

Since the sequence is decreasing, it is enough to show that it does not converge uniformly to $-\infty$ to show that it converges in the Hartogs' sense (see Proposition \ref{cor_decreasing_sqp}). In \cite{DS6}, the authors prove that fact in the algebraically $q$-stable case with an additional assumption on the size of $\C^+$ (that fails in our case) using the fact that the sequence is bounded from below by the super-potential of any weak limit of the Cesar\`o mean  of $(L^m(\omega^q))$. Here the idea  is to show that the convergence holds at the point $[I^-]_{\text{nor}}$. We also prove that Hypothesis \ref{distancefort} implies Hypothesis \ref{distance}.

We need the following estimate of $\U_{L(\omega)}$ for that. It is similar to Proposition 1.3 in \cite{BD1} though our proof is simpler taking advantage of the fact that we are in $\P^k$.
\begin{lemme}\label{estimate}
Let $U_{L(\omega)}$ be a quasi-potential of $L(\omega)$. Then there exist constants $A>0$, $B$, $A'>0$, $B'$ such that:
\begin{equation}\label{encadrementL}
  A\log\text{dist}(x,I^+)-B \leq U_{L(\omega)}(x)\leq A'\log\text{dist}(x,I^+)+B', 
 \end{equation}
for all $x$. 
\end{lemme}
\emph{Proof of the lemma.} Let $P_1,\dots, P_{k+1}$ be homogeneous polynomials of degree 
$d$ with no common factors such that $f=[P_1:\dots:P_{k+1}]$. That way, $I^+=\{P_1=\dots=P_{k+1}=0\}$.
For an element $Z=(z_1,\dots,z_{k+1})\in \mathbb{C}^{k+1}$, we write $|Z|^2=|z_1|^2+\dots+|z_{k+1}|^2$. Let $\pi:\mathbb{C}^{k+1}\to \P^k$ denote the canonical projection and we write $F=(P_1,\dots,P_{k+1})$. Then, we have that:
$$ \pi^*( L(\omega)-\omega)= dd^c (\frac{1}{d}\log |F|^2-\log |Z|^2).$$   
Observe that the qpsh function $d^{-1}\log |F|^2-\log |Z|^2$ is well defined on $\P^k$ since it does not depend on the choice of coordinates. So we can write that $U_{L(\omega)}= d^{-1}\log |F|^2-\log |Z|^2$. So in $\P^k$, the singularities of $U_{L(\omega)}$ come from the terms in  $\log |F|^2$. In the open set of $\mathbb{C}^{k+1}$ defined by $1-\varepsilon<|Z|<1+\varepsilon$, we have that the map $F(Z)$ is equal to $(0,\dots,0)\in \mathbb{C}^{k+1}$ exactly in $\pi^{-1}(I^+)$. Using Lojasiewicz Theorem (Chapter IV Theorem 7 in \cite{Loj}), that provide us two constants $\alpha>0$ and $C>0$ such that on $|Z|=1$ we have:
$$ |F(Z)| \geq C (\text{dist}(Z,\pi^{-1}(I^+)))^{\alpha}.$$
Now from the fact that the projection $\pi$ is Lipschitz in $|Z|=1$ and the above bound, we have constants $A>0$, $B$ such that:
 $$U_{L(\omega)} \geq A\log \text{dist}(. \ ,I^+)-B .$$
For the other inequality, we work in a chart of $\P^k$ where we let $z$ be the coordinate. Let $V$ be a relatively compact open set in the chart. Observe that it is sufficient to prove the upper bound in $V$. For $y\in I^+$ in the chart, we have that $|F(z)|^2=||F(z)|^2-|F(y)|^2|$ is less than $C' \text{dist}(z,y)$. Taking the infimum over all such $y$, we get that $|F(z)|^2$ is less than $C' \text{dist}(z,I^+)$. 
Taking the logarithm gives the estimate in $V$ and the lemma follows. \hfill $\Box$ \hfill 

\begin{theorem}\label{convergencecurrent}
The sequence $(L^m(\omega^s))$ converges in the Hartogs' sense 
to the \emph{Green current of order $s$ of $f$} that we 
denote by $T^+_s$. More precisely, for an application satisfying (\ref{weak}), the first half 
of Hypothesis \ref{distance}:
$$\sum_{n} \left(\frac{1}{d^s}\right)^n (\text{deg}(I^-))^{-1}\int_{f^n(I^-)} U_{L(\omega)}L(\omega^{s-1}) >-\infty$$
is equivalent to the fact that the sequence:
$$\left(\sum_{n= 0}^{m} \left(\frac{1}{d}\right)^n \U_{L(\omega^s)}(\Lambda^{n}([I^-]_{\text{nor}} ))\right),$$
converges. That is to say $\U_{T^+_s}([I^-]_{\text{nor}}) > -\infty$. 

Finally, any map satisfying Hypothesis \ref{distancefort} satisfies Hypothesis \ref{distance}. For those maps, we have for $q\leq s$ that $(L^m(\omega^q))$ also converges in the Hartogs' sense to the \emph{Green current of order $q$ of $f$} that we denote by $T^+_q$.
\end{theorem} 
\emph{Proof of the theorem.} By hypothesis $f$ satisfies (\ref{weak}). So, $\Lambda^n([I^-]_{\text{nor}})\in \C_{k-s+1}$ is smooth in a neighborhood of $I^+$ for all $n$ and  $\Lambda^n([I^-]_{\text{nor}})=[f^n(I^-)]_{\text{nor}}$ (counting the multiplicity). So we have that a super-potential of $L^m(\omega^q)$ is given on $[I^-]_{\text{nor}}$ by (\ref{sum_green}):
 $$\U_{L^m(\omega^s)}([I^-]_{\text{nor}})=\sum_{n= 0}^{ m-1} \left(\frac{1}{d}\right)^n \U_{L(\omega^s)}(\Lambda^{n}([I^-]_{\text{nor}})).$$
In other words:
 $$\U_{L^m(\omega^s)}([I^-]_{\text{nor}})=\sum_{n= 0}^{ m-1} \left(\frac{1}{d}\right)^n \U_{L(\omega^s)}([f^n(I^-)]_{\text{nor}}).$$
 
Recall again that $L(\omega^s)=L(\omega^{s-1})\wedge L(\omega)$ in the sense of current by Corollary 4.11 in \cite[Chapter III]{dem2}. So, in particular by Lemma \ref{dec_wedge}, a super-potential of $\U_{L(\omega^s)}$ is given by:
\begin{align*}
\U_{L(\omega^s)}(R)=\U_{L(\omega)}(L(\omega^{s-1})\wedge R)+\U_{L(\omega^{s-1})}(\omega \wedge R).
\end{align*}
on currents  $R\in \C_{k-s+1}$ such that $L(\omega^{s-1})$ and $R$ are wedgeable. 
A straightforward induction gives that a super-potential of $L(\omega^s)$ is given by:
\begin{align}\label{decomposition_L}
 \sum_{0\leq j\leq s-1} \U_{L(\omega)}( \omega^j\wedge L(\omega)^{s-1-j}\wedge R),
\end{align} 
on currents  $R\in \C_{k-s+1}$ such that $L(\omega^{s-1})$ and $R$ are wedgeable (since  $\omega^j\wedge L(\omega^{s-1-j})$ is more H-regular than $L(\omega^{s-1})$, we have that $\omega^j\wedge L(\omega)^{s-1-j}$ and $R$ are wedgeable by Lemma \ref{lemma_wedge_regular}).
In particular, $L(\omega^{s-1})$ and $\Lambda^n([I^-]_{\text{nor}})$ are wedgeable by Hypothesis \ref{distance} ((\ref{weak}) is enough for that since $L(\omega^{s-1})$ is smooth near $f^n(I^-)$) so we can take $R=\Lambda^n([I^-]_{\text{nor}})$ in the previous formula. \\ 

We want to show that for $0\leq j\leq s-1$, the following series is convergent:
$$\sum_{n= 0}^{m} \left(\frac{1}{d}\right)^n  \U_{L(\omega)}( \omega^j\wedge L(\omega)^{s-1-j}\wedge \Lambda^n([I^-]_{\text{nor}})). $$   
The term of the series can be rewritten as:
\begin{equation}\label{termseries}
a_{j,n}= \text{deg}(I^-)^{-1} \frac{1}{d^{sn}}\int_{f^n(I^-)} U_{L(\omega)}\omega^j\wedge L(\omega)^{s-1-j}, 
\end{equation}
since $\Lambda^n([I^-]_{\text{nor}})= [f^n(I^-)]_{\text{nor}}=\text{deg}(I^-)^{-1} d^{-(s-1)n} [f^n(I^-)]$ (observe that the form $U_{L(\omega)}  \omega^j\wedge L(\omega)^{s-1-j}$ is smooth on $f^{n}(I^-)$ so this integral makes sense). So in particular, Hypothesis \ref{distance} is equivalent to the fact that the series converges for $j=0$. We prove that the series converges for $j>0$ by induction.

So let $j>0$ be given such that the above series converges for $j-1$. Using $L(\omega)=dd^c U_{L(\omega)}+\omega$, we write: 
$$\omega^{j-1}\wedge L(\omega)^{s-j}=\omega^j\wedge L(\omega)^{s-1-j}+ dd^c U_{L(\omega)}\wedge \omega^{j-1} \wedge L(\omega)^{s-1-j}.$$   
 So replacing in (\ref{termseries}), we see that:
 $$a_{j-1,n}=a_{j,n}+ \text{deg}(I^-)^{-1} \frac{1}{d^{sn}}\int_{f^n(I^-)}  U_{L(\omega)} dd^c U_{L(\omega)}\wedge \omega^{j-1} \wedge L(\omega)^{s-1-j}. $$
 By Stokes, we have that the last integral is equal to:
 $$-\int_{f^n(I^-)}  dU_{L(\omega)} \wedge d^c U_{L(\omega)}\wedge \omega^{j-1} \wedge L(\omega)^{s-1-j}, $$ 
 which is non positive  because $\omega^{j-1} \wedge L(\omega)^{s-1-j}$ is positive. That means that:
 $$ a_{j-1,n}\leq a_{j,n}.$$
 Since $a_{j,n}<0$ (because $U_{L(\omega)}<0$), we have the convergence of the series for $j$. That concludes the induction.
 
By Proposition \ref{cor_decreasing_sqp}, we obtain the convergence in the Hartogs' sense to $T^+_s$. Furthermore, the convergence of the series giving $\U_{T^+_s}([I^-]_{\text{nor}})$ is indeed equivalent to the first half of Hypothesis \ref{distance}. \\

 Let $f$ satisfying Hypothesis \ref{distancefort} and let $q\leq s$. Then, we consider $R_{I^-}\in \C_{k-q+1}$ any positive closed current with support in $I^-$ (for example $\omega^{s-q}\wedge [I^-]_{\text{nor}}$). Then $\Lambda^j(R)$ is smooth near $I^+$ for all $j\leq m-1$, so we can apply (\ref{sum_green}):
$$\U_{L^m(\omega^q)}(R_{I^-})=\sum_{n= 0}^{ m-1} \left(\frac{1}{d}\right)^n \U_{L(\omega^q)}(\Lambda^{n}(R_{I^-})).$$
Using (\ref{decomposition_L}), we see that 
$$\U_{L(\omega^q)}( \Lambda^{n}(R_{I^-}))= \langle U_{L(\omega)},S_q\wedge \Lambda^{n}(R_{I^-}) \rangle,$$ 
where $S_q=  \sum_{0\leq j\leq q-1} \omega^j\wedge L(\omega)^{q-1-j}$ is smooth near $f^n(I^-)$ and is of mass $q$.
The measure $S_q\wedge \Lambda^{n}(R_{I^-})$ is of mass $q$ with support in $f^n(I^-)$. By Lemma \ref{estimate}, the function  $U_{L(\omega)}$ is greater than $A\log \text{dist}(I^+,f^n(I^-))-B$ on $f^n(I^-)$. Hypothesis \ref{distancefort} implies exactly the convergence of the series giving $\U_{L^m(\omega^q)}(R_{I^-})$. That concludes the proof by Proposition \ref{cor_decreasing_sqp}. \hfill $\Box$ \hfill \\

We will see in the next section how to construct the Green current of order $q$ ($q\leq s$) using only Hypothesis \ref{distance}. 
\begin{Remark}\label{otherpoints} \rm Using the same argument for  $f^n(I^-)$ instead of $I^-$ shows that the super-potentials of the current $T^+_s$ are in fact finite at every $[f^n(I^-)]_{\text{nor}}$. 

Observe also that the Green current $T^+_s(f^n)$ of $f^n$ is well defined and equal to $T^+_s$.
\end{Remark}
\begin{theorem}\label{invariance}
The current $T^+_s$ is $f^*$-invariant, that is $L(T_s^+)$ is well defined and equal to $T^+_s$. Furthermore, $T^+_s$ is the most $H$-regular current which is $f^*$-invariant in $\C_s$. In particular, $T_s^+$ is extremal in the set of $f^*$-invariant currents of $\mathcal{C}_s$.
\end{theorem}
\emph{Proof.} Recall from the appendix that a current $T$ is $f^*$-admissible if there exists a current $R_0$ which is smooth on a neighborhood of $I^+$ such that the super-potentials of $T$ are finite at $\Lambda(R_0)$. For such $T$, $f^*(T)$ is well defined and if $(T_n)$ is  a sequence of current converging in the Hartogs' sense to $T$ then $T_n$ is $f^*$-admissible and $f^*(T_n)$ converges to $f^*(T)$ in the Hartogs' sense.\\

In our case, we take for $R_0$ the current $[I^-]_{\text{nor}}$ which is smooth near $I^+$. Then by Remark \ref{otherpoints}, the super-potentials of $T^+_s$ are finite at $\Lambda([I^-]_{\text{nor}})= [f(I^-)]_{\text{nor}}$. So the current $L(T_s^+)$ is well defined. Now, we have that $(L^{n+1}(\omega^s))_n=(L(L^n(\omega^s)))_n$ converges in the Hartogs' sense to $T^+_s$ and to $L(T^+_s)$ so that $T^+_s$ is indeed invariant. \\

Let $\U_{T^+_s}$ be the super-potential of $T^+_s$ defined as:
\begin{equation}\label{correctmean}
\U_{T^+_s}= \sum_{n= 0}^{\infty} \left(\frac{1}{d}\right)^n \U_{L(\omega^s)}\circ \Lambda^{n},
\end{equation}
on smooth forms in $\C_{k-s+1}$. Then composing (\ref{correctmean}) with $\Lambda$, we have that on smooth forms in $\C_{k-s+1}$:
$$\U_{T^+_s}=d^{-1}\U_{T^+_s}\circ \Lambda +\U_{L(\omega^s)}.$$ 
By iteration, we have that on smooth forms in $\C_{k-s+1}$:
$$\U_{T^+_s}= \sum_{n= 0}^{m-1} \left(\frac{1}{d}\right)^n \U_{L(\omega^s)}\circ \Lambda^{n}  + \left(\frac{1}{d}\right)^m \U_{T^+_s}\circ \Lambda^{m}. $$ 
By Theorem \ref{convergencecurrent}, that implies by difference that:
 $$\left(\frac{1}{d}\right)^m \U_{T^+_s}\circ \Lambda^{m}$$
goes to zero on smooth forms in $\C_{k-s+1}$.\\

Now, let $S$ be a $f^*$-invariant current in $\mathcal{C}_s$ such that there are constants $A>0$ and $B$ satisfying $A\U_{T_s^+} +B \leq \U_S<0$ for some super-potentials $\U_S$ and $\U_{T_s^+}$.  Then on smooth forms in $\C_{k-s+1}$, a super-potential $\U_{L^m(S)}$ of $L^m(S)=S$ is given by:
$$ \sum_{n= 0}^{m-1} \left(\frac{1}{d}\right)^n \U_{L(\omega^s)}\circ \Lambda^{n}  + \left(\frac{1}{d}\right)^m \U_{S}\circ \Lambda^{m}.$$ 
Since $\left(\frac{1}{d}\right)^m \U_{T_s^+}\circ \Lambda^{m}$ goes to zero on smooth forms in $\C_{k-s+1}$, our hypothesis implies that $\left(\frac{1}{d}\right)^m \U_{S}\circ \Lambda^{m}$ also goes to zero on smooth forms in $\C_{k-s+1}$. In particular, a super-potential of $S$ is given on smooth forms  in $\C_{k-s+1}$ by:
$$ \sum^{\infty}_{n=0}\left(\frac{1}{d}\right)^n \U_{L(\omega^s)}\circ \Lambda^{n}. $$
Now, using the fact that two currents having the same super-potential on smooth forms are in fact equal we deduce that $T^+_s=S$. 

In particular, for $A=1$, we obtain that $T^+_s$ is the most $H$-regular current which is $f^*$-invariant. It is extremal in the set of $f^*$-invariant currents of $\mathcal{C}_s$ since if not we could write $T^+_s=tS_1+(1-t)S_2$ with $S_1$ and $S_2$ two $f^*$-invariant currents in $\mathcal{C}_s$. Take $M$ small enough so that the super-potentials $\U_1$, $\U_2$ and $\U_{T^+_s}$ of $S_1$, $S_2$ and $T^+_s$ of same mean $M$ are negative. Observe that $\U_{T^+_s}=t\U_1+(1-t)\U_2$ so that $t^{-1}\U_{T^+_s}\leq \U_1$. Then we can apply the previous result for $A=t^{-1}$ and it follows that $S_1=T$ (similarly $S_2=T$). \hfill $\Box$ \hfill \\
 
In the previous proof, we have obtained:
\begin{corollaire}\label{greennorma}
Let $\U_{T^+_s}$ be the super-potential of $T^+_s$ defined on smooth forms by:
$$\U_{T^+_s}= \sum_{n= 0}^{\infty} \left(\frac{1}{d}\right)^n \U_{L(\omega^s)}\circ \Lambda^{n}.$$
Then we have that:
$$\left(\frac{1}{d}\right)^m \U_{T^+_s}\circ \Lambda^{m}$$
goes to zero on smooth forms.
\end{corollaire} 
The current $L(\omega)^{s+1}$ is a well defined element of $\C_{s+1}$ and we have that 
$$dd^c U_{L(\omega)} \wedge L(\omega)^{s}+\omega \wedge L(\omega)^{s}= L(\omega)^{s+1}$$ 
by Corollary 4.11 Chapter III in \cite{dem2}. Observe that it is not true though that $L(\omega^{s+1})=L(\omega)^{s+1}$.
Indeed, $f^*(\omega^{s+1})$ is a well defined form (of mass $\delta^{k-s-1}$) which is in $L^1$ hence that does not give mass to algebraic sets of dimension $\leq k-1$. But $f^*(\omega)^{s+1}$ is a smooth form outside $I^+$ which coincides with $f^*(\omega^{s+1})$ there. So we have by Siu's Theorem that:
$$f^*(\omega)^{s+1}=\sum_{i} a_i [I^+_i] +f^*(\omega^{s+1}) $$ 
where the sum goes over all the irreducible components $I^+_i$ of $I^+$ and where the $a_i$ are positive numbers. Observe that 
this formula is related to King's formula (see \cite{dem2} Chapter III proposition 8.18). In particular, one has that:
$$ f^*(\omega)^{s+1} \leq  C[I^+] +f^*(\omega^{s+1}) $$  
for $C>0$ large enough. 
 Similarly, one has that $f_*(\omega)^{k-s+1}$ is well defined and satisfies:
 $$f_*(\omega)^{k-s+1} \leq  C[I^-] +f_*(\omega^{k-s+1}).$$
The following proposition is useful in the construction of the equilibrium measure.
\begin{proposition}\label{Hregular}
The super-potentials of $T^+_s$ are finite at $\omega^j\wedge\Lambda(\omega)^{k-s+1-j}$ for all $k-s+1\geq j\geq 0$. 
 \end{proposition}
\emph{Proof.} If two currents $S_1$ and $S_2$ in $\C_r$ satisfies $S_1\leq c S_2$ for a constant $c>0$ then the super-potentials of $S_1$ and $S_2$ of mean $0$ satifies $\U_{S_1}\geq c\U_{S_2}+c'$ for some constant $c'$.  In particular, the super-potentials of $S_1$ are finite wherever $\U_{S_2}$ is. 

Now, we have that the super-potentials of $T^+_s$ are finite at $[I^-]_{\text{nor}}$.  Since $T^+_s$ is $f^*$-admissible, its super-potential are finite at every point of the form $\Lambda(S)$ for $S\in \C_{k-s+1}$ smooth near $I^+$. So they are also finite at $\Lambda(\omega^{k-s+1})$. The affinity of the super-potentials of $T^+_s$ implies that they are finite at $ (C[I^-] +f_*(\omega^{k-s+1}))_\text{nor}$. So the super-potential of $T^+_s$ are finite at $(f_*(\omega)^{k-s+1})_\text{nor}$. Since for $j\geq 0$, the current $\omega^j\wedge \Lambda(\omega)^{k-s+1-j}$ is more H-regular than $\Lambda(\omega)^{k-s+1}$, we have that the super-potentials of $T^+_s$ are finite at $\omega^j\wedge \Lambda(\omega)^{k-s+1-j}$ (we use the symmetry of the super-potential : $\U_{T^+_s}(\omega^j\wedge \Lambda(\omega)^{k-s+1-j})=\U_{\omega^j\wedge \Lambda(\omega)^{k-s+1-j}}(T^+_s)$).  \hfill $\Box$ \hfill 
\begin{corollaire}
The current $T^+_s$ gives no mass to $I^-$ (nor $I^+$ by dimension's argument). 
\end{corollaire}
\emph{Proof.} From above, the super-potentials of $T^+_s$ are finite at $\Lambda(\omega)\wedge\omega^{k-s}\in \C_{k-s+1}$. Observe that for two currents $R$ and $S$ in $\C_r$ and $\C_s$ with $r+s \leq k$, then:
$$\U_{R}(S\wedge \omega^{k+1-r-s})=\U_{R\wedge \omega^{k+1-r-s}}(S)= \U_{S}(R\wedge \omega^{k+1-r-s})$$
where all the super-potentials are of same mean. 

So for $R=T^+_s$ and $S=\Lambda(\omega)$, we get that the super-potentials of $\Lambda(\omega)$ are finite at $T_s^+\wedge \omega^{k-s}$. This is equivalent to the fact that the quasi-potential  $U_{\Lambda(\omega)}$ is integrable with respect to the measure  $T_s^+\wedge \omega^{k-s}$. In other words:
$$\int U_{\Lambda(\omega)} \omega^{k-s}\wedge T^+_s$$
is finite. Applying Lemma \ref{estimate} to $f^{-1}$ shows that the singularities of $U_{\Lambda(\omega)}$ are in $\log \text{dist}(x,I^-)$ so $T^+_s$ gives no mass to $I^-$.  \hfill $\Box$ \hfill 

\begin{Remark}\label{Lelong} \rm The quantity $\U_{T^+_s}([I^-]_{\text{nor}})$
is related to a generalized Lelong number (\cite{dem3}). Let us explain this point.  
From the previous proposition, we have that the super-potentials of $T^+_s$ are finite at $\Lambda(\omega)^{k-s+1}$.  

We define the Lelong number of $T^+_s$ associated to the function $U_{\Lambda(\omega)}$ as:
$$\nu(T^+_s, U_{\Lambda(\omega)})=\lim_{r\to -\infty} \int_{\{U_{\Lambda(\omega)}<r\}} T^+_s \wedge \Lambda(\omega)^{k-s}.$$
The previous current is well defined by the theory of super-potentials: the super-potentials of $T^+_s$ are finite at $\omega\wedge \Lambda(\omega)^{k-s}$ which means that $T^+_s$ and $\Lambda(\omega)^{k-s}$ are wedgeable and their wedge product is well defined by Definition \ref{def_wedge}.

As in formula (\ref{decomposition_L}), we have that a super-potential of $\Lambda(\omega)^{k-s+1}$ is given by:
\begin{align*}
 \sum_{0\leq j\leq k-s} \U_{\Lambda(\omega)}( \omega^j\wedge \Lambda(\omega)^{k-s-j}\wedge R),
\end{align*} 
on currents  $R\in \C_{s}$ such that $\Lambda(\omega)^{k-s}$ and $R$ are wedgeable, so in particular for $R=T^+_s$ by the previous proposition. Since the super-potentials of $T^+_s$ are finite at $\Lambda(\omega)^{k-s+1}$, that implies that every term in the previous sum is finite. So we have in particular that:
$$\U_{\Lambda(\omega)}(  \Lambda(\omega)^{k-s}\wedge T^+_s) $$  
is finite. That means that the quasi-potential $U_{\Lambda(\omega)}$ is integrable with respect to the measure $\Lambda(\omega)^{k-s}\wedge T^+_s$. 
So we can use the bound:
\begin{align*}
\int_{\{U_{\Lambda(\omega)}<r\}} T^+_s \wedge \Lambda(\omega)^{k-s}  &\leq \frac{1}{r} \int_{\{U_{\Lambda(\omega)}< r\}} U_{\Lambda(\omega)} T^+_s \wedge \Lambda(\omega)^{k-s}\\
																																 &\leq \frac{1}{ r} \int_{\P^k} U_{\Lambda(\omega)}  T^+_s \wedge\Lambda(\omega)^{k-s}.
\end{align*}
So we have that:
$$\nu(T^+_s,U_{\Lambda(\omega)})=0.$$
This is a generalization of the fact 
that a psh function finite at the point $x$ has zero Lelong number at $x$. 
\end{Remark}

A classical question in complex dynamics is to ask by what $\omega^s$ can be replaced.
In other words, what are the currents $T$ in $\C_s$ such that $L^m(T)\to T^+_s$? 
The following proposition and theorem give partial results toward this direction. 
\begin{proposition}
Let $(T_m) \in \C_s$ be a sequence of currents such that 
a super-potential $\U_{T_m}$ of $T_m$ satisfies 
$\| \U_{T_m} \|_{\infty} = o(d^{m})$. 
Then $L^m(T_m)\to T^+_s$ in the Hartogs' sense.
\end{proposition} 
\emph{Proof.} Recall first that if $T\in \C_s$ has bounded super-potential it is $(f^n)^*$-admissible (its super-potential is in particular bounded at the point $\Lambda^n(\omega^{k-s+1})$). So the sequence of pull-back is well defined by Lemma \ref{stability}. Using Proposition \ref{prop_pull_general} and (\ref{sum_green}), a super-potential of $L^m(T_m)$ is given on smooth currents in $\C_{k-s+1}$ by:
$$ \U_{L^m(T_m)} = \sum_{n= 0}^{m-1} \left(\frac{1}{d}\right)^n \U_{L(\omega^s)}\circ \Lambda^{n}  + \left(\frac{1}{d}\right)^m \U_{T_m}\circ \Lambda^{m}. $$ 
By Theorem \ref{convergencecurrent} we know that the series $\sum_{n= 0}^{m-1} \left(\frac{1}{d}\right)^n \U_{L(\omega^s)}\circ \Lambda^{n}$ converges in the Hartogs' sense to $\U_{T^+_s}$. The hypothesis on $(T_m)$ implies that $\left(\frac{1}{d}\right)^m \U_{T_m}\circ \Lambda^{m}=o(1)$ goes to 0 uniformly on smooth form. Since the control is uniform we have that ($|\U_{L^m(T_m)}- \U_{L^m(\omega^s)}|\to 0$), and the convergence is in the Hartogs' sense and we can conclude by Proposition \ref{lemma_cv_sp_pointwise}. \hfill $\Box$ \hfill \\

We also have the following result. We believe the proof is of interest although the result is essentially already known. We refer the reader to \cite{Sib} for the case $q=1$ and also \cite{DS2} for the general case. See the Appendix for the notion of super-polarity.
\begin{theorem}
There exists a super-polar set $P$ of $\C_s$ such that for $S \in \C_s \backslash P$, $L^m(S)$ is well defined and converges to $T^+_s$. 
\end{theorem}
\emph{Proof.} The set of currents $S\in \C_s$ such that $S$ is not $(f^m)^*$-admissible is super-polar since it is contained in the set of currents such that $\U_S(\Lambda^m(\omega^{k-s+1}))= -\infty$. Now, since a countable union of super-polar set is super-polar, we have that outside a super-polar set of $\C_s$, $S$ is $(f^m)^*$-admissible and so $L^m(S)$ is well defined by Lemma \ref{stability}.
 
As above, a super-potential $\U_{L^m(S)}$ of $L^m(S)$ is given on smooth forms by:
$$ \sum_{n=0}^{m-1} \left(\frac{1}{d}\right)^n \U_{L(\omega^s)}\circ \Lambda^{n}  + \left(\frac{1}{d}\right)^m \U_{S}\circ \Lambda^{m}, $$
where $\U_S$ is a super-potential of $S$. For $\Omega \in \C_{k-s+1}$ smooth, consider the current $R(\Omega)\in \C_{k-s+1}$ defined by $R(\Omega):= (\sum_m \left(\frac{1}{d}\right)^m \Lambda^{m}(\Omega ))_{\text{nor}}$. Let $P$ be the set of currents $S$ in $\C_s$ such that $\U_S(R(\omega^{k-s+1}))=-\infty$, then $P$ is super-polar by definition. Observe that for $\Omega \in \C_{k-s+1}$ smooth, we have a constant $c_\Omega>0$ such that $R(\Omega) \leq c_\Omega R(\omega^{k-s+1})$. In particular, for $S \notin P$, we have that $\U_S(R(\Omega))>-\infty$.

That implies that for any $\Omega$ smooth and $S \notin P$, the sequence $\U_{L^m(S)}(\Omega)$ converges to the value $\U_{T^+_s}(\Omega)$. 
Indeed, the fact that $\U_S(R(\Omega))$ is finite gives that:
$$\left(\frac{1}{d}\right)^m \U_{S}\circ \Lambda^{m}(\Omega) $$
goes to $0$. So $\U_{L^m(S)}(\Omega)$ converges to $\U_{T^+_s}(\Omega)$.
Then Proposition \ref{lemma_cv_sp_pointwise} gives us that the sequence $L^m(S)$ converges in fact to $T^+_s$ (maybe not in the Hartogs' sense).   \hfill $\Box$ \hfill \\

Of course, the above theorem is not optimal and it is conjectured that for $T$ the current of integration on a (very) generic algebraic set of dimension $k-s$, $L^m(T)$ goes to $T^+_s$ (see in the endomorphisms case \cite{BDu} and \cite{DS3} for the case of measures  and see \cite{DS7}  for the case of bidegree $(1,1)$). That is beyond the scope of this study. \\

Recall that we consider the critical sets $\mathcal{C}^+$ (or $\C(f)$) and $\mathcal{C}^-$ (or $\C(f^{-1})$) defined by:
\begin{align*}
 \mathcal{C}^+&:= f^{-1}(I^-) \\
 \mathcal{C}^-&:=(f^{-1})^{-1}(I^+).
 \end{align*}
We have the proposition:
\begin{proposition}\label{pushforwardgreen}
The current $T^+_s$ does not give mass to the critical sets $\C^-$ and  $\C^+$. In particular, the current $T^+_s$ satisfies the equation:
$$f_*(T^+_s)=\frac{1}{d_s} T^+_s,$$ 
in $\P^k \setminus \C^-$.
\end{proposition} 
We first need the following lemma:
\begin{lemme}
Let $\varphi$ be a smooth function. Then $f_*(\varphi)$ is in $L^1(T^+_s\wedge \Lambda(\omega^{k-s}))$ and we have the identity:
$$\int  \varphi T^+_s\wedge \omega^{k-s}     =    \int f_*(\varphi) T^+_s\wedge \Lambda(\omega^{k-s})    $$
\end{lemme}
\emph{Proof of the lemma.} First $\varphi$ is in $L^1(T^+_s\wedge \omega^{k-s})$ and the quantity:
$$\int  \varphi T^+_s\wedge \omega^{k-s} $$
depends continuously on $ T^+_s$ in the sense of currents.

On the other hand, $f_*(\varphi)$ is in $L^1(T^+_s\wedge \Lambda(\omega^{k-s}))$ since it is smooth and bounded outside $I^-$ which has no mass for $T^+_s\wedge \Lambda(\omega^{k-s})$ (see Remark \ref{Lelong}).
Recall that $T^+_s \wedge \omega^{k-s}$ is $f^*$-admissible by Remark \ref{Lelong} (we proved that $\U_{\Lambda(\omega)}(  \Lambda(\omega)^{k-s}\wedge T^+_s) $  
is finite). So we have that 
$$ \int f_*(\varphi) T^+_s\wedge \Lambda(\omega^{k-s})= \int \varphi   L(T^+_s\wedge \Lambda(\omega^{k-s})),$$
as this stands if $T^+_s\wedge \Lambda(\omega^{k-s})$ was smooth and we can conclude by Hartogs' convergence. Now, applying Lemma \ref{pushpull} to $f^{-1}$ and the invariance of $T^+_s$, we have that $L(T^+_s\wedge \Lambda(\omega^{k-s}))= T^+_s\wedge \omega^{k-s}$.  \hfill $\Box$ \hfill \\

\noindent \emph{Proof of the proposition.} We consider $\C^+$ first. Let $V_\varepsilon$ be a small neighborhood of $I^+$. Since $T^+_s$ gives no mass to $I^+$ we can assume that the mass of $V_\varepsilon$ for $T^+_s$ is arbitrarily small. Let $W_\alpha$ be a small neighborhood of $\C^+$. 
Let $0\leq \varphi \leq 1$ be a smooth function such that $\varphi=1$ in $W_\alpha \backslash V_\varepsilon$, $\varphi=0$ in $ V_{2^{-1}\varepsilon}$ and $\varphi=0$ outside $W_{2\alpha}$.

Then by the previous lemma:
 \begin{align*}
\|W_\alpha \backslash V_\varepsilon \|_{T^+_s}& \leq \int \varphi T^+_s \wedge \omega^{k-s}\\
                                              & \leq \int f_*(\varphi) T^+_s\wedge \Lambda(\omega^{k-s})\\
                                              & \leq \int_{f(W_{2\alpha} \backslash V_{2^{-1}\varepsilon})} T^+_s\wedge \Lambda(\omega^{k-s})\\
\end{align*}
Now, $f(W_{2\alpha} \backslash V_{2^{-1}\varepsilon})$ is a neighborhood $W$ of $I^-$ that can be taken arbitrarily small by taking $\varepsilon$ and $\alpha$ small enough. We have seen in Remark \ref{Lelong} that the quantity $\int_{W} T^+_s \wedge \Lambda(\omega)^{k-s}$ goes to the Lelong number $\nu(T^+_s,U_{\Lambda(\omega)})$ which is equal to zero.  Thus $\C^+$ has no mass for $T^+_s$. \\

 For $\C^-$, we write $T^+_s=T_{\C^-}+T'$ where $T_{\C^-}$ is a positive closed current with support in $\C^-$ and $T'$ is a positive closed current with no mass on $\C^-$ (\cite{Sko1}). Both currents are $f^*$-admissible since $T^+_s$ is. But $f^*(T_{\C^-})$ has support in $\C^+$ which means it is equal to zero since $T_s^+=d_s^{-1}(f^*(T_{\C^-})+f^*(T'))$ gives no mass to $\C^+$. That implies that $T_{\C^-}=0$ hence $T^+_s$ has no mass on $\C^-$.\\

 Now $f:\P^k  \setminus   \mathcal{C}^+ \rightarrow \P^k \setminus  \mathcal{C}^-$ is a proper biholomorphism that we will denote by $f_1$. If $\Theta$ is a smooth form in $\P^k \setminus  \mathcal{C}^-$ then using the invariance of $T^+_s$:
 \begin{align*}
 \langle (f_1)_*(T^+_s)       ,   \Theta          \rangle &=  \langle T^+_s ,  (f_1)^*( \Theta)          \rangle \\
                                                          &=  \langle  \frac{1}{d^s}f^*(T^+_s) ,  (f_1)^*( \Theta)          \rangle \\
                                                          &=  \langle  \frac{1}{d^s}T^+_s ,  f_*((f_1)^*( \Theta))         \rangle.
 \end{align*} 
 But $f_*(f_1)^*=(f_1)_*(f_1)^*=\text{id}$ so $f_*((f_1)^*( \Theta))=\Theta$ and the result follows.   \hfill $\Box$ \hfill 

\begin{Remark} \rm In order to define $\Lambda(T^+_s)$, we need to add to the equation 
$$(f_1)_*(T^+_s)= \frac{1}{d_s} T^+_s$$
a term of mass $d_{k-s}-d_s^{-1}$ and support in $\C^-$ in order to obtain a current of mass $d_{k-s}$. For example, in the case of Hénon maps, we need to add a multiple of the current of integration on the line at infinity. In general, such choice is not clear and they might be no way to add a current to the equation  $(f_1)_*(T^+_s)= \frac{1}{d_s} T^+_s$ in a continuous way.
\end{Remark} 

The previous corollary implies that the Green current is meaningful. For example, in the case of Hénon maps in $\mathbb{P}^2$, the set $\C^-$ and $\C^+$ are in fact $L_\infty$ (the line at infinity) which is totally invariant and the Green current $T^+$ gives no mass to $L_\infty$.\\

We can now prove the following stronger result of extremality which implies strong ergodic properties (see Theorem \ref{mixing}). The inequalities between currents in $\C_s$ have to be understood in the strong sense.
\begin{theorem}\label{extremality}
The current $T^+_s$ is extremal in $\C_s$, that is if there exists a $c>0$ and $S\in \C_s$ such that $S\leq cT^+_s$ then $S= T^+_s$. 
\end{theorem} 
\emph{Proof.} Applying the previous results to $f^n$ gives that $T^+_s$ gives no mass to the indeterminacy sets $I(f^{\pm n})$ and critical sets $\C^{\pm n}$ of $f^n$ and $f^{-n}$ (recall that $T_s^+$ is also the Green current of $f^n$). Let $S$ be as above, in particular $S$ gives no mass to the sets $I(f^{\pm n})$ and $\C(f^{\pm n})$ and $S$ is $(f^n)^*$-admissible for all $n$. By Lemma \ref{stability}, $L^n(S)$ is well defined for all $n$ and equal to $L_n(S)$ ($L_n$ and $\Lambda_n$ are the normalized pull-pack and push-forward associated to $f^n$).

For $n>0$, we denote by $\Lambda_n'$ the push forward operator $(f^n)_*$ from $\mathbb{P}^k \C(f^n)$ to $\mathbb{P}^k\setminus \C(f^{-n})$. 

The operator $(\Lambda_n')$ is positive. That and the previous proposition applied to $f^n$ imply that $(d_s)^n(\Lambda_n')(S) \leq c T_s^+$ in $\mathbb{P}^k\setminus \C(f^{-n})$. We denote by $S_n$ the trivial extension of $(d_s)^n(\Lambda_n')(S)$ over $\mathbb{P}^k$ : it does exist since the current $(d_s)^n(\Lambda_n')(S)$ is of bounded mass. We have $S_n \leq c T_s^+$  in $\P^k$. In particular, $S_n$ is $(f^*)^n$-admissible.
Using the same argument as in the previous proposition, we see that:
$$ (f^n)^*(S_n)=  (d_s)^n S,$$
outside $\C(f^n)$. Now $S$ has no mass on $\C(f^n)$ and  $(f^n)^*(S_n)$ is less than $c(d^s)^nT_s^+$ (by positivity of the operator $(f^n)^*$) which implies that $(f^n)^*(S_n)$ also has no mass on $I(f^n)\cup\C(f^n)$. So we have:
$$ (f^n)^*(S_n)=  (d_s)^n S,$$
on $\P^k$.
In particular, $S_n$ has mass 1. We have that $L_n(S_n)=S$ and  since $S_n$ is $(f^*)^n$-admissible we have $L_n(S_n)=L^n(S_n)$ by Lemma \ref{stability}.\\

For $R\in C^{k-s+1}$ smooth, let $\U_{S_n}$, $\U_{T^+_s,0}$ and $\U_{\Lambda^j(R)}$ be the super-potential of $S_n$, $T^+_s$ and $\Lambda^j(R)$ of mean $0$. We have from Proposition \ref{prop_pull_general} and (\ref{sum_green}) that a super-potential $\U_{L^n(S_n)}$ of $L^n(S_n)=S$ on smooth forms is given by:
$$\sum_{j=0}^{n-1} \left(\frac{1}{d}\right)^j\U_{L(\omega^s)} \circ \Lambda^j+ \left(\frac{1}{d}\right)^n\U_{S_n} \circ \Lambda^n.$$
 Recall that there exists a $M>0$ that does not depend on $R$ and $n$ such that $\U_{\Lambda^n(R)}-M$ is negative and $\U_{S_n}\leq M$. More precisely, by Proposition \ref{resolutionquasipotential}, there exists a quasi-potential $U_{\Lambda^n(R)}$ of $\Lambda^n(R)$ such that $U_{\Lambda^n(R)}-M\omega^{k-s}$ is negative ($U_{\Lambda^n(R)}$ is a quasi-potential of $\Lambda^n(R)$ of mean $0$). Then, we have that:
 $$(\U_{\Lambda^n(R)}-M)(S_n) \geq  c(\U_{\Lambda^n(R)}-M)(T^+_s).$$
 Indeed if $S_n$ and $T^+_s$ were smooth, it would follow from the fact that $U_{\Lambda^n(R)}-M\omega^{k-s}$ is negative and that $S_n\leq c T^+_s$. The result follows then by Hartogs' convergence: observe that the regularization is obtained by a mean of the composition over the automorphisms of $\P^k$ thus the approximations $S'_n$ and $T'^+_s$ of $S_n$ and $T^+_s$ also satisfy  $S'_n\leq c T'^+_s$. 

So we have the estimates:
\begin{align*}
\U_{S_n}(\Lambda^n(R))&= \U_{\Lambda^n(R)}(S_n)\\
                      &=(\U_{\Lambda^n(R)}-M)(S_n)+\langle S_n,M \omega^{k-s} \rangle\\
                      &\geq c(\U_{\Lambda^n(R)}-M)(T^+_s)+\langle S_n,M \omega^{k-s} \rangle\\
                      &\geq c\U_{T^+_s,0}(\Lambda^n(R))+\langle S_n-c T^+_s, M \omega^{k-s} \rangle \\
                      &\geq c\U_{T^+_s,0}(\Lambda^n(R))+M(1-c).
\end{align*}
 The last term multiplied by $d^{-n}$ goes to zero by Corollary \ref{greennorma}. So the super-potential $\U_{L^n(S_n)}$ converges to a super-potential of $T^+_s$ on smooth forms. By Proposition \ref{lemma_cv_sp_pointwise}, that implies that $S=T^+_s$. \hfill $\Box$ \hfill \\
 
 Now we prove that the mappings satisfying Hypothesis \ref{distance} are generic. Our statement is similar to Proposition 4.5 in \cite{BD1}. In addition, we show here that the genericity stands in any orbit. The idea of the proof is to construct an element in any orbit satisfying Hypothesis \ref{distance} and then to show that the series giving Hypothesis \ref{distance} varies as a difference of psh functions (dsh) along the orbit.
\begin{theorem}\label{generic}
Let $E_s$ be the set of birational maps $f:\P^k\to \P^k$ such that $I^+$ and $I^-$ satisfies $\text{dim}(I^+)=k-s-1$ and $\text{dim}(I^-)=s-1$. Consider the group action:
\begin{align*}
 \Phi: \text{PGL}(k+1,\mathbb{C})\times\text{PGL}(k+1,\mathbb{C})\times E_s &\to E_s \\
                                                                    (A,B,f) &\mapsto B\circ f \circ A^{-1}.
\end{align*}
Then outside a pluripolar set of the orbit $Orb(f)$ of $f\in E_s$, the maps of $Orb(f)$ satisfy Hypothesis \ref{distance}.
 \end{theorem}
\emph{Proof.} We change the definition of $\Phi$ and we define $\Phi(A,B,f)=B\circ f \circ A$ instead. This is not a group action but it is sufficient to prove the statement for such $\Phi$ since taking the inverse is a biholomorphism on $\text{PGL}(k+1,\mathbb{C})$. We define $Orb'(f)=\{B\circ f \circ A\}$ and we still speak of the orbit of $f$. Observe first that the algebraic degree of the elements of $Orb'(f)$ is constant equal to $d$.  \\

\noindent{\bf Construction of an example stable by perturbations satisfying the first
half of Hypothesis \ref{distance}.}

 We have that:
$$I(B\circ f \circ A)=A^{-1}I(f) $$
and
$$I(A^{-1}\circ f^{-1} \circ B^{-1})=BI(f^{-1}) .$$
In particular, for $A$ and $B$ generic, we can assume that $I^+ \cap I^-=\varnothing$. Remark also that the dimension of these indeterminacy sets is constant on the orbit of $f$. \\

Consider a projective linear subspace $E$ of $\P^k$ of dimension $s$ such that $E\cap I(f)=\varnothing$. Let $V$ be a neighborhood of $E$ such that $\overline{V}\cap I(f)=\varnothing$.

Choose the coordinates $[z_0: \dots:z_k]$ in $\P^k$ such that $E$ is given by $z_0=\dots=z_{k-s-1}=0$ and so that if $E'=\{z_{k-s}= \dots = z_k=0\}$ then  $E' \cap I(f^{-1})=\varnothing$, $E'\cap \overline{V}=\varnothing $ and $E'\cap \overline{f(V)}=\varnothing $ (this is always possible).
Consider the element $B_0$ of $\text{PGL}(k+1,\mathbb{C})$ defined by $B_0([z_0: \dots:z_k])=[\lambda z_0: \dots: \lambda z_{k-s-1}: z_{k-s}: \dots :z_k]$ with $\lambda >0$. We take $\lambda$ small enough so that:
\begin{itemize} 
\item $B_0(I(f^{-1})) \subset V$;
\item $B_0(f(V)) \subset V$.
\end{itemize}
Consider the element $f_{B_0}$ of $Orb'(f)$ defined by $f_{B_0}={B_0}\circ f$. Now, $I(f_{B_0})=I(f)$ and $I((f_{B_0})^{-1})={B_0}(I(f^{-1})) \subset V$. We have the inclusion:
\begin{align*}
f_{B_0}(I((f_{B_0})^{-1}))=f_{B_0}{B_0}(I(f^{-1}))\subset ({B_0}\circ f)(V) \subset V.
\end{align*}
An immediate induction gives that $(f_{B_0})^n(I((f_{B_0})^{-1})) \subset V$. So the element $f_{B_0}$ satisfies the first half of Hypothesis \ref{distancefort} since 
$$\text{dist}(I(f_{B_0}),f_{B_0}^n(I(f_{B_0}^{-1})))\geq  \text{dist}(I(f),V)>0.$$ 

For each $n, m \in \mathbb{N}$, the condition $f^n(I(f^{-1})) \cap f^{-m}(I(f)) \neq \varnothing$ is algebraic (and not always satisfied by the above), so (\ref{weak}) is satisfied outside a countable union of subvarieties of $\text{PGL}(k+1,\mathbb{C})^2$. Wherever all these conditions are satisfied, namely wherever (\ref{weak}) holds the dynamical degrees are given by Proposition \ref{s-stability} and are thus constant. \\   

Now, we show that we can find a small open set $W'_0$ in $Orb'(f)$
where the first part of Hypothesis \ref{distance} is satisfied. Fix $E$ and $\overline{V}$ as above. If $\Phi(f)=B\circ f \circ A$ with $A$ close to the identity and $B$ closed to $B_0$, then we still have $I(\Phi(f)^{-1})= BI(f^{-1}) \subset V$ and $B\circ f \circ A(V)\subset V$ since $A(V)$ is close to $V$. Thus $\Phi(f)^n(I(\Phi(f)^{-1}))\subset V$ and $I(\Phi(f))=A^{-1}I(f)$ is close to $I(f)$. 

This implies that there exists some $\alpha>0$ such that for every
$(A,B)$ in a small neighborhood $W'_0$ of $(Id, B_0)$ and every $n\in \mathbb{N}$, we have:
$$\text{dist}(I(\Phi(f)),\Phi(f)^n(I(\Phi(f)^{-1})))\geq \alpha. $$  

Now we prove the genericity. In what follows, $C$ denotes a constant independent of $n$ and $N$ that may change from line to line. 
Let $W := \text{PGL}(k+1,\mathbb{C})^2$. It is a Zariski dense open
set in the projective space $\widetilde{W}=\mathbb{P}^N \times
\mathbb{P}^N$. Let $l$ be the complex dimension of $W$
($l=2(k+1)^2-2$). Let $c$ denote the homogeneous coordinate on
$\widetilde{W}$. When $c \in W$, we write $f_c$ instead of
$\Phi(c,f)$. We can extend this notation for $c \in \widetilde{W}$. Of
course, in this case $f_c$ is not a birational map.  

Consider the rational map:
\begin{align*}
\widetilde{F}_n: \widetilde{W} \times \P^k &\to \widetilde{W} \times \P^k\\
(c,z) &\mapsto  (c,f^n_c(z)).
\end{align*}
Let $\Pi_i$ denote the canonical projections of $\widetilde{W}\times
\P^k$ to its factor for $i=1,2,3$ (recall that $\widetilde{W}=\P^N
\times \P^N$) and, in $\widetilde{W}\times \P^k$, let $\omega_i:=\Pi_i^*(\omega_{FS})$ be the pull-back of the Fubini-Study form by the projection for $i=1,2,3$. That way, $\omega_1+\omega_2+ \omega_3$ is a Kähler form on $\widetilde{W}\times \P^k$.\\

\noindent{\bf Action of $\widetilde{F}_n^*$ on the cohomology.}

We study the action of $\widetilde{F}_n^*$ on $\omega_3$. Write $c=(c_1,c_2)=([c_{1,1}:\dots:c_{1,N+1}],[c_{2,1}:\dots:c_{2,N+1}])$. First we have $\widetilde{F}_n(c,z)=(c,f^n_c(z))$ where the second coordinate is a polynomial of degree $d^n$ in the $z_i$, of degree $\leq C d^{n}$ in the $c_{1,i}$ and in the $c_{2,i}$. We compute the mass of $\widetilde{F}_n^*(\omega_3)$. For that, we test against $(\omega_1+\omega_2+ \omega_3)^{k+l-1}$. Write $\Omega:=\omega_1+\omega_2$. We developp $(\omega_1+\omega_2+ \omega_3)^{k+l-1}$:
$$(\omega_1+\omega_2+ \omega_3)^{k+l-1}= \sum_{i=0}^{k+l-1}\binom{k+l-1}{i} \Omega^i \wedge \omega_3^{k+l-1-i}.$$
We have that $\Omega^i=0$ for $i> l$ and $\omega_3^{k+l-1-i}=0$ for $k+l-1-i>k$ thus $i<l-1$. So there are only two terms in the previous sum: for $i=l-1$ and for $i=l$.

There are two terms to control:
\begin{align*}
\langle  \widetilde{F}_n^* (\omega_3), \Omega^{l-1} \wedge \omega_3^k  \rangle  \quad \text{and} \quad \langle  \widetilde{F}_n^* (\omega_3), \Omega^{l} \wedge \omega_3^{k-1}  \rangle. 
\end{align*}
By Bézout's theorem, those two terms are $\leq C d^n$ (the terms can
be computed in cohomology so we replace their factors by analytic
sets). Here, we use that $\widetilde{F}_n(c,z)$ is a polynomial of degree $d^n$ in
the $z_i$ and of degree $\leq C d^{n}$ in the $c_{1,i}$ and in the
$c_{2,i}$.

As a result, we have that:
\begin{equation*}
\| \widetilde{F}_n^*(\omega_3) \| \leq C d^n. 
\end{equation*} 
and consequently:
\begin{equation*}
\| \widetilde{F}_n^*(\omega_3^s) \| \leq C d^{sn}. 
\end{equation*} 
We also remark that:
$$ \widetilde{F}_n^*(\omega_1)=\omega_1$$
and
$$ \widetilde{F}_n^*(\omega_2)=\omega_2$$
since $\widetilde{F}_n$ acts as the identity on $\widetilde{W}$. \\

\noindent {\bf Construction of a function $g$ that tests Hypothesis \ref{distance}}

We can write $\widetilde{F}_1^*(\omega_3^s)$ in cohomology:
\begin{equation*}
\widetilde{F}_1^*(\omega_3^s)=\sum_{i_1+i_2+i_3=s} a_{i_1,i_2,i_3} \omega_1^{i_1}\wedge \omega_2^{i_2}\wedge \omega_3^{i_3} + dd^c \mathcal{U} 
\end{equation*}
where $\mathcal{U}$ is a negative $(s-1,s-1)$ current, which is $C^1$
where $\widetilde{F}_1^*(\omega_3^s)$ is smooth (see Proposition
2.3.1 in \cite{DS6} and observe that $\widetilde{W}\times \P^k$ is
homogeneous). We also denote $\widetilde{\Omega}= \sum_{i_1+i_2+i_3=s}
a_{i_1,i_2,i_3} \omega_1^{i_1}\wedge \omega_2^{i_2}\wedge
\omega_3^{i_3} $ the smooth form cohomologuous to
$\widetilde{F}_1^*(\omega_3^s)$. Testing against $\omega_1^{a}\wedge
\omega_2^{b} \wedge \omega_3^{b}$ for $a+b+c=l+k-s$ gives that
$a_{i_1,i_2,i_3} \geq 0$. In what follows, we take a particular
$\mathcal{U}$. We explain now its construction.

The indeterminacy set of $\widetilde{F}_1$ has dimension $l+k-s-1$ (it
is obvious in $W \times \P^k$ and in $(\widetilde{W}
\setminus W) \times \P^k$, use a stratification with the dimension of the kernel of $c_1$
and $c_2$). In particular, by Theorem 4.5 in \cite[Chapter III]{dem2}, we have that
$\widetilde{F}_1^*(\omega_3^s)=(\widetilde{F}_1^*(\omega_3))^s$.  Let
$u$ be a quasi-potential of $\widetilde{F}_1^*(\omega_3)$ and $\beta$
be a Kähler form cohomologuous to $\widetilde{F}_1^*(\omega_3)$. In
other words, $\widetilde{F}_1^*(\omega_3)=\beta+dd^c u $. We can write $\mathcal{U}$ as in the proof of Theorem \ref{convergencecurrent}, that is:
$$\mathcal{U}= \sum_{j=0}^{s-1} u \widetilde{F}_1^*(\omega_3)^{s-1-j} \wedge
\beta^j .$$
In this case $\widetilde{\Omega}=\beta^s$.

 Consider the rational map:
 \begin{align*}
\sigma: \widetilde{W} \times \P^k &\to \widetilde{W} \times \P^k \\
                            (c,z) &\mapsto  (c,c_2(z)).
\end{align*}
We define $\mathcal{I}^-:=\sigma(\widetilde{W} \times I^-) $.
 It is an analytic set of $\widetilde{W} \times \P^k$ of dimension $l+s-1$ 
 such that for $c \in W$, $\mathcal{I}^- \cap \{c\}\times \P^k=I^-(f_c)=c_2(I^-)$.

Let $[\mathcal{I}^-]$ denote the current of integration on $\mathcal{I}^-$, 
it is of bidimension $(l+s-1,l+s-1)$. Consider the set
 $$V_n := \{c \in W, f^n_c(I^-(f_c))\cap I^+(f_c)=\varnothing  \},$$ 
it is a Zariski open set in $\widetilde{W}$. 

Now, consider 
$$\varphi_n:=\Pi_*( \sum_{j=0}^{s-1} \widetilde{F}_n^* u  \widetilde{F}_{n+1}^*(\omega_3)^{s-1-j} \wedge
\widetilde{F}_n^* (\beta)^j \wedge   [\mathcal{I}^-]  ) $$
where $\Pi$ is the canonical projection from $\widetilde{W} \times
\P^k$ to $\widetilde{W}$.

This function and its $dd^c$ are well defined since the dimension of the indeterminacy set of
the $\widetilde{F}_n$ is $l+k-s-1$, the dimension of $\mathcal{I}^-$
is $l+s-1$ and the dimension of the intersection of these sets is less
than $l-1$ (we use again the Theorem 4.5 in \cite[Chapter III]{dem2}
and a stratification of $\widetilde{W} \setminus W$ with the dimension
of the kernel of $c_1$ and $c_2$).

On the set $\cap_{i=0}^{n+1} V_i$,
$\varphi_n$ is continuous since $\Pi$ restricted to $\mathcal{I}^-$ is a submersion so the
push-forward of a continuous form is continuous (we can remove to $\cap_{i=0}^{n+1} V_i$ the fibers
of $\Pi$ which are contained in the singular locus of $\mathcal{I}^-$
because it is an analytic subset in $W$). Finally, we define on
$\cap_{n=0}^{N+1} V_n$, 
$$g_N:= \sum_{n=0}^{N} d^{-sn} \varphi_n.$$\\

\noindent {\bf Computation of $dd^c g_N$}

 We have that:
$$ dd^cg_N=\Pi_*\left( \sum_{n=0}^{N} d^{-sn}   \sum_{j=0}^{s-1} dd^c
 \widetilde{F}_n^* u \wedge  \widetilde{F}_{n+1}^*(\omega_3)^{s-1-j} \wedge
\widetilde{F}_n^* (\beta)^j \wedge   [\mathcal{I}^-] \right).$$

Recall that $dd^c u = \widetilde{F}_1^*(\omega_3)-\beta$. So,

$$dd^c \widetilde{F}_n^* u=
\widetilde{F}_{n+1}^*(\omega_3)-\widetilde{F}_{n}^*(\beta)$$

since is is true outside a set of dimension $l+k-s-1 \leq l+k-2$.

We obtain
\begin{align*}
&dd^cg_N=\\
&\Pi_*\left( \sum_{n=0}^{N} d^{-sn}   \sum_{j=0}^{s-1} (\widetilde{F}_{n+1}^*(\omega_3)^{s-j} \wedge
\widetilde{F}_n^* (\beta)^j -   \widetilde{F}_{n+1}^*(\omega_3)^{s-1-j} \wedge
\widetilde{F}_n^* (\beta)^{j+1})   \wedge   [\mathcal{I}^-]     \right)
\end{align*}
which is equal to
$$ dd^cg_N=\Pi_*\left( \sum_{n=0}^{N} d^{-sn}
(\widetilde{F}_{n+1}^*(\omega_3^{s}) - \widetilde{F}_n^* (\beta^{s}))   \wedge   [\mathcal{I}^-]     \right).$$

Recall that $\beta^s=\widetilde{\Omega}= \sum_{i_1+i_2+i_3=s}
a_{i_1,i_2,i_3} \omega_1^{i_1}\wedge \omega_2^{i_2}\wedge
\omega_3^{i_3} $. We show now that $a_{0,0,s}=d^s$. First we have that:
$$\langle \widetilde{\Omega}, \omega_3^{k-s} \wedge \omega_1^{N} \wedge \omega_2^N \rangle = a_{0,0,s}. $$ 
Since $\widetilde{\Omega}$ and $\widetilde{F}^*_1(\omega_3^s)$ are cohomologuous, we deduce:
$$\langle \widetilde{F}^*_1(\omega_3^s), \omega_3^{k-s} \wedge \omega_1^{N} \wedge \omega_2^N \rangle = a_{0,0,s}. $$ 
So, we want to compute:
$$ \int \widetilde{F}^*_1(\omega_3^s)\wedge \omega_3^{k-s} \wedge \omega_1^{N} \wedge \omega_2^N .$$
This can be done in cohomology. If $L_{k-s}$ is a generic analytic
subspace of dimension $k-s$ in $\P^k$ and $L_s$ is a generic analytic
subspace of dimension $s$ in $\P^k$ and $\{c\} \times \P^k$ is a
slice, then the previous quantity is the number of intersections of
$f^{-1}_c(L_{k-s}) \cap L_s$ on the slice. This is equal to $d^s$
since the degree $d_s$ of $f_c$ is $d^s$ on $W$ which is a Zariski
open set in $\widetilde{W}$, so we have indeed that $d_s=a_{0,0,s}$. \\

We have the equality:
$$\widetilde{F}_n^*(\widetilde{\Omega})=d^s \widetilde{F}_n^*(\omega_3^s)+ \sum_{i_1+i_2+i_3=s, \ i_3 \neq s} a_{i_1,i_2,i_3} \widetilde{F}_n^*(\omega_1^{i_1}\wedge \omega_2^{i_2}\wedge \omega_3^{i_3}). $$
We denote the second term on the right-hand side by $S_n$. Since $ \widetilde{F}_n^*(\omega_1)=\omega_1$ and $ \widetilde{F}_n^*(\omega_2)=\omega_2$, we can bound the mass of $S_n$ :
\begin{align*}
\|S_n \| &\leq \sum_{i_1+i_2+i_3=s, \ i_3 \neq s} a_{i_1,i_2,i_3} \|\omega_1^{i_1}\wedge \omega_2^{i_2}\wedge  \widetilde{F}_n^*(\omega_3^{i_3})\| \\
         &\leq C d^{n(s-1)}
\end{align*}
since $\|\widetilde{F}_n^*(\omega^j)\|\leq C d^{jn}$ for $j \leq s$. So replacing in $dd^c g_N$, we have:
$$ dd^cg_N=\Pi_*\left( \sum_{n=0}^{N} d^{-sn} (\widetilde{F}_{n+1}^*(\omega_3^s)-d^s\widetilde{F}_{n}^*(\omega_3^s)    ) \wedge   [\mathcal{I}^-]   \right)- \Pi_*\left( \sum_{n=0}^{N} d^{-sn} S_n \wedge   [\mathcal{I}^-]   \right) .$$
The second term in the right-hand side is a positive closed currents
with mass uniformly bounded in $N$ by the above. We control the first term. Reorganizing the sum, we see that it is equal to:
$$ \Pi_*( (d^{-sN} \widetilde{F}^*_{N+1}(\omega_3^s) -d^sw_3^s)\wedge  [\mathcal{I}^-]).$$
Using the fact that the mass of the positive closed current $\widetilde{F}^*_{N+1}(\omega_3^s)$ is bounded by $Cd^{s(N+1)}$ gives that:
\begin{align*}
dd^cg_N= \widetilde{\Omega^+_{1,N}} - \widetilde{\Omega_{2,N}^+}
\end{align*}
where $ \widetilde{\Omega^+_{i,N}}$ is a positive closed current of bidegree
$(1,1)$ with $\| \widetilde{\Omega^+_{i,N}} \| \leq C$ where $C$ is independent of
$N$. We can write

$$\widetilde{\Omega^+_{i,N}}= a_{i,N} \omega_1 + b_{i,N} \omega_2 +
dd^c \psi_{i,N}$$
with $ a_{i,N}$ and $b_{i,N}$ smaller than $C$. We explain now what is
the normalization on the qpsh functions $\psi_{i,N}$ that we take.\\

We say that a measure is $PLB$ if the qpsh
functions are integrable for the measure. Any measure given by a smooth distribution is $PLB$. In
particular, we can find a $PLB$ probability measure that we denote
$\nu$ with support in the $W'_0$ defined previously. We have the following lemma (see Proposition 2.4 in \cite{DS2}):
\begin{lemme}
 The family of qpsh functions in $\widetilde{W}$ such that
 $dd^c \psi \geq - \Omega= - (\omega_1 + \omega_2)$ and one of the two following conditions:
 $$\max_{\widetilde{W}} \psi =0 \ \text{or} \ \int \psi d\nu=0 $$
is bounded  in $L^1(\nu)$ and is bounded from above.
\end{lemme} 
When, we write 
$$\widetilde{\Omega^+_{i,N}}= a_{i,N} \omega_1 + b_{i,N} \omega_2 +
dd^c \psi_{i,N}$$
we suppose that we take the normalization $\int \psi_{i,N} d\nu=0$.\\ 

\noindent {\bf Link with Hypothesis \ref{distance}}

Let $c\in \cap_{i\leq n+1} V_i$, then we want to show that:
$$\varphi_n(c)=\int_{f_c^n(I^-(f_c))} U_c$$ 
where $(f_c)^*(\omega^s)=d^s\omega^s+dd^c U_c$. 

First, when $c\in \cap_{i\leq n+1} V_i$, we have:

$$\varphi_n(c)=\int_{I^-(f_c)} F_n^{*}( \mathcal{U})_{| \{c \} \times
  \P^k}.$$

Here, $F_n^{*}( \mathcal{U})_{| \{c \} \times
  \P^k}$ is the restriction of $F_n^{*}( \mathcal{U})$ on $ \{c \} \times
  \P^k$ which is well defined because $F_n^{*}( \mathcal{U})$ is
  continous near $\{ c \} \times I^-(f_c)$. But $F_n^{*}( \mathcal{U})_{| \{c \} \times
  \P^k}$ is equal to $(f_c^n)^{*}( \mathcal{U}_{|\{c \} \times
  \P^k})$, so

$$\varphi_n(c)=\int_{f_c^n(I^-(f_c))}  \mathcal{U}_{| \{c \} \times
  \P^k}.$$

Recall that

$$\mathcal{U}= \sum_{j=0}^{s-1} u \widetilde{F}_1^*(\omega_3)^{s-1-j} \wedge
\beta^j .$$

In particular, $\mathcal{U}_{| \{c \} \times  \P^k}$ near
$f_c^n(I^-(f_c))$ can be written
$$ U_c=\sum_{j=0}^{s-1} u_{| \{c \} \times  \P^k} f_c^*(\omega)^{s-1-j} \wedge
d^j \omega^j$$
because the coefficient of $\omega_3$ in $\beta$ is $d$.

The singularities of $u_{| \{c \} \times  \P^k}$ are in $I^+(f_c)$, so
by Theorem 4.5 in \cite[Chapter III]{dem2}, we have that $U_c$ and $dd^c U_c= \sum_{j=0}^{s-1} dd^c (u_{| \{c \}
  \times  \P^k}) \wedge f_c^*(\omega)^{s-1-j} \wedge
d^j \omega^j$ are well defined in all $\P^k$. But, if we take the restriction of the equation
$\widetilde{F}_1^*(\omega_3)=\beta+dd^cu $ on $\{c \} \times  \P^k$, we obtain

$$f_c^{*} (\omega)= d  \omega + (dd^c u)_{|\{c \} \times  \P^k}= d
\omega + dd^c (u_{|\{c \} \times  \P^k})$$

since it is true outside $I^+(f_c)$ and $f_c^{*} (\omega)$ or $dd^c
(u_{|\{c \} \times  \P^k})$ have no mass on this set of dimension $k-s-1$. Moreover $u$ is
a qpsh function, so it takes a value at every point.

Finally,

$$dd^c U_c=\sum_{j=0}^{s-1} (f_c^{*} (\omega)- d  \omega) \wedge f_c^*(\omega)^{s-1-j} \wedge
d^j \omega^j=(f_c)^*(\omega^s)-d^s\omega^s. $$

\noindent  {\bf Proof of the genericity}

Recall that $\varphi_n$ is continuous on $\cap_{i \leq n+1} V_i$. This
implies  that $g_N$ is continuous on $\cap_{i \leq N+1} V_i$ and it
decreases to a function $g$ on $\cap_{i\geq 0} V_i$ with $g$ usc on
$\cap_{i\geq 0} V_i$. It means that for every point $x$ in
$\Lambda=\cap_{i\geq 0} V_i$, we have $\limsup_{y \rightarrow x \mbox{, }
  y \in \Lambda} g(y) \leq g(x)$.\\ 

Let $m_N =\int g_N d\nu$. We can write on $\widetilde{W}$, $g_N-m_N = \psi_{1,N}-\psi_{2,N}$. Here the equality is true on
a set of full Lebesgue measure in $\widetilde{W}$. But, since $g_N$ is continuous
on $\cap_{i \leq N+1} V_i$ and the $\psi_{i,N}$ are qpsh, the equality
is true for every point in $\cap_{i \leq N+1} V_i$ (see below the
proof of the inequality $g -m \geq  \psi_1 - \psi_2$).

We apply the previous lemma to the sequences $\psi_{i,N}$ and we have
that these sequences are uniformly bounded from above and bounded in $L^1( \nu)$. So we can extract converging
subsequences to some limit points $\psi_1$ and $\psi_2$ in $L^1$. The sequence
$m_N$ is bounded thanks to the definition of $W'_0$ and of
$\nu$. So $m_N$ converges to $m$ by monotone convergence. In particular, $g-m=\psi_1 - \psi_2$ up to a set
of zero Lebesgue measure in $\widetilde{W}$. We want to show now that
we have

$$g -m \geq  \psi_1 - \psi_2$$ 
for every point in $\Lambda=\cap_{i\geq 0} V_i$. Indeed, assume there
is a point $x \in \Lambda$ such that $(g  +\psi_2)(x) < m + \psi_1(x)
-\varepsilon$. On a chart which contains $x$, we can write $\psi_1=
\lambda_1 + \xi_1$ with $\lambda_1$ smooth and $\xi_1$ psh.

Since $g$ and $\psi_2$ are usc on $\Lambda$, so is
their sum and so $(g +\psi_2)(y) < m + \psi_1(x) - \varepsilon /2$
on a small ball $B(x,r)$ centered at $x$ and of radius $r$ (for $y \in
\Lambda$). For a function $h$, we denote by $m_{B(x,r)}(h)$ the mean
value of $h$ on the ball $B(x,r)$. We have that that $m_{B(x,r)}(g
+\psi_2)=m_{B(x,r)}(m + \lambda_1 + \xi_1)$ since both functions are equal
a.e. and $m_{B(x,r)}(\xi_1) \geq \xi_1(x)$ since $\xi_1$ is psh, so
$$m+ \psi_1(x)- \varepsilon /2\geq m_{B(x,r)}(g +\psi_2) \geq
m_{B(x,r)}(\lambda_1) +m+ \xi_1(x)$$
which is false if we take $r$ small enough to have
$m_{B(x,r)}(\lambda_1)$ near $\lambda_1(x)$.

  In particular, the set of points where $g=-\infty$ is pluripolar
  since it is included in the set of points where $\psi_1$ is $-\infty$. By the proof of Theorem \ref{convergencecurrent}, we see that $g \neq -\infty$ is equivalent to the fact that the first half of Hypothesis \ref{distance} is satisfied.

 We do the same thing for the second half of Hypothesis \ref{distance}
 and we conclude since the intersection of two pluripolar sets is
 pluripolar. \hfill $\Box$ \hfill   \\

The results of this section remain valid for $f^{-1}$.
So we can construct the Green current of order $k-s$ for $f^{-1}$ that we denote by $T^-_{k-s}$.

\section{The equilibrium measure} 
\subsection{Construction of the measure}
We want to define the equilibrium measure $\mu$ as $T^+_{s}\wedge T^-_{k-s}$. In \cite{BD1}, the authors used an approach based on the energy. More precisely, they show that the potential of the Green current is in the Hilbert space $H_{T^-}$ defined by the closure of the smooth forms for the norm $\sqrt{ \int d\varphi \wedge d^c\varphi \wedge T^-}$. They deduced from that fact that the measure $T^+ \wedge T^-$ is well defined and that the potential of the Green current is integrable with respect to that measure. \\

Such an approach cannot be adapted here since the super-potential is not a function defined on $\mathbb{P}^k$ but a function defined on $\C_{k-s+1}$. Instead, we will use the formalism of super-potential. See Definition \ref{def_wedge} for the definition of wedgeability. We prove the theorem:
\begin{theorem}\label{constructionmeasure} 
The current  $T^+_{s}$ and $T^-_{k-s}$ are wedgeable. So the intersection $T^+_{s}\wedge T^-_{k-s}$ is a well defined measure $\mu$ and the quasi-potential of the Green current of order $1$ is integrable with respect to this measure.
\end{theorem} 
Recall that $T^+:=T^+_1$ is a well defined invariant current in $\C_1$ (\cite{Sib}) and that it admits the quasi-potential:
$$G=\sum_n (\frac{1}{d}f^*)^n u$$ 
where $u<0$ is a quasi-potential of the current $d^{-1}f^*(\omega)$ (we write $u$ instead of $U_{L(\omega)}$ in order to simplify the notations).
We denote as before $L_n$ and $\Lambda_n$ the normalized pull-pack and push-forward associated to $f^n$.
 In what follows, for $q\leq s$,  $\U_{L^m(\omega^q)}$ denotes the super-potential of $L^m(\omega^q)$ of the previous section, that is:
$$ \U_{L^m(\omega^q)}=\sum_{n=0}^{ m-1} \frac{1}{d^n }\U_{L(\omega^q)} \circ \Lambda^n,$$
on smooth forms in $\C_{k-q+1}$, where $\U_{L(\omega^q)}$ is a negative super-potential of $L(\omega^q)$. By Corollary \ref{composition_smooth}, we can write it as:
$$ \U_{L^m(\omega^q)}=\sum_{n=0}^{ m-1} \frac{1}{d^n }\U_{L(\omega^q)} \circ \Lambda_n $$ 
on smooth forms in $\C_{k-q+1}$. Then Lemma \ref{stability} assures us that if $S\in \C_{k-q+1}$ is $(f^m)_*$-admissible, it is also $(f^n)_*$-admissible for $n\leq m$. So by definition of super-potentials and by Hartogs' convergence we have that 
$$ \U_{L^m(\omega^q)}=\sum_{n=0}^{ m-1} \frac{1}{d^n }\U_{L(\omega^q)} \circ \Lambda_n $$ 
on $(f^m)_*$-admissible currents in $\C_{k-q+1}$. Again, Lemma \ref{stability} gives that on $(f^m)_*$-admissible currents in $\C_{k-q+1}$, we have that $\Lambda_n=\Lambda^n$ hence:
 \begin{equation}\label{Greenadm}
 \U_{L^m(\omega^q)}=\sum_{n=0}^{ m-1} \frac{1}{d^n }\U_{L(\omega^q)} \circ \Lambda^n,
 \end{equation} 
on $(f^m)_*$-admissible currents in $\C_{k-q+1}$. \\

 We need the following lemma to construct the measure.
\begin{lemme}\label{defwedgepull}
The current $\omega^{s-q}\wedge L^n(\omega^q)$ and $T^-_{k-s}$ are wedgeable for all $n\geq 0$ and $0\leq q \leq s$. Furthermore, for all integers $n$ and $n'$ and $1\leq q\leq s+1$ we have that
$\U_{L^n(\omega)^q}(L^{n'}(\omega)^{s-q+1}\wedge T^-_{k-s})$ is finite.   
\end{lemme}
\emph{Proof.} We have seen in Proposition \ref{Hregular} that $\omega^{s-q}\wedge L(\omega)^{q}$ and $T^-_{k-s}$ are wedgeable for $q\leq s$ and that the super-potentials of $T^-_{k-s}$ are finite at $L(\omega)^{s+1}$. So applying that to $f^n$ instead of $f$, we have that $\omega^{s-q}\wedge L^n(\omega)^q$ and $T^-_{k-s}$ are wedgeable for all $n\geq 0$ and that the super-potentials of $T^-_{k-s}$ are finite at $L^n(\omega)^{s+1}$. \\

The case where $q=s+1$ is already known so we assume $1 \leq q \leq s$. 
The current $L(\omega)^{s+1-q}$ and $L(\omega)^{q}$ are wedgeable and their wedge-product is $L(\omega)^{s+1}$ (it follows from Corollary 4.11 Chapter III in \cite{dem2} and Lemma \ref{lemma_wedge_equi}). So using Lemma \ref{dec_wedge} we have that a super-potential of $L(\omega)^{s+1}$ is given by:
$$ \U_{L(\omega)^{q}}(L(\omega)^{s-q+1}\wedge R) +  \U_{L(\omega)^{s-q+1}}(\omega^{q}\wedge R),$$
on current $R\in \C_{k-s}$ such that $R$ and $L(\omega)^{s-q+1}$ are wedgeable. In particular, we can take $R=T^-_{k-s}$ at which point the super-potential of $L(\omega)^{s+1}$ is finite.
A super-potential of $L(\omega)^{s-q+1}\wedge \omega^{q}$ is given by:
$$\U_{L(\omega)^{s-q+1}}(\omega^{q}\wedge \star).$$
So by difference,
$$\U_{L(\omega)^{q}}(L(\omega)^{s-q+1}\wedge T^-_{k-s})$$
 is well defined in the sense of super-potentials (that is it is continuous for the Hartogs' convergence) and is finite. 

So we have proved the lemma for $n=n'=1$. \\

Applying the result to $f^n$ gives the lemma for $n=n'$. Now, let $n\leq n'$. Then  $L^n(\omega)^q$ is more H-regular than $L^{n'}(\omega)^q$. The super-potentials of $L^{n'}(\omega)^q$ are finite at $L^{n'}(\omega)^{s-q+1}\wedge T^-_{k-s}$ so the super-potentials of $L^{n}(\omega)^q$ are also finite at $L^{n'}(\omega)^{s-q+1}\wedge T^-_{k-s}$.

Similarly, let $n\geq n'$. Then $L^{n'}(\omega)^{s-q+1}$ is more H-regular than $L^{n}(\omega)^{s-q+1}$ and so Lemma \ref{lemma_wedge_regular} implies that $L^{n'}(\omega)^{s-q+1}\wedge T^-_{k-s}$ is more H-regular than $L^{n}(\omega)^{s-q+1}\wedge T^-_{k-s}$. The super-potentials of $L^{n}(\omega)^q$ are finite at $L^{n}(\omega)^{s-q+1}\wedge T^-_{k-s}$, which means by symmetry of the super-potentials that the super-potentials of $L^{n}(\omega)^{s-q+1}\wedge T^-_{k-s}$ are finite at  $L^{n}(\omega)^q$. Hence the super-potentials of $L^{n'}(\omega)^{s-q+1}\wedge T^-_{k-s}$ are finite at  $L^{n}(\omega)^q$ which means that the super-potentials of $L^{n}(\omega)^q$ are finite at $L^{n'}(\omega)^{s-q+1}\wedge T^-_{k-s}$.  That gives the lemma. \hfill $\Box$ \hfill \\

\noindent \emph{Proof of Theorem \ref{constructionmeasure}}
By the above lemma, we have that $L^n(\omega^s)\wedge T^-_{k-s}$ is $(f^n)_*$-admissible since it is finite at $\U_{L^n(\omega)}$. Hence by Lemma \ref{stability}, we have that $\Lambda^n(L^n(\omega^s)\wedge T^-_{k-s})$ is well defined and equal to $\Lambda_n(L_n(\omega^s)\wedge T^-_{k-s})$ (recall that Corollary \ref{composition_smooth} gives $L_n(\omega^s)=L^n(\omega^s)$). 

We consider:
$$ \frac{1}{d^n}\U_{L(\omega)} (\Lambda_n(L_n(\omega^s)\wedge T^-_{k-s})).$$
It is finite since by Lemma \ref{pushpull} applied to $f^n$ and the invariance ot $T^-_{k-s}$, it is equal to  
$$\frac{1}{d^n}\U_{L(\omega)} (\omega^s\wedge T^-_{k-s}),$$ 
and the previous lemma assures us that this is finite.

Using Lemma \ref{decpushpull} for $f^n$ instead of $f$, we see that it is equal to:
$$\U_{L^n(\omega^s)}(L^{n+1}(\omega)\wedge T^-_{k-s})-\U_{L^n(\omega^s)}(L^{n}(\omega)\wedge T^-_{k-s})+ \left(\frac{1}{d}\right)^n \U_{L(\omega)}(\Lambda_n(\omega^s\wedge T^-_{k-s})).$$
We now perform some sort of Abel transform. We sum from $0$ to $N$ and we regroup the terms in $L^{n}(\omega)$ (observe for the first term that $\U_{\omega^s}=0$):
\begin{align}\label{key}
\sum_{n=0}^{N} \frac{1}{d^n}\U_{L(\omega)} (\omega^s\wedge T^-_{k-s}) &= \sum_{n=1}^{N} (-\U_{L^n(\omega^s)}+\U_{L^{n-1}(\omega^s)})(L^{n}(\omega)\wedge T^-_{k-s}) \\
 &+\U_{L^N(\omega^s)}(L^{N+1}(\omega)\wedge T^-_{k-s})+ \sum_{n=0}^{N}  \frac{1}{d^n}\U_{L(\omega)}(\Lambda_n(\omega^s\wedge T^-_{k-s})) \nonumber
\end{align} 

 Now, $\U_{L^{n-1}(\omega^s)}-\U_{L^{n}(\omega^s)}=-d^{-n+1} \U_{L(\omega^s)}\circ \Lambda^{n-1}$  on smooth forms. By Corollary \ref{composition_smooth}, we can write it as:
 $$\U_{L^{n-1}(\omega^s)}-\U_{L^{n}(\omega^s)}=-d^{-n+1} \U_{L(\omega^s)}\circ \Lambda_{n-1}.$$ 
 on smooth forms. Let $T \in \C_{k-s+1}$ be $(f^n)_*$-admissible, then $T$ is $(f^{n-1})_*$-admissible by Lemma \ref{stability}. Taking a sequence of smooth currents converging in the Hartogs' sense to $T$ and using that $\Lambda_{n-1}$ is continuous for the Hartogs' convergence (Theorem \ref{continuiteL}), we have that:
  $$\U_{L^{n-1}(\omega^s)}-\U_{L^{n}(\omega^s)}=-d^{-n+1} \U_{L(\omega^s)}\circ \Lambda_{n-1},$$  
 on  $(f^n)_*$-admissible currents (observe that $\U_{L^{n-1}(\omega^s)}$ and $\U_{L^{n}(\omega^s)}$ are finite on $(f^n)_*$-admissible currents). 
 In particular, we consider the current  $L^{n}(\omega)\wedge T^-_{k-s}$ which is $(f^{n})_*$-admissible by the previous lemma. So using again Lemma \ref{pushpull} for $f^{n-1}$ gives:
\begin{align*}
(\U_{L^{n-1}(\omega^s)}-\U_{L^{n}(\omega^s)})(L^{n}(\omega)\wedge T^-_{k-s}) &= -d^{-n+1} \U_{L(\omega^s)}(\Lambda_{n-1}(L^{n}(\omega)\wedge T^-_{k-s})) \\
																																						&= -d^{-n+1} \U_{L(\omega^s)}(\Lambda_{n-1}(L_{n-1}(L(\omega))\wedge T^-_{k-s})) \\
																																					 &=-d^{-n+1} \U_{L(\omega^s)}(L(\omega)\wedge T^-_{k-s}).
\end{align*}
So the series $\sum_{n=1}^{N} (-\U_{L^n(\omega^s)}+\U_{L^{n-1}(\omega^s)})(L^{n}(\omega)\wedge T^-_{k-s})$ is also convergent thanks to the previous lemma.  We also have that $\U_{L^N(\omega^s)}(L^{N+1}(\omega)\wedge T^-_{k-s})$ is negative since $\U_{L^N(\omega^s)}$ is negative.
Thus, letting $N$ go to $\infty$:
\begin{align*}
\left(\sum_{n\geq 0}\frac{1}{d^n} \right) \U_{L(\omega)} (\omega^s\wedge T^-_{k-s}) &+ \left(\sum_{n\geq 1} d^{-n+1}\right) \U_{L(\omega^s)}(L(\omega)\wedge T^-_{k-s}) \\
                                                                        & \qquad \leq \sum_{n \geq 0}  \frac{1}{d^n}\U_{L(\omega)}(\Lambda^n(\omega^s\wedge T^-_{k-s})).
\end{align*}
We recognize by (\ref{Greenadm}) that the right-hand side is in fact $\U_{T^+}(\omega^s\wedge T^-_{k-s})$ which in term of quasi-potential is  $\int G \omega^s\wedge T^-_{k-s}$ (recall that $T^+$ is the Green current of order $1$). Thus by Hartogs' convergence, we have that $\U_{T^+}(\omega^s\wedge T^-_{k-s})$ is finite (we could also conclude by monotone convergence that $G\in L^1(\omega^s\wedge T^-_{k-s})$, both properties being equivalent).  \\
 
 Observe now that in (\ref{key}) every term converge. In particular, 
 $$(\U_{L^N(\omega^s)}(L^{N+1}(\omega)\wedge T^-_{k-s}))_N $$ 
 converges to a finite value. Using Lemma \ref{byparts}, we have the identity:
 \begin{align*}
\U_{L^N(\omega^s)}(L^{N+1}(\omega)\wedge T^-_{k-s})&= \U_{L^{N+1}(\omega)}( L^N(\omega^s) \wedge T^-_{k-s})  \\
 &-  \U_{L^{N+1}(\omega)}(\omega^s \wedge T^-_{k-s})   +\U_{L^N(\omega^s)}(\omega \wedge T^-_{k-s}). 
 \end{align*}
 On the right-hand side, the first and third terms are negative, the third term is decreasing and we just proved that the second term converges to $\U_{T^+}(\omega^s \wedge T^-_{k-s})$ which is finite. That implies that every term is in fact convergent. 
 
 In particular, we have the convergence of $ \U_{L^N(\omega^s)}(\omega\wedge T^-_{k-s})$. Since $L^N(\omega^s) \to T^+_s$ in the Hartogs' sense, that means that $ \U_{T_s^+}(\omega\wedge T^-_{k-s})$ is finite. Hence the current $T^+_s$ and $T^-_{k-s}$ are wedgeable and their intersection is a well defined probability measure $\mu$ (we could also have deduced that from the convergence of the first term but this is more natural).   

Recall that the  function $(R,S)\to\U(R,S):=\U_R(S)=\U_{S}(R)$ for $R$ and $S$ in $\C_q$ and $\C_{k-q}$ ($\U_{R}$ and $\U_{S}$ being the super-potentials of mean 0) is upper semi-continuous. The convergence of $\U(L^{N+1}(\omega), L^N(\omega^s)\wedge T^-_{k-s})$ implies that $\U(T^+, T^+_s\wedge T^-_{k-s})$ is finite which means exactly that the quasi-potential of the Green current is integrable with respect to $\mu$. \hfill $\Box$ \hfill \\

Of course, the potential of the Green current of order $1$ of $f^{-1}$ is  also integrable for the measure $\mu$.
\begin{corollaire}
The measure $\mu$ gives no mass to the indeterminacy sets $I^+$ and $I^-$. Furthermore, $L(\mu)=f^*(\mu)$ and $\Lambda(\mu)=f_*(\mu)$ are well defined in the sense of super-potentials.
\end{corollaire} 
\emph{Proof.} The fact that $\mu$ is $f_*$-admissible follows from Theorem \ref{constructionmeasure} since its super-potentials are finite at the point $L(T^+)=T^+$ and so they are finite at the point $L(\omega)$ which is more H-regular than $L(T^+)$. Since the potential of $T^+$ is equal to $-\infty$ on $I^+$ and is in $L^1(\mu)$ (in fact $\log \text{dist}(x, I^+)\in L^1(\mu)$) we have that $\mu$ gives no mass to the indeterminacy set $I^+$, similarly for $I^-$.  \hfill $\Box$ \hfill 
\begin{proposition}\label{invariancemeasure}
The measure $\mu$ is invariant, that is $f^*(\mu)$ and $f_*(\mu)$ are equal to $\mu$.
\end{proposition}
\emph{Proof.} The currents $L^n(\omega^s)$ and $\Lambda^m(\omega^{k-s})$ are wedgeable for $m$ and $n$ in $\mathbb{N}$ since they are more H-regular than $T^+_s$ and $T_{k-s}^-$. So let $\mu_n:=L^n(\omega^s)\wedge \Lambda^n(\omega^{k-s})$ (resp. $\mu'_n:=L^{n-1}(\omega^s)\wedge \Lambda^{n+1}(\omega^{k-s})$). Now since $L^n(\omega^s)$ and $\Lambda^n(\omega^{k-s})$ converge in the Hartogs' sense to $T^+_{s}$ and $T^-_{k-s}$ which are wedgeable, we have that $\mu_n$ (resp. $\mu'_n$) converges to $\mu$ in Hartogs' sense (Proposition  \ref{prop_cv_wedge_hartogs}). 

By Lemma \ref{pushpull}, we have that $\mu'_n= \Lambda(L^n(\omega^s)\wedge \Lambda^n(\omega^{k-s}))=\Lambda(\mu_n)$ (observe that $L^n(\omega^s)\wedge \Lambda^n(\omega^{k-s})$ is $f_*$ admissible since it is more H-regular than $T_s^+\wedge T_{k-s}^-$ which is $f_*$-admissible). 

So, since $\mu$ is $f_*$-admissible, we have that $\mu'_n$ converges in the Hartogs' sense to $\Lambda(\mu)=\mu$ which is what we wanted. 
\hfill $\Box$ \hfill 

\begin{corollaire}
The measure $\mu$ gives no mass to the indeterminacy sets $I(f^{\pm n})$ and the critical sets $\C(f^{\pm n})$.
\end{corollaire}
\emph{Proof.} We already know that the indeterminacy sets have no mass for $\mu$ so using the invariance of $\mu$, we have that $\mu(\C(f))=\mu(f^{-1}(I^-))=\mu(I^-)=0$.\hfill  $\Box$  \hfill

\subsection{Green currents of order $1\leq q\leq s$}
The purpose of this paragraph is to construct the Green currents of order $q$ for $q\leq s$. This will allow us to prove that $T^+_s$ can be written as $(T^+)^s$. As an application, we show that the equilibrium measure gives no mass to the pluripolar sets.

Using the same arguments than in Theorem  \ref{constructionmeasure}, we construct the Green currents $T^+_q$ of order $q$ for $q\leq s$:
\begin{proposition}\label{extension}
For $1\leq q\leq s$, the sequence $(L^n(\omega^q))_n$ converges in the Hartogs' sense to $T^+_q$ \emph{the Green current of order $q$} and the Green currents $T^+_{q}$ and $T^-_{k-s}$ are wedgeable. Furthermore, any super-potential $\U_{T^+_q}$ of $T^+_q$ satisfies 
$$ \U_{T^+_q}(T^+_{s-q+1}\wedge T^-_{k-s})>-\infty.$$ 
\end{proposition}
\emph{Proof.} Observe that the roles of $q$ and $s-q+1$ are symmetric, so anything proved for $q$ stands for $s-q+1$.
The current $L^n(\omega^{s-q+1})\wedge T^-_{k-s}$ is $(f^n)_*$ admissible by Lemma \ref{defwedgepull}. Lemma \ref{stability} implies that $\Lambda^n(L^n(\omega^{s-q+1})\wedge T^-_{k-s})$ is well defined and equal to $\Lambda_n(L^n(\omega^{s-q+1})\wedge T^-_{k-s})$. So we consider this time:
$$ \frac{1}{d^n}\U_{L(\omega^q)} (\Lambda_n(L^n(\omega^{s-q+1})\wedge T^-_{k-s})).$$
By Lemma \ref{pushpull} and the invariance ot $T^-_{k-s}$, it is equal to  
$$\frac{1}{d^n}\U_{L(\omega^q)} (\omega^{s-q+1}\wedge T^-_{k-s}),$$ 
and Lemma \ref{defwedgepull} assures us that this is finite.

Using Lemma \ref{decpushpull}, performing the same Abel transform and using again Lemma \ref{pushpull}, we obtain similarly that: 
\begin{align}\label{key2}
\left(\sum_{n=0}^{N} \frac{1}{d^n}\right)\U_{L(\omega^q)} (\omega^{s-q+1}\wedge T^-_{k-s}) &= \left(\sum_{n=1}^{N} -d^{-n+1}\right) \U_{L(\omega^{s-q+1})}(L(\omega^q)\wedge T^-_{k-s}) \\
 &+\U_{L^N(\omega^{s-q+1})}(L^{N+1}(\omega^q)\wedge T^-_{k-s}) \nonumber \\
 &+ \sum_{n=0}^{N}  \frac{1}{d^n}\U_{L(\omega^q)}(\Lambda^n(\omega^{s-q+1}\wedge T^-_{k-s})) \nonumber
\end{align} 
We have again that $\U_{L^N(\omega^{s-q+1})}(L^{N+1}(\omega^q)\wedge T^-_{k-s})$ is negative since $\U_{L^N(\omega^{s-q+1})}$ is negative.
Thus, letting $N$ go to $\infty$:
\begin{align*}
\left(\sum_{n\geq 0} d^{-n} \right)\U_{L(\omega^q)} (\omega^{s-q+1}\wedge T^-_{k-s}) &+ \left(\sum_{n\geq1}d^{-n+1}\right) \U_{L(\omega^{s-q+1})}(L(\omega^q)\wedge T^-_{k-s})  \\
                                                                        & \qquad \leq \sum_{n\geq 0}  \frac{1}{d^n}\U_{L(\omega^q)}(\Lambda_n(\omega^{s-q+1}\wedge T^-_{k-s})).
\end{align*}
Again, by Proposition \ref{cor_decreasing_sqp}, we have that the sequence of super-potential of $L^n(\omega^q)$ is decreasing thus to have the convergence in the Hartogs' sense, it is sufficient to have the convergence at one point. We recognize  by (\ref{Greenadm}) that the right-hand side gives in fact the convergence at the point $\omega^{s-q+1}\wedge T^-_{k-s}$ (again $\omega^{s-q+1}\wedge T^-_{k-s}$ is $(f^n)_*$-admissible so $\Lambda_n(\omega^{s-q+1}\wedge T^-_{k-s})=\Lambda^n(\omega^{s-q+1}\wedge T^-_{k-s})$). So we have that $\U_{T_q^+}(\omega^{s-q+1}\wedge T^-_{k-s})$ is finite and $L^n(\omega^q)$ converges to $T^+_q$ in the Hartogs' sense. \\
 
In (\ref{key2}) every term converges. In particular, 
 $$(\U_{L^N(\omega^{s-q+1})}(L^{N+1}(\omega^q)\wedge T^-_{k-s}))_N $$ 
converges to a finite value. Using Lemma \ref{byparts}, we have the identity:
 \begin{align*}
\U_{L^N(\omega^{s-q+1})}(L^{N+1}(\omega^q)\wedge T^-_{k-s})&= \U_{L^{N+1}(\omega^q)}( L^N(\omega^{s-q+1}) \wedge T^-_{k-s})  \\
 &-  \U_{L^{N+1}(\omega^q)}(\omega^{s-q+1} \wedge T^-_{k-s})   \\
 &+\U_{L^N(\omega^{s-q+1})}(\omega^q \wedge T^-_{k-s}). 
 \end{align*}
As above, every term converges to a finite value. In particular, that means that $ \U_{T_{s-q+1}^+}(\omega^q\wedge T^-_{k-s})$ is finite (which is already known by exchanging the role of $q$ and $s-q+1$). Finally the convergence of $ \U_{L^{N+1}(\omega^q)}( L^N(\omega^{s-q+1}) \wedge T^-_{k-s}) $ implies that $\U(T_q^+, T^+_{s-q+1}\wedge T^-_{k-s})$ is finite. \hfill $\Box$ \hfill

We prove that $T^+_q$ is invariant.
\begin{lemme}\label{q-invariant}
For $1\leq q\leq s$, the Green current $T^+_q$ is $f^*$-admissible and satisfies $T^+_q=L(T^+_q)$. Furthermore, $T^+_q$ is the most $H$-regular current which is $f^*$-invariant in $\C_q$. In particular, $T_q^+$ is extremal in the set of $f^*$-invariant currents of $\mathcal{C}_q$.
\end{lemme}
\emph{Proof.} For $q=s$, this is Theorem \ref{invariance}. So take $q< s$. We have that $L^n(\omega^{q})$ converges in the Hartogs' sense to $T^+_q$. So this means that at the point $\omega^{k-q+1}$ we have the convergence of the series:
$$  \sum_{n \geq 0}  d^{-n} \U_{L(\omega^q)}(\Lambda^n(\omega^{k-q+1})).  $$
In particular, dropping the first term and multiplying by $d$, we have the convergence of the series:
$$  \sum_{n \geq 0} d^{-n}\U_{L(\omega^q)}(\Lambda^n( \Lambda(\omega^{k-q+1}))).$$ 
We recognize $\U_{T^+_q}(\Lambda(\omega^{k-q+1}))>-\infty$ hence $T^+_q$ is $f^*$-admissible 

By Theorem \ref{continuiteL}, we see that $L(L^n(\omega^q))$ converges to $L(T^+_q)$ and to $T^+_q$. So we have proved the first part of the lemma. The rest is exactly as in Theorem \ref{invariance}. \hfill $\Box$ \hfill \\

Now, we also want to consider the intersection $T^+_q \wedge T^-_{k-s}$ for $q < s$. First, we have that these intersections are well defined elements of $\C_{k-s+q}$ from Proposition \ref{extension} ($T^+_q$ and $T^-_{k-s}$ are wedgeable). Furthermore, it is $f_*$-admissible since we have by Proposition \ref{extension} that 
$$ \U_{T^+_{s-q+1}}(T^+_{q}\wedge T^-_{k-s})>-\infty$$ 
and since $L(\omega^{s-q+1})$ is more H-regular than $T^+_{s-q+1}=L(T^+_{s-q+1})$, we see that:
$$ \U_{L(\omega^{s-q+1})}(T^+_{q}\wedge T^-_{k-s})>-\infty,$$ 
which means that $T^+_{q}\wedge T^-_{k-s}$ is $f_*$-admissible since by symmetry of the super-potential, its super-potentials are finite at the point $L(\omega^{s-q+1})$.

 Using the same argument than in the proof of the invariance of the measure $\mu$, one has:
\begin{proposition}
The current $T^+_q \wedge T^-_{k-s} \in \C_{k-s+q}$ is $f_*$-invariant, that is:
$$\Lambda(T^+_q \wedge T^-_{k-s})=T^+_q \wedge T^-_{k-s}.$$
\end{proposition} 
\emph{Proof.} This follows from the fact that $L^n(\omega^q)$ and $\Lambda^m(\omega^{k-s})$
 converge in the Hartogs' sense to $T^+_q$ and $T^-_{k-s}$ which are wedgeable and we use Proposition \ref{pushpull}. \hfill $\Box$ \hfill \\

Now, we use the same arguments than in the proof of Theorem \ref{constructionmeasure}, but we replace $T^-_{k-s}$ by $T^+_q \wedge T^-_{k-s}$. Our purpose is to show that $T^+_s=(T^+)^s$. We need the following lemma first:
\begin{lemme}\label{lem_wedge_green}
Let $q_1\geq 1$ and $q_2$ such that $q_1+q_2=s-q+1$. Then the current $T^+_{q_2}$ and $T^+_q \wedge T^-_{k-s}$ are wedgeable and we have that a super-potential $\U_{T^+_{q_1}}$ of $T^+_{q_1}$ satisfies:
$$\U_{T^+_{q_1}}(T^+_{q_2}\wedge T^+_q \wedge T^-_{k-s}) >-\infty.$$  
\end{lemme}
The proof is essentially the same as the one of Theorem \ref{constructionmeasure}. We need the equivalent of Lemma \ref{defwedgepull} first:
\begin{lemme}
Let $q_1 \geq 1$ and $q_2$ such that $q_1+q_2=s-q+1$ and let $n\in \mathbb{N}$. Then the currents $L^n(\omega^{q_2})$ and  $T^+_q \wedge T^-_{k-s}$ are wedgeable. Furthermore, for $n' \in \mathbb{N}$:
$$ \U_{L^{n'}(\omega^{q_1})}(L^n(\omega^{q_2}) \wedge T^+_q \wedge T^-_{k-s}) >-\infty.$$
\end{lemme}
\emph{Proof.} We can assume that $ q_2 \geq 1$ (else it is just Proposition \ref{extension}). The super-potentials of the current $T^+_q \wedge T^-_{k-s}$ are finite at $L(\omega^{q_1+q_2})=L(\omega^{q_1})\wedge L(\omega^{q_2})$ which is less H-regular than $\omega^{q_1}\wedge L(\omega^{q_2})$. Hence the super-potentias of the current $T^+_q \wedge T^-_{k-s}$ are finite at $\omega^{q_1}\wedge L(\omega^{q_2})$. This means that the currents $L(\omega^{q_2})$ and  $T^+_q \wedge T^-_{k-s}$ are wedgeable. 

On the other hand, $\U_{L(\omega^{q_1+q_2})}(T^+_q \wedge T^-_{k-s})$ is finite. We can use Lemma \ref{dec_wedge} and we have that:
$$\U_{L(\omega^{q_1+q_2})}(T^+_q \wedge T^-_{k-s})= \U_{L(\omega^{q_2})}(\omega^{q_1}\wedge T^+_q \wedge T^-_{k-s}) +\U_{L(\omega^{q_1})}(L(\omega^{q_2})\wedge T^+_q \wedge T^-_{k-s}). $$
Again taking the difference with $\U_{L(\omega^{q_2})}(\omega^{q_1}\wedge T^+_q \wedge T^-_{k-s})$, we have that:
$$\U_{L(\omega^{q_1})}(L(\omega^{q_2})\wedge T^+_q \wedge T^-_{k-s})$$
is well defined in the sense of super-potentials and is finite. We have proved the lemma for $n=n'=1$. The rest follows as in Lemma \ref{defwedgepull}. \hfill $\Box$ \hfill \\

\noindent \emph{Proof of Lemma \ref{lem_wedge_green}.} We replace $T^-_{k-s}$  by $T^+_q \wedge T^-_{k-s}$ and we do the same computations. Lemma \ref{pushpull}, \ref{decpushpull} and \ref{byparts} still apply for $T^+_q \wedge T^-_{k-s}$. \hfill $\Box$ \hfill \\

We can now prove the following corollary. Observe that if a sequence $S_n$ converges in the Hartogs' sense to $S$ and a sequence $R_n$ converges in the Hartogs' sense to $R$ with $R_n \wedge S_n$ wedgeable converging in the Hartogs' sense to a current $T$, we cannot claim a priori that $S$ and $R$ are wedgeable and that $T=R\wedge S$. But if $S$ and $R$ are wedgeable, then we do have $T=R\wedge S$.
\begin{corollaire}
The current $T^+_s$ satisfies $T^+_s=(T^+)^s$. Consequently, one has $\mu=(T^+)^s \wedge (T^-)^{k-s}$ where $T^\pm$ is the Green current of order 1 of $f^{\pm}$. 
\end{corollaire} 
\emph{Proof.} Applying the previous lemma to $q=1$, $q_1=1$ and $q_2=s-1$  gives that 
$$\U_{T^+}(T^+_{s-1}\wedge T^+ \wedge T^-_{k-s}) >-\infty.$$
Since $\omega^{k-s+1}$ is more  H-regular than $T^+ \wedge T^-_{k-s}$ that implies that:
$$\U_{T^+}(T^+_{s-1}\wedge \omega^{k-s+1}) >-\infty.$$
In particular that $T^+$ and $T^+_{s-1}$ are wedgeable. Since $L^n(\omega)$ and $L^n(\omega^{s-1})$ converges in the Hartogs' sense to $T^+$ and $T^+_{s-1}$ and $L^n(\omega^s)$ converges in the Hartogs' sense to $T^+_s$, Proposition  \ref{prop_cv_wedge_hartogs} implies that $T^+\wedge T^+_{s-1}=T^+_s$. An easy induction gives the result for $T^+_s$ and the result follows for $\mu$. \hfill $\Box$ \hfill \\ 

\begin{Remark} \rm We do not know how to prove the previous 
result without constructing $T^-_{k-s}$  first.
In the case where $f$ satisfies Hypothesis \ref{distancefort}, 
the result was proved directly (see Theorem \ref{convergencecurrent}). 
This illustrate the difference between Hypotheses \ref{distancefort} 
and \ref{distance}. For a map satisfying Hypotheses \ref{distancefort}, 
we have that the potential of $T^+$ is finite at every point of $I^-$, 
if it only satisfies Hypotheses \ref{distance} we can only say that 
$\int_{I^-} U_{T^+} \omega^{s-1}$ is finite since 
$T^+ \wedge \omega^{s-1}$ is more H-regular than $T^+_s=(T^+)^s$. 
\end{Remark}

Now, we improve the previous results and we show that the measure $\mu$ gives no mass to pluripolar sets (hence analytic sets). The proof relies on a space of test functions introduced by Dinh and Sibony in \cite{DS4} and studied by the second author in \cite{moi2}. 
Recall that the space $W^{1,2}(\P^k)$ is the set of functions in $L^2$ whose differential in the sense of currents can be represented by a form in $L^2$. The space $W^*(\P^k)$ is the set of functions $\varphi$ in $W^{1,2}(\P^k)$ such that there exists a positive closed current $S_\varphi$ of bidegree $(1,1)$ satisfying:
	\begin{equation}
	 \label{definitionW^*}
	d \varphi \wedge d^c \varphi \leq S_\varphi.
	\end{equation}
For $\varphi \in W^{*}$, we define the norm: 
$$ \|\varphi \|^2_{*}=  \|\varphi \|_{L^2}^2 + \inf \Big\{ m(S), \  S \ \text{closed, satisfying (\ref{definitionW^*})} \Big\}.$$ 
 Let $\psi$ be a qpsh function in $W^*(\P^k)$. Consider the regularization $\psi_n$ of $\psi$ obtained through an approximation of the identity in $\text{PGL}(k+1,\mathbb{C})$. Let $S$ be minimal in (\ref{definitionW^*}) for $\psi$ and let $S_n$ be the smooth regularization of $S$ obtained through the same approximation of the identity. Using Lemma 5 in \cite{moi2}, we have
\begin{itemize}
 \item $\psi_n$ ``decreases'' to $\psi$. 
 \item $d\psi_n\wedge d^c\psi_n \leq S_n$, and $m(S_n) \to m(S)$ thus  $\lim \| \psi_n \|_* = \| \psi \|_*$. 
\end{itemize}
If $\varphi$ is a qpsh function in $\P^k$ with $\varphi \leq -2$, then $\psi:=-\log-\varphi$ is in $W^*(\P^k)$, thus for every pluripolar set in $\P^k$ there exists a qpsh function in $W^*(\P^k)$ equal to $-\infty$ on that set (see Example 1 p. 253 in \cite{moi2}). In particular, if the qpsh functions in $W^*(\P^k)$ are integrable for a measure, the measure cannot give mass to the pluripolar sets. We can now state the theorem:
\begin{theorem}
The measure $\mu$ gives no mass to pluripolar sets (hence analytic sets). More precisely, there exists $C>0$ such that for $\psi<0$ a qpsh function in $W^*(\P^k)$, we have that:
$$ |\mu(\psi)| \leq C \| \psi \|_*.$$
\end{theorem}
\emph{Proof.} Let $\psi$ and $\psi_n$ be as above. Recall that $G$ is the potential of $T^+$. Let $T_m^+$ and $T^-_m$ be sequence of smooth currents in $\C_1$ converging to $T^+$ and $T^-$ in the Hartogs' sense. Then $\mu_m=(T_m^+)^s \wedge (T_m^-)^{k-s}$ converges to $\mu$ in the Hartogs' sense by Proposition \ref{prop_cv_wedge_hartogs}. Let $G_m$ be the associated potential of $T^+_m$.
Using Stokes' formula and Cauchy-Schwarz inequality, we have that:
\begin{align*}
\left| \int \psi_n d\mu_m \right|& = \left|\int \psi_n (dd^cG_m+\omega)\wedge (T_m^+)^{s-1} \wedge (T_m^-)^{k-s} \right| \\
                    &\leq \left|\int d\psi_n \wedge d^cG_m\wedge (T_m^+)^{s-1} \wedge (T_m^-)^{k-s} \right| \\
                    &+\left|\int \psi_n \omega\wedge (T_m^+)^{s-1} \wedge (T_m^-)^{k-s}\right|  \\
                    &\leq \left(\int d\psi_n \wedge d^c\psi_n \wedge (T_m^+)^{s-1} \wedge (T_m^-)^{k-s} \right)^{\frac{1}{2}}\\
                    &\times \left(\int dG_m\wedge d^cG_m \wedge (T_m^+)^{s-1} \wedge (T_m^-)^{k-s} \right)^{\frac{1}{2}}\\
                    &+\left|\int \psi_n \omega\wedge (T_m^+)^{s-1} \wedge (T_m^-)^{k-s}\right|.
\end{align*}
Let $S_n$ be the positive closed current of bidegree $(1,1)$ such that $d\psi_n \wedge d^c\psi_n\leq S_n$. Using again Stokes' formula for the second term of the product yields:
\begin{align*}
\left| \int \psi_n d\mu_m \right|  &\leq \left(\int S_n \wedge (T_m^+)^{s-1} \wedge (T_m^-)^{k-s} \right)^{\frac{1}{2}}\\
                    							 &\times \left(-\int G_m dd^cG_m \wedge (T_m^+)^{s-1} \wedge (T_m^-)^{k-s} \right)^{\frac{1}{2}}\\
                   								 &+\left|\int \psi_n \omega\wedge (T_m^+)^{s-1}\wedge (T_m^-)^{k-s}\right|. 
\end{align*}
We let $m$ go to $\infty$, we have that $\left| \int \psi_n d\mu_m \right|$ converges to $\left| \int \psi_n d\mu \right|$,  
$$\left(\int S_n \wedge (T_m^+)^{s-1} \wedge (T_m^-)^{k-s} \right)$$
 converges to $\left(\int S_n \wedge T^+_{s-1} \wedge T^-_{k-s} \right)$, and $\left|\int \psi_n \omega\wedge (T_m^+)^{s-1} \wedge (T_m^-)^{k-s}\right|$ converges to $\left|\int \psi_n \omega\wedge T^+_{s-1} \wedge T^-_{k-s}\right|$. The term:
\begin{align*}
\int G_m dd^cG_m \wedge (T_m^+)^{s-1} \wedge (T_m^-)^{k-s}= &\int G_m (T_m^+)^{s} \wedge (T_m^-)^{k-s} \\
                                                            &-\int G_m \omega \wedge (T_m^+)^{s-1} \wedge (T_m^-)^{k-s}
\end{align*}
can be rewritten as:
$$\U_1({T^+_m},(T_m^+)^{s} \wedge (T_m^-)^{k-s})-\U_1({T^+_m},\omega\wedge (T_m^+)^{s-1} \wedge (T_m^-)^{k-s}) $$
which by Hartogs' convergence goes with $m$ to:
 $$\U_1({T^+},\mu)-\U_1({T^+},\omega\wedge T^+_{s-1} \wedge T^-_{k-s})$$
which is finite by Theorem \ref{constructionmeasure} (observe that $\omega\wedge T^+_{s-1} \wedge T^-_{k-s}$ is more H-regular than $\mu$). So we have that:
\begin{align*}
\left| \int \psi_n d\mu \right| & \leq C \left(\int S_n \wedge T_{s-1}^+ \wedge T^-_{k-s} \right)^{\frac{1}{2}}
                                  +\left|\int \psi_n \omega\wedge T^+_{s-1}\wedge T^-_{k-s}\right|,
\end{align*}
where $C^2=\U_1({T^+},\mu)-\U_1({T^+},\omega\wedge T^+_{s-1} \wedge T^-_{k-s})$ is a constant. 

The term $\left(\int S_n \wedge T^+_{s-1} \wedge T^-_{k-s}
\right)^{\frac{1}{2}}$ is controlled by $\|\psi\|_* +
\varepsilon$ for $n$ large enough because $S_n$ is smooth so that wedge-product is well defined and the mass can be computed in cohomology. 

We use an induction to control in the same way the term $\left|\int \psi_n \omega\wedge T^+_{s-1} \wedge T^-_{k-s}\right|$ (at the last step, we have a term in $\int -\psi_n\omega^k$). Since for $n$ large enough we have $\|\psi_n\|_* \leq \| \psi \|_*+\varepsilon$ ($\varepsilon>0$), we have proved that:
 $$\left| \int \psi_n d\mu \right| \leq C (\| \psi \|_* +\varepsilon).$$
By monotone convergence and letting $\varepsilon \to 0$, we have the theorem. \hfill $\Box$ \hfill

\subsection{Mixing, entropy and hyperbolicity of $\mu$}
We now prove that $\mu$ is mixing, that is $\lim_{n\to \infty} \mu(\varphi \psi\circ f^n)=\mu(\varphi)\mu(\psi)$ for $\varphi$ and $\psi$ smooth functions on $\P^k$. Here the function $\psi\circ f^n$ is not smooth, so by definition $\mu(\varphi \psi\circ f^n)$ is the integral of  $\varphi \psi\circ f^n$ on $\P^k \setminus I(f^n)$ for the measure $\mu$ which gives no mass to $I(f^n)$. Recall that $I(f^n) \subset \C(f^n)$.

We need the classical lemma (\cite{Sib} and \cite{Gu}):
\begin{lemme}
Let $\psi$ be smooth function on $\P^k$, then the sequence of currents  $(\psi\circ f^n T^+_s)_n$ converges to $cT^+_s$ where $c=\mu(\psi)$. Moreover, we have that $\|d( \psi\circ f^n T^+_s) \|$ and $\|dd^c( \psi\circ f^n T^+_s) \|$ go to zero.
\end{lemme}  
\emph{Proof.} The norm  $\|d( \psi\circ f^n T^+_s) \|$ is the operator's norm on the space of smooth forms.

We can assume that $0\leq\psi\leq 1$. Then, the sequence $(\psi\circ
 f^n T^+_s)_n$ is bounded so we can extract a subsequence converging
 in the sense of currents to $S \geq 0$ which satisfies $S\leq
 T^+_s$. In order to show that $S$ is closed and that $\|d( \psi\circ
 f^n T^+_s) \|\to0$, we only need to show that for every smooth
 $(0,1)$-form $\theta$ we have that $| \langle \psi \circ f^n T^+_s,
 \partial (\theta \wedge \omega^{k-s-1}) \rangle |$ goes to $0$ uniformly on $\theta$ (see \cite{DT3} p. 3 for details). In other words, we want to compute the limit of:
$$ \int_{\P^k \setminus I(f^n) } \psi \circ f^n T^+_s \wedge \partial(\theta) \wedge \omega^{k-s-1}. $$
We are going to use the technics of \cite{Sibony2}. Let $v<0$ be a qpsh function equal to $-\infty$ on $ \C(f^n) $ and smooth outside $ \C(f^n) $. Let $\max'$ be a smooth convex increasing function approximating the function $\max^+:=\max(x,0)$ such that its derivative is less than 1. Let $v_j=\max'(v/j+1)$. Then $(v_j)$ is an increasing sequence of smooth qpsh ($i\partial \bar{\partial} v_j +\omega \geq 0$) functions with $0\leq v_j \leq 1$ converging uniformly to $1$ on the compact sets of $\P^k \backslash  \C(f^n)$ and equal to $0$ on some neighborhood of $\C(f^n) $. Let $\alpha:[0,1]\to [0,1]$ be a smooth function equal to $0$ in $[0,1/3]$ and to $1$ in $[2/3,1]$. Then the sequence of functions $v'_j := \alpha\circ v_j$ is equal to $1$  on the compact sets of $\P^k \backslash  \C(f^n)$ for $j$ large enough and is equal to $0$ on some neighborhood of $ \C(f^n)$.  

Since $T^+_s$ gives no mass to $\C(f^n)$, the previous quantity is the limit when $j$ goes to $\infty$ of:
 $$ \langle v'_j \psi \circ f^n T_s^+, \partial(\theta) \wedge \omega^{k-s-1} \rangle. $$
 By Stokes' formula, it is equal to:
 \begin{align*}
 -\langle  v'_j\partial(\psi \circ f^n)\wedge T^+_s, \theta \wedge
  \omega^{k-s-1}\rangle -\langle   \psi \circ f^n\partial(v'_j) \wedge
  T^+_s , \theta \wedge \omega^{k-s-1} \rangle .
\end{align*}
We apply Cauchy-Schwarz inequality for the first term of the sum, we bound the absolute
 value of the first term of the previous quantity by:
\begin{align*}
\langle (v'_j)^2 i\partial \psi\circ f^n\wedge \bar{\partial}\psi\circ f^n \wedge  T^+_s, \omega^{k-s-1} \rangle^\frac{1}{2} \times \langle i\theta\wedge \overline{\theta} \wedge   T^+_s, \omega^{k-s-1}\rangle^\frac{1}{2} .
\end{align*}
The second term of the product is bounded and does not depend on $j$
and $n$ (uniformly in $\| \theta \|$). For the first term, observe that :
$$ i\partial \psi\circ f^n\wedge \bar{\partial}\psi\circ f^n \wedge  T^+_s =  d^{-sn} (f^n)^* (i\partial \psi \wedge \bar{\partial}\psi \wedge  T^+_s)$$ 
in the integral since $f^n$ is smooth on the support of $v'_j$ and one can multiply a positive closed current by a smooth form and take the pull-back by a smooth function. So, assuming that $i\partial \psi \wedge \bar{\partial}\psi\leq \omega$, we have that the first term is less than:
$$  \langle d^{-sn}  (f^n)^*( \omega \wedge  T^+_s) , \omega^{k-s-1}
\rangle^\frac{1}{2} = (\delta^{-(k-s)n} \delta^{(k-s-1)n})^{\frac{1}{2}}= \delta^{-n/2}$$
which goes to $0$ when $n$ goes to $\infty$ independtly of $j$.\\

Now we have to control the term:
$$ \langle  \partial(v'_j) \psi \circ f^n\wedge T^+_s, \theta \wedge \omega^{k-s-1} \rangle. $$
We have that  $\partial(v'_j)= \alpha'(v_j)\partial v_j$ and observe that the sequence of functions $(\alpha'(v_j))$ is bounded and converges uniformly to $0$ on the compact sets of $\P^k \backslash \C(f^n)$. 
We apply Cauchy-Schwarz inequality and we get that:
\begin{align*}
 \langle  \partial(v'_j)\psi \circ f^n\wedge T^+_s, \theta \wedge \omega^{k-s-1} \rangle^2 \leq \\
  \langle  i\partial(v_j) \wedge \bar{\partial}(v_j) \wedge T^+_s,
  \omega^{k-s-1} \rangle   \langle i (\alpha'(v_j))^2\theta \wedge
  \overline{\theta} \wedge T^+_s, \omega^{k-s-1} \rangle .
\end{align*}
The first term of the product is equal by Stokes' formula to:
 \begin{align*}
  \langle  -v_j \wedge i\partial\bar{\partial}(v_j) \wedge T^+_s, \omega^{k-s-1} \rangle
 \end{align*}
 Since $0\leq v_j \leq 1$ and $ i\partial\bar{\partial}v_j +\omega \geq 0$, it is less than:
 \begin{align*}
  \langle v_j \omega  \wedge T^+_s, \omega^{k-s-1} \rangle \leq   \langle  \omega  \wedge T^+_s, \omega^{k-s-1} \rangle 
 \end{align*}
which is bounded independently of $n$ and $j$. The second term of the product goes to $0$ when $j \to \infty$ uniformly on $\theta$ by dominated convergence since $T^+_s$ gives no mass to $ \C(f^n)$. 
So letting $ j \to \infty$ first, we see that:
$$ \langle \psi \circ f^n T^+_s , \partial(\theta) \wedge \omega^{k-s-1} \rangle $$
goes to $0$ when $n \to \infty$ uniformly on $\|\theta\|$.\\

By Theorem \ref{extremality}, this shows that $S=cT^+_s$. To compute $c$, consider $\langle \psi\circ f^n T^+_s,\omega^{k-s} \rangle$. It is equal to $\langle  T^+_s \wedge\Lambda^n(\omega^{k-s}),\psi \rangle$: replace $T^+_s$ by a smooth approximation $T^+_m$, then $\psi\circ f^n L^n(T^+_m)=d^{-ns}(f^n)^*(\psi T^+_m)$, so
  $$\langle \psi\circ f^n L^n(T^+_m),\omega^{k-s} \rangle=\langle  T^+_m \wedge\Lambda^n(\omega^{k-s}),\psi \rangle$$
 and let $m$ go to $\infty$.  So we have $\langle \psi\circ f^n T^+_s,\omega^{k-s} \rangle=\langle  T^+_s \wedge\Lambda^n(\omega^{k-s}),\psi \rangle$ because $T^+_s\wedge \omega^{k-s}$ gives no mass to $I(f^n)$.
 
 By Theorem \ref{constructionmeasure}, we have that $T^+_s \wedge\Lambda^n(\omega^{k-s})$ converges (in the Hartogs' sense hence in the sense of currents) to $\mu$ which means that $c=\mu(\psi)$. In particular, $c$ does not depend on the choice of $S$ and the first part of the lemma follows. \\

Now we show  that $\|dd^c ( \psi\circ f^n T^+_s) \|$ goes to zero. Let $\Theta$ be a test form of bidegree $(k-s-1,k-s-1)$. Again, we consider a smooth approximation of $T^+_s$ that we denote $T^+_m$. Using the fact that $(\psi\circ f^n) L^n(T^+_m)=d^{-sn}(f^n)^*(\psi T^+_m)$, we compute:
\begin{align*}
\langle \psi\circ f^n L^n(T^+_m),dd^c \Theta \rangle &= \langle  d^{-sn}(f^n)^*(\psi T^+_m),dd^c \Theta) \rangle \\
                                                     &= \langle d^{-sn}(f^n)^*(dd^c(\psi) \wedge T^+_m), \Theta \rangle. 
\end{align*}
 Writing $\Theta =\Theta^+-\Theta ^-$ we can assume that $\Theta$ is
 positive (so $\Theta \leq B\omega^{k-s-1}$ with $B>0$ large enough
 which depends only on $\| \Theta\|$). Let $A>0$ be such that $-A\omega\leq dd^c\psi\leq A\omega$. It is sufficient to control:
 $$ \langle d^{-sn}(f^n)^*(\omega \wedge T^+_m), \omega^{k-s-1} \rangle.$$
 We recognize that this is equal by definition to
 $d_{s+1}^n d^{-sn}= \delta^{-n}$. We let $m$ go to $\infty$
 and we have that $\langle dd^c(\psi\circ f^n \wedge T^+_s),\Theta
 \rangle $ goes to $0$ with $n$ uniformly on $\| \Theta\|$.  \hfill $\Box$ \hfill

\begin{theorem}\label{mixing}
The measure $\mu$ is mixing.
\end{theorem} 
\emph{Proof.} Let $\psi$ and $\varphi$ be real smooth functions on $\P^k$. 
We can assume without loss of generality that $0\leq \psi,\varphi\leq 1$. 
Then for $S$ in $\C_{k-s}$ smooth, we have by the above lemma that:
$$ \langle  (\varphi \psi\circ f^n ) T^+_s, S \rangle $$
 converges  to: 
$$  \mu(\psi) \langle  \varphi T^+_s, S\rangle. $$ 
We consider a sequence $(T^-_m)$ of smooth currents in $\C_{1}$ converging in 
the Hartogs' sense to $T^-$ (the Green current of order 1 of $f^{-1}$). Then 
let $m=(m_1,m_2,\dots, m_{k-s})$ and $m'=(m'_1,m'_2,\dots, m'_{k-s})$ in $\mathbb{N}^{k-s}$. 
We have that $T^-_{m_1}\wedge \dots  \wedge T^-_{m_{k-s}}$ converges 
to $T^-_{k-s}$ in the Hartogs' sense when the $m_i$ go to $\infty$. We decompose: 
$$T^-_{m_1}\wedge \dots  \wedge T^-_{m_{k-s}}- T^-_{m'_1}\wedge \dots  \wedge T^-_{m'_{k-s}}$$
as:
\begin{align*}
(T^-_{m_1}- T^-_{m'_1})\wedge T^-_{m_2}\wedge \dots  \wedge T^-_{m_{k-s}} + \\
T^-_{m'_1}\wedge (T^-_{m_2}-T^-_{m'_2}) \wedge \dots \wedge T^-_{m_{k-s}}+ \\
\dots \\
T^-_{m'_1}\wedge \dots  \wedge T^-_{m'_{k-s-1}} \wedge (T^-_{m_{k-s}}- T^-_{m'_{k-s}}).
\end{align*} 
As in the previous lemma, let  $(v_j)$ be an increasing sequence of 
smooth qpsh ($i\partial \bar{\partial} v_j +\omega \geq 0$) functions 
with $0\leq v_j \leq 1$ converging uniformly to $1$ on the compact sets 
of $\P^k \backslash \C(f^n)$ and equal to $0$ on some 
neighborhood of $ \C(f^n)$.

 We also define $v'_j:=\alpha \circ v_j$ with 
$\alpha:[0,1]\to [0,1]$ a smooth function equal to $0$ in $[0,1/3]$ and to $1$ in $[2/3,1]$ so that the sequence of functions $v'_j$ is equal to $1$ on the compact sets of $\P^k \backslash   \C(f^n)$ for $j$ large enough and is equal to $0$ on some neighborhood of $ \C(f^n)$.  \\

We consider the quantity $\langle v'_j \varphi \psi\circ f^n  T^+_s, (T^-_{m_1}- T^-_{m'_1})\wedge \dots  \wedge T^-_{m_{k-s}}  \rangle $.
 Write $T^-_i=\omega+dd^c g_i$ where the $g_i$ are decreasing.  By Stokes' formula, we have that:
\begin{align*}
\langle v'_j \varphi \psi\circ f^n  T^+_s, (T^-_{m_1}- T^-_{m'_1})\wedge \dots  \wedge T^-_{m_{k-s}}  \rangle = \\
-\langle  (\varphi \psi\circ f^n dv'_j+v'_j \psi\circ f^n d\varphi + v'_j\varphi d\psi\circ f^n)\wedge T^+_s, d^c(g_{m_1}- g_{m'_1})\wedge \dots  \wedge T^-_{m_{k-s}}  \rangle.
\end{align*}
 Write the last sum $I+II+III$ with obvious notations. Using Cauchy Schwarz inequality for the first term, we have that:
\begin{align*}
 |I|^2 \leq \langle dv_j \wedge d^cv_j \wedge   T^+_s, T^-_{m_{2}}\wedge \dots  \wedge T^-_{m_{k-s}}  \rangle  \times \\
 \langle (\alpha'(v_j))^2d(g_{m_1}- g_{m'_1})\wedge d^c(g_{m_1}- g_{m'_1}) \wedge  T^+_s, T^-_{m_{2}}\wedge \dots  \wedge T^-_{m_{k-s}}  \rangle.
\end{align*}
As in the proof of the previous lemma, we have that this term goes to zero when $j \to \infty$ since $\alpha'(v_j)$ converges uniformly to $0$ on the compact sets of $\P^k \backslash \C(f^n)$.\\ 

Now for $II$, we use Cauchy Schwarz inequality  and we have that:
\begin{align*}
|II|^2 \leq &\langle d\varphi \wedge d^c \varphi \wedge   T^+_s, T^-_{m_{2}}\wedge \dots  \wedge T^-_{m_{k-s}}  \rangle \\
       &\langle d(g_{m_1}- g_{m'_1})\wedge d^c(g_{m_1}- g_{m'_1}) \wedge  T^+_s, T^-_{m_{2}}\wedge \dots  \wedge T^-_{m_{k-s}}  \rangle.
 \end{align*}
 The first term of the product is bounded as it converges to $\int  d\varphi \wedge d^c \varphi \wedge T^+_s \wedge (T^-)^{s-1}$. By Stokes, we recognize that the second term is equal to 
\begin{align*} 
\langle -(g_{m_1}- g_{m'_1})\wedge dd^c(g_{m_1}- g_{m'_1}) \wedge  T^+_s, T^-_{m_{2}}\wedge \dots  \wedge T^-_{m_{k-s}}  \rangle &= \\
\langle -(g_{m_1}- g_{m'_1})\wedge (T^-_{m_1}- T^-_{m'_1}) \wedge  T^+_s, T^-_{m_{2}}\wedge \dots  \wedge T^-_{m_{k-s}}  \rangle&= \\
\U_{T^-_{m_1}}( T^+_s\wedge T^-_{m'_1}\wedge T^-_{m_{2}}\wedge \dots  \wedge T^-_{m_{k-s}})-\U_{T^-_{m'_1}}( T^+_s\wedge T^-_{m'_1}\wedge T^-_{m_{2}}\wedge \dots  \wedge T^-_{m_{k-s}}) &+ \\
\U_{T^-_{m'_1}}( T^+_s\wedge T^-_{m_1}\wedge T^-_{m_{2}}\wedge \dots  \wedge T^-_{m_{k-s}})-\U_{T^-_{m_1}}( T^+_s\wedge T^-_{m_1}\wedge T^-_{m_{2}}\wedge \dots  \wedge T^-_{m_{k-s}}).
\end{align*}
Observe that this term goes to $0$ when the $m_i, m'_i$ are large enough. Indeed recall that $\U_1(S,T)$ is continuous for the Hartogs' convergence (Lemma \ref{lemma_decreas_sf}), so:
$$\U_{T^-_{m_1}}( T^+_s\wedge T^-_{m'_1}\wedge T^-_{m_{2}}\wedge \dots  \wedge T^-_{m_{k-s}})$$\\
converges to $\U_{T^-}( \mu)$ which is finite and so does the other terms in the majoration of $II$ (the convergence is uniform else we could extract a subsequence which does not converge). \\

Now we bound $III$. Applying Cauchy-Schwarz inequality gives:
\begin{align*}
|III|^2 \leq &\langle |v'_j|^2 d\psi \circ f^n  \wedge d^c \psi \circ f^n \wedge   T^+_s, T^-_{m_{2}}\wedge \dots  \wedge T^-_{m_{k-s}}  \rangle \\
       &\langle d(g_{m_1}- g_{m'_1})\wedge d^c(g_{m_1}- g_{m'_1}) \wedge  T^+_s, T^-_{m_{2}}\wedge \dots  \wedge T^-_{m_{k-s}}  \rangle
 \end{align*}
Observe that the second integral is the same than in the bound of $II$ so it goes to zero. For the first term of the product, we use that $f^n$ is smooth in the support of $v'_j$ and thus $ d\psi\circ f^n \wedge d^c \psi\circ f^n =(f^n)^*(d \psi\wedge d^c \psi)$ in the integral. We can assume that $d \psi\wedge d^c \psi \leq \omega$. Using the invariance of $T^+_s$ and the fact that $v'_j$ is equal to $0$ near $\C(f^n)$, we have that $(f^n)^*(\omega)\wedge T^+_s =d^{-sn} (f^n)^*(\omega \wedge T^+_s)$ in the integral, so the first term in the bound of $III$ is less than:
 $$ \frac{1}{d^{sn}} \langle  (f^n)^* (\omega \wedge  T^+_s), T^-_{m_2} \wedge \dots  \wedge T^-_{m_{k-s}} \rangle.$$
That last term can be computed cohomologically and is equal to $\frac{\delta^{n(k-s-1)}}{d^{ns}}<1$. So as for $II$, we have that $III$ goes to $0$ uniformly in $n$. \\

\noindent Letting $j$ go to $\infty$, we have that 
$$\langle  \varphi \psi\circ f^n  T^+_s, T^-_{m_1}\wedge \dots  \wedge T^-_{m_{k-s}} - T^-_{m'_1}\wedge \dots  \wedge T^-_{m'_{k-s}}\rangle$$
converges uniformly to 0. 
In particular, we can interchange the limit in:
$$ \lim_m \lim_n \langle  (\varphi \psi\circ f^n ) T^+_s, T^-_{m_1}\wedge \dots  \wedge T^-_{m_{k-s}}  \rangle $$
which gives $ \lim_n \mu(\varphi  \psi\circ f^n)= \mu(\varphi)\mu(\psi)$ hence the mixing.  \hfill $\Box$ \hfill 

\vspace{1cm}

We now show that the measure $\mu$ satisfies the hypothesis of Chapter \ref{Henry} and we deduce from that a bound of its entropy. Recall that we denote by $\mu_n$ the sequence of probabilities:
$$\mu_n:= \frac{1}{n} \sum_{i=0}^{n-1} f^i_{*}\left( \frac{(f^n)^* \omega^{s}
  \wedge \omega^{k-s}}{\lambda_s(f^n)} \right).$$
In our case, using Lemma \ref{pushpull}, we can write it as:
$$\mu_n = \frac{1}{n} \sum_{i=0}^{n-1} L^{n-i}( \omega^{s})\wedge \Lambda^{i}(\omega^{k-s}).$$
We consider the hypothesis $(H)$: there exists a subsequence $\mu_{\psi(n)}$ of $\mu_n$ converging to a measure $\mu'$ such that:
$$(H) \mbox{  :   }  \lim_{n \rightarrow + \infty}
\int \log d(x,I)  d \mu_{\psi(n)} (x) = \int \log d(x,I)  d \mu'(x) > - \infty.$$
In here, we do not need to take a subsequence:
\begin{proposition}
The sequence $(\mu_n)$ converges to $\mu$ and satisfies the hypothesis $(H)$.
\end{proposition}
\emph{Proof.} Let $\varphi$ be a smooth test function. Choose $\varepsilon>0$. By Theorem \ref{constructionmeasure}, since $L^{n-i}( \omega^{s})$ and $\Lambda^i(\omega^{k-s})$ converge in the Hartogs' sense, Proposition  \ref{prop_cv_wedge_hartogs} assures us that  $L^{n-i}( \omega^{s})\wedge \Lambda^{i}(\omega^{k-s})$ converges in the Hartogs' sense to $\mu$. So we have for $\sqrt{n}\leq i\leq n-\sqrt{n}$ and $n$ large enough  that$|L^{n-i}( \omega^{s})\wedge \Lambda^{i}(\omega^{k-s})(\varphi)-\mu(\varphi)| \leq \varepsilon$. The fact that $(\mu_n)$ goes to $\mu$ follows since they are $o(n)$ terms for which the estimation does not stand. \\

Now, by Lemma \ref{lemma_wedge_regular2}, we see that there exist constants $A_{i,n}\geq 0$  such that $\U_{L^{n-i}( \omega^{s})\wedge\Lambda^i(\omega^{k-s})} \geq \U_\mu -A_{i,n}$ with $A_{i,n}$  uniformly bounded from above by $C$
and arbitrarily close to zero for $i$ and $n$ large enough. We
consider super-potentials of mean $0$. In particular:

$$\U_{\mu_n} \geq \U_\mu - \frac{1}{n} \sum_{i=0}^{n-1} A_{i,n}.$$

So we have that the sequence $\mu_n$ is more H-regular than $\mu$ for
all $n$. We also have the convergence in the Hartogs' sense to $\mu$
since $\frac{1}{n} \sum_{i=0}^{n-1} A_{i,n}$ goes to $0$ when $n \rightarrow + \infty$.

Thus $\mu_n(G)\to \mu(G)$ which is finite by Theorem \ref{constructionmeasure} where $G$ is a negative potential of the Green current of order $1$ that we denote $T^+$. Since $T^+$ is less H-regular than $L(\omega)$, we have that if $U_{L(\omega)}$ is a quasi-potential of $L(\omega)$ then
$\mu_n(U_{L(\omega)})\to \mu(U_{L(\omega)})$ which is also finite. By Lemma \ref{estimate}, we have that:
$$AU_{L(\omega)}(x) < \log\text{dist}(x, I^+)$$
for $A>0$ large enough. We denote $\varphi:= \log\text{dist}(x, I^+)$. Since $\mu$ gives no mass to $I^+$ that means that $AU_{L(\omega)}(x) \leq \varphi \leq 0$ for $\mu$ a.e point, so we have that $\varphi \in L^1(\mu)$. We have the classical lemma:
\begin{lemme}
Let $\nu_n$ be a sequence of measures converging to $\nu$ in the sense of measures. Then for $v$ an upper semi-continuous function, we have that $$\limsup \nu_n(v) \leq \nu(v).$$ 
\end{lemme}
\emph{Proof.} Recall that an usc function can be written as the limit of a decreasing sequence of continuous functions. So for some small $\alpha>0$ 
we can take $v' \geq v $ a continuous function such that $ \int v' d \nu \leq \int v d \nu +\alpha$ by monotone convergence. In particular:
$$\int v d\nu_n \leq \int v' d\nu_n \to \int v' d\nu \leq \int v d \nu +\alpha.$$  
And the result follows by letting $\alpha \to 0$.  \hfill $\Box$ \hfill \\

\noindent \emph{End of the proof of the proposition.} Now, $\varphi$ is upper semi-continuous, so:
$$ \limsup \mu_n(\varphi) \leq  \mu(\varphi)$$
We also have that $(A+1)U_{L(\omega)}-\varphi$ is  upper semi-continuous (we use the fact that it is equal to $-\infty$ on $I^+$). That and $\mu_n(U_{L(\omega)})\to \mu(U_{L(\omega)})$ give:
$$ \liminf \mu_n(\varphi) \geq  \mu(\varphi).$$
This is exactly the fact that $\mu$ satisfies Hypothesis (H). \hfill $\Box$ \hfill \\

We can now apply Theorem \ref{ENTROPY} to get the proposition:
 \begin{theorem}
The topological entropy of $f$ is greater than $\log d_s=s\log d$. More precisely, the  entropy of $\mu$ is greater than $s\log d$.
 \end{theorem}
On the other hand, the topological entropy is always bounded by $\max_{0\leq s\leq k}{\log d_s}$ (see \cite{DS9} for the projective case and  \cite{DS5} for the Kähler case). So we have the fundamental result:
 \begin{theorem}\label{entropy}
The topological entropy of $f$ is equal to $\log d_s$. Moreover, the entropy of $\mu$ is equal to $s\log d$ so $\mu$ is a measure of maximal entropy.
\end{theorem}

This allows us to use the first author's estimate of the Lyapunov exponents (Corollary 3 in \cite{DT1}). To apply that result, we need to have that $\log(\text{dist}(x,\C^+))$ is integrable with respect to $\mu$. For that observe that the function $U_{L(\omega)}$ is integrable with respect to $\mu$. By invariance, $f_*(U_{L(\omega)})$ is also integrable.  Write $U_{L(\omega)}$ as in Lemma \ref{estimate}:
$$U_{L(\omega)}= d^{-1}\log |F|^2-\log |Z|^2,$$ 
where $f=[P_1:\dots:P_{k+1}]$ and $F=(P_1,\dots,P_{k+1})$. Write $f^{-1}=[Q_1:\dots: Q_{k+1}]$ where the $Q_i$ are homogeneous polynomials of degree $\delta$ and write  $F^{-1}=(Q_1,\dots,Q_{k+1})$. There is of course an abuse of notation since $F \circ F^{-1} \neq Id $ instead, we have that:
$$F \circ F^{-1} = P(z_1,\dots,z_{k+1}) \times \left(z_1,\dots,z_{k+1}\right),$$
where $P$ is an homogeneous polynomial of degree $d \delta-1$ equal to $0$ in $\pi^{-1}(\C^-) $ and $\pi : \mathbb{C}^{k+1} \to \P^k$ is the canonical projection. 
Then, we have that:
$$ f_*(\frac{1}{d}\log |F|^2-\log |Z|^2)= \frac{1}{d}\log |F\circ F^{-1}|^2-\log |F^{-1}(Z)|^2 .$$ 
We recognize $d^{-1} \log |F\circ F^{-1}|^2  -\delta \log |Z|^2+\delta \log |Z|^2-\log |F^{-1}(Z)|^2$.  But 
$$\delta \log |Z|^2-\log |F^{-1}(Z)|^2=-\delta U_{\Lambda(\omega)}$$
 is in $L^1(\mu)$ and by difference, so is $d^{-1} \log |F\circ F^{-1}|^2 -\delta \log |Z|^2$. As in Lemma \ref{estimate}, we then have that $\log \text{dist}( . , \C^-)$ is in $L^1(\mu)$. Similarly, so is $\log \text{dist}( . , \C^+)$ is in $L^1(\mu)$.
\begin{theorem}\label{hyperbolic}
The Lyapunov exponents $\chi_1 \geq \chi_2 \geq \dots \geq \chi_k$ of $\mu$ are well defined and we have the estimates:
\begin{align*}
\chi_1\geq\dots\geq\chi_s\geq \frac{1}{2}\log \frac{d_s}{d_{s-1}}= \frac{1}{2} \log d >0\\
0>-\frac{1}{2} \log \delta=\frac{1}{2}\log \frac{d_{s+1}}{d_{s}} \geq \chi_{s+1}\geq\dots\geq\chi_k. 
\end{align*}
In particular, the measure $\mu$ is hyperbolic. 
 \end{theorem}

\appendix

\chapter{Super-potentials}
\section{Definitions and properties of super-potentials}
We recall here the facts and definitions we use on super-potentials. 
Everything in this section was taken from \cite{DS6} 
so we refer the reader to this paper for proofs and details. \\

Recall that $\C_s$ is the convex compact set of (strongly) positive closed currents
$S$ of bidegree $(s,s)$ on $\P^k$ and of mass 1. To develop the calculus, we have to consider $\C_s$ 
as an infinite dimensional space with special families of currents
that we parametrize by the unit disc $\Delta$ in $\mathbb{C}$. We call these families
{\it special structural discs of currents}. The notion of structural varieties of $\C_s$ was introduced
in \cite{DS11}. In some sense, we consider
$\C_s$ as a space of infinite dimension admitting "complex
subvarieties" of finite dimension. For $S$ in $\C_s$, 
it is always possible to construct a special structural variety 
$\varphi:\Delta\to \C_s$ such that $\varphi(0)=S$ and $\varphi(z)$ is a smooth form for $z\neq 0$.\\

Let $S$ be a current in $\C_s$ with $s\geq 1$. 
If $U$ is a  $(s-1,s-1)$-current such
that $dd^c U=S-\omega^s$, we say that $U$ is a
{\it quasi-potential} of $S$. 
The integral $\langle U,\omega^{k-s+1}\rangle$ is {\it the mean} of $U$. Observe that such quasi-potential is defined up to a $dd^c$-closed current. For $s=1$ such functions are constant a.e., but in the general case, they can be singular currents. Nevertheless, we have the proposition:
\begin{proposition}\label{resolutionquasipotential}
Let $S$ be a current in $\C_s$. 
Then, there is a  negative  quasi-potential $U$ of $S$ depending
linearly on $S$ such that for every $r$ with $1\leq r<k/(k-1)$ and for
$1\leq \rho < 2k/(2k-1)$  
$$\|U\|_{\mathcal{L}^r} \leq c_r\quad \mbox{and} \quad \|dU\|_{\mathcal{L}^\rho}\leq c_\rho$$
for some positive constants $c_r, c_\rho$  independent of $S$. Moreover, $U$ depends
continuously on $S$ with respect to the $\mathcal{L}^r$ topology on $U$ and the
weak topology on $S$.
\end{proposition}
We are going to
introduce a {\it super-potential} associated to $S$. It is an  affine upper
semi-continuous function $\U_S$ defined on $\C_{k-s+1}$
with values in $\mathbb{R}\cup\{-\infty\}$. For $R\in\C_{k-s+1}$ smooth, we define the super-potential of mean $M$ of $S$ by $\U_S(R):=\langle S ,  U_R\rangle$ where $U_R$ is a quasi-potential of $R$ of mean $\langle U_R,\omega^s\rangle=M$. The integral $\langle S ,  U_R\rangle$ does not depend on the choice of $U_R$ with a fixed mean $M$. If $S$ is smooth, we have $\U_S(R)=\langle U_S , R\rangle$ where  $U_S$ is a quasi-potential of $S$ of mean $M$. Now assume that $R$ is not smooth. 
Consider the above special structural variety $\varphi:\Delta\to \C_{k-s+1}$ associated to $R \in \C_{k-s+1}$ and write $R_\theta$ for $\varphi(\theta)$. Then the function $u(\theta):=\U_S(R_\theta)$ defined on $\Delta^*$
can be extended as a quasi-subharmonic function on $\Delta$. Let $(S_\theta)$ and let $(R_\theta)$ be special structural disks associated to $S\in \C_s$ and $R\in C_{k-s+1}$. Then we have the proposition:
\begin{proposition}\label{upper_bound}
The function $\U_S$ can be extended in a unique way to an affine upper semi continuous function on
$\C_{k-s+1}$ with values in $\mathbb{R}\cup\{-\infty\}$, also denoted by $\U_S$, such that
$$\U_S(R)=\lim_{\theta\rightarrow 0}\U_{S_\theta}(R)=\lim_{\theta\rightarrow 0}\U_{S}(R_\theta).$$
In particular, we have 
$$\U_S(R)=\limsup_{R' \rightarrow R} \U_S(R') \mbox{ with } R' \mbox{ smooth}.$$
Moreover, there is a constant $c\geq 0$ independent of $S$ such
that if $\U_S$ is the super-potential of mean $m$ of $S$, then 
$\U_S\leq m+c$ everywhere. 
\end{proposition}
Super-potentials determine the current, more precisely, we have the proposition:
\begin{proposition}\label{prop_unique_sp}
Let $I$ be a compact subset in $\P^k$ with $(2k-2s)$-dimensional
Hausdorff measure $0$.
Let $S$, $S'$ be currents in $\C_s$ and $\U_S$, $\U_{S'}$ be
super-potentials of $S$, $S'$. If $\U_S=\U_{S'}$ on smooth forms in $\C_{k-s+1}$ with compact support in  $\P^k\setminus
  I$, then $S=S'$. 
\end{proposition}
For $I=\varnothing$, this tells us that the values of the super-potential on smooth forms determine uniquely the current.

A crucial notion to prove the convergence of currents is the following: 
\begin{defi} 
Let $(S_n)$ be a sequence in $\C_s$  converging to a current $S$. Let $\U_{S_n}$
(resp. $\U_S$) be the super-potential of mean $M_n$ (resp. $M$) of $S_n$
(resp. $S$).
Assume that $M_n$ converge to $M$. If $\U_{S_n}\geq \U_S$ for every
$n$, we say that \emph {$(S_n)$ converge to $S$ in the Hartogs'
  sense}.  If a current $S'$ in $\C_s$ admits a
super-potential $\U_{S'}$ such that $\U_{S'}\geq \U_S$ we say that
$S'$ is \emph{more H-regular than $S$}.
\end{defi} 
Smooth currents are dense for the Hartog's convergence, more precisely:
\begin{proposition} \label{prop_regularization_sqp}
Let $S \in \C_s$ and let $\U$ be a super-potential of $S$ of mean $M$. There is a sequence
of smooth forms $(S_n)$ in $\C_s$ 
with super-potentials $\U_n$ of mean $M_n$ 
such that
\begin{enumerate}
\item[$\bullet$] $\text{supp}(S_n)$ converge to $\text{supp}(S)$;
\item[$\bullet$] $S_n$ converge to $S$ and $M_n\rightarrow M$;
\item[$\bullet$] $(\U_n)$ decreases to $\U$.
\end{enumerate}
\end{proposition}
We have the following convergence theorem:
\begin{proposition} \label{prop_hartogs_sqp}
Let $(S_n)$ be a sequence in $\C_s$  converging to a current $S$. Let $\U_{S_n}$
(resp. $\U_S$) be the super-potential of mean $M_n$ (resp. $M$) of $S_n$
(resp. $S$). Assume that $M_n$ converge to $M$. Let $\U$ be a continuous function
on a compact subset $K$ of $\C_{k-s+1}$ such that $\U_S<\U$ on $K$. Then, for $n$ large enough we have 
$\U_{S_n}< \U$ on $K$. In particular, we have $\limsup \U_{S_n}\leq
\U_S$ on $\C_{k-s+1}$. Furthermore, if $S_n\rightarrow S$ in the Hartogs' sense, then $\U_{S_n}\rightarrow\U_S$ pointwise.
\end{proposition}
In $\C_{k-s+1}$, they are points which are more ``regular'' than other, namely smooth forms. This is a difference with psh functions. In particular, it is often easier to obtain the convergence at such points:
\begin{proposition} \label{lemma_cv_sp_pointwise}
Let $(S_n)$ be a sequence in $\C_s$ and $\U_{S_n}$
 be super-potentials of mean $M_n$  of $S_n$. Assume that $(\U_{S_n})$
 converges to a finite function $\U$ on smooth forms in $\C_{k-s+1}$. Then, $(M_n)$ converges to a
 constant $M$, $(S_n)$ converges to a current $S$ and $\U$ is equal to
 the super-potential of mean $M$ of $S$ on smooth forms in $\C_{k-s+1}$. 
\end{proposition}  
The following is the main argument to get the convergence of the Green current:
\begin{proposition} \label{cor_decreasing_sqp}
Let $\U_{S_n}$ be super-potentials of mean $M_n$ of $S_n$. Assume
that $\U_{S_n}$ decrease to a function $\U\not=-\infty$. Then, $(S_n)$
converges to a current $S$, $(M_n)$ converges to a constant $M$ and $\U$
is the super-potential of mean $M$ of $S$.
\end{proposition}
In particular, the convergence at one point of the 
super-potentials gives the convergence of the currents in the Hartogs' sense in the case of decreasing super-potentials. \\

An interesting symmetry result is that if  $\U_S$ and $\U_R$ are super-potentials 
of the same mean $M$ of $R$ and $S$ respectively, then $\U_S(R)=\U_R(S)$. \\

There is a notion of \emph{super-polarity} for Borel subsets $E$ of $\C_{k-s+1}$.
This notion does not describe ``small'' sets $E$ but rather how singular are the currents in
$E$. 
\begin{defi}
\rm
We say that $E$ is {\it super-polar} in $\C_{k-s+1}$  if there is a
super-potential $\U_S$ of a current $S$ in $\C_s$ such that $E\subset \{\U_S=-\infty\}$. 
\end{defi}
Denote by $\widetilde E$ the
set of currents $cR+(1-c)R'$ with $R\in \widehat E$, $R'\in \C_{k-s+1}$ and $0<c\leq
1$, and $\widehat E$ the {\it barycentric hull} of $E$, i.e. the set of
currents $\int R d\nu(R)$ where $\nu$ is a probability measure on
$\C_{k-s+1}$ such that $\nu(E)=1$. Then, $\widetilde E$ and $\widehat
E$ are convex.
\begin{proposition} \label{prop_convex_hull}
The following properties are equivalent
\begin{enumerate}
\item $E$ is super-polar in $\C_{k-s+1}$.
\item $\widehat E$ is super-polar in $\C_{k-s+1}$.
\item  $\widetilde E$ is super-polar in $\C_{k-s+1}$.
\end{enumerate}
Moreover, a countable union of super-polar sets is super-polar.
\end{proposition}

One of the purposes of super-potentials is to define the
 wedge product of current (see Section 4 in \cite{DS6}). We define a universal function $\U_s$ on
$\C_s\times\C_{k-s+1}$ by
$$\U_s(S,R):=\U_S(R)=\U_R(S)$$ 
where $\U_S$ and $\U_R$ are super-potentials of mean 0 of $S$ and
$R$. The function $\U_s$ is
 is u.s.c. on $\C_s\times\C_{k-s+1}$. It even enjoys a nice continuity for the Hartogs' convergence: 
 \begin{lemme} \label{lemma_decreas_sf}
Let $(S_n)_{n\geq 0}$ and  $(R_n)_{n\geq 0}$ 
be sequences of currents in $\C_s$ 
and $\C_{k-s+1}$ converging in the Hartogs' sense to $S$ and $R$ respectively.
Then, $\U_s(S_n,R_n)$ converge to $\U_s(S,R)$.
Moreover, if $\U_s(S,R)$ is finite, then
$\U_s(S_n,R_n)$ is finite for every $n$.
\end{lemme}
We have the proposition:
\begin{proposition} \label{prop_hyp_wedge}
Let $s_1\in \mathbb{N}^*$ and $s_2\in \mathbb{N}^*$ with $s_1+s_2 \leq k$. The following conditions are equivalent and are symmetric on $R_1 \in \C_{s_1}$ and
$R_2 \in \C_{s_2}$:
\begin{enumerate}
\item $\U_{s_1}(R_1,R_2\wedge \Omega)$ is finite for at least one smooth
      form $\Omega$ in $\C_{k-s_1-s_2+1}$.
\item $\U_{s_1}(R_1,R_2\wedge \Omega)$ is finite for every smooth
      form $\Omega$ in $\C_{k-s_1-s_2+1}$.
\item There are sequences $(R_{i,n})_{n\geq 0}$ in
      $\C_{s_i}$ converging to 
      $R_i$  and a smooth form $\Omega$ in
      $\C_{k-s_1-s_2+1}$ such that $\U_{s_1}(R_{1,n},R_{2,n}\wedge\Omega)$ is bounded.
\end{enumerate}
\end{proposition}

\begin{defi}\label{def_wedge}
We say that $R_1$ and $R_2$ are {\it wedgeable} if they satisfy the
conditions in Proposition \ref{prop_hyp_wedge}.
\end{defi}
Assume that  $R_1\in \C_{s_1}$  and $R_2\in \C_{s_2}$  are wedgeable.
For every  smooth real form $\varphi$ of bidegree 
$(k-s_1-s_2,k-s_1-s_2)$, write
$dd^c\varphi=c(\Omega^+-\Omega^-)$
where $\Omega^\pm$ are smooth forms in $\C_{k-s_1-s_2+1}$
and $c$ is a positive constant.
We define the wedge-product (or the intersection)
$R_1\wedge R_2$  by its action on the smooth forms by:
\begin{equation} \label{eq_wedge}
\langle R_1\wedge R_2,\varphi\rangle  :=   \langle R_2,\omega^{s_1}\wedge\varphi\rangle +
c\U_{s_1}(R_1,R_2\wedge\Omega^+)-c\U_{s_1}(R_1,R_2\wedge\Omega^-).
\end{equation}
The right hand
side of (\ref{eq_wedge}) is
independent of the choice of $c$, $\Omega^\pm$ and depends
linearly on $\varphi$. Moreover,
$R_1\wedge R_2$ defines
a positive closed $(s_1+s_2,s_1+s_2)$-current of mass $1$ with support in
$\text{supp}(R_1)\cap\text{supp}(R_2)$ 
which depends linearly on each $R_i$ and is symmetric with respect
to the variables. The notion of wedgeability behave well with the notion of H-regularity: 
\begin{lemme} \label{lemma_wedge_regular}
Let $R_i$ and $R_i'$ be currents in $\C_{s_i}$. Assume
  that $R_1$ and $R_2$ are wedgeable.
If $R_i'$ is  more H-regular than $R_i$ for $i=1,2$, then $R_1'$ and $R_2'$ are wedgeable and
$R_1'\wedge R_2'$ is more H-regular than $R_1\wedge R_2$. 
\end{lemme}
We will use the following proposition in the construction of the equilibrium measure: 
\begin{proposition} \label{prop_cv_wedge_hartogs}
Let $R_1$, $R_2$ be wedgeable currents as
  above and $R_{i,n}$ be currents in $\C_{s_i}$ converging to $R_i$
  in the Hartogs' sense. Then, $R_{1,n}$, $R_{2,n}$ are wedgeable
  and $R_{1,n}\wedge R_{2,n}$ converge to $R_1\wedge R_2$ in the
  Hartogs' sense.
\end{proposition}
For several currents (more than 2), the notion of wedgeability is defined by induction: that is $R_1$, $R_2$ and $R_3$ are wedgeable if $R_1$ and $R_2$ are wedgeable and $R_1\wedge R_2$  and $R_3$ are wedgeable. One shows that this definition is in fact symmetric in the $R_i$ and we have Proposition \ref{prop_cv_wedge_hartogs} for several currents.\\

An interesting subcase is when we consider currents $R_1,\ldots, R_l$ such that $R_i$ is of bidegree $(1,1)$ for $i\geq 2$.
 For $2\leq i\leq l$,
there is a quasi-psh function $u_i$ on $\P^k$ such that
$$dd^c u_i=R_i- \omega.$$
\begin{lemme} \label{lemma_wedge_equi}
The currents $R_1,\ldots, R_l$ are wedgeable if and only if 
for every $2\leq i\leq l$, $u_i$ is integrable with respect to the trace measure
      of $R_1\wedge\ldots\wedge R_{i-1}$. 
In particular, the last condition is symmetric with respect to
$R_2,\ldots, R_l$. 
\end{lemme}
If $R_2$ has a quasi-potential integrable with respect to $R_1$, it is classical to
define the wedge-product $R_1\wedge R_2$ by
$$R_1\wedge R_2:=dd^c (u_2 R_1)+ \omega\wedge R_1.$$
One defines $R_1\wedge \ldots\wedge R_l$ by
induction. These two definitions coincide. \\

The other use of super-potentials is to define pull-back and push-forward of current by meromorphic maps (see section 5.1 in \cite{DS6}).  We state the result in the case where $f$ is birational although the results are true in the case where $f$ is just meromorphic. Recall that pull-back and push-forward of a current are defined formally by formulae (\ref{eq_pullback_def}) and (\ref{eq_pushforward_def}) of the previous section: 
\begin{align*} 
f^*(S)&:=(\pi_1)_*\big(\pi_2^*(S)\wedge [\Gamma]\big) \\
f_*(R)&:=(\pi_2)_*\big(\pi_1^*(R)\wedge [\Gamma]\big),
\end{align*}
where $[\Gamma]$ is the current of integration of $\Gamma$. We denote by $I^+:=I(f)$ and $I^-=I'(f)=I(f^{-1})$ the indeterminacy sets of $f$ and $f^{-1}$.

In particular, for a current in $R\in\C_{k-s+1}$ smooth near $I^+$ the push-forward is a well defined positive closed $(k-s+1,k-s+1)$-current and the  mass $\lambda_{s-1}$ of $f_*(R)$ does not depend on $R$. Similarly, for a current $S$ in $S\in\C_{s}$ smooth near $I^-$ the pull-back is a well defined positive closed $(s,s)$-current and the mass of $f^*(S)$ is equal to $\lambda_{s}$. So as above we define for these currents $\Lambda(R)= \lambda_{s-1}^{-1}f_*(R)$ and $L(S)= \lambda_s^{-1}f^*(S)$ (the normalized push-forward and pull-back). 

Using the theory of super-potentials we can extend these definitions to other currents. Namely, we say that a current $S\in\C_s$ is \emph{$f^*$-admissible} if there exists a current $R_0\in \C_{k-s+1}$ which is smooth on a neighborhood of $I^+$ such that the super-potentials of $S$ are finite at $\Lambda(R_0)$. For such $S$, if $(S_n)$ is  a sequence of currents converging in the Hartogs' sense to $S$ then $S_n$ is $f^*$-admissible and $(\lambda_s)^{-1}f^*(S_n)$ converges in the Hartogs' sense to a limit independent on the choice of $(S_n)$ that we denote $(\lambda_s)^{-1}f^*(S)$ (in particular $f^*(S)$ is of mass $\lambda_s$).  
In other words, we have the continuity result:
\begin{theorem}\label{continuiteL}
Let $S$ be an $f^*$ admissible current. Let $S_n$ be a sequence converging to $S$ in the Hartogs' sense, then $S_n$ is $f^*$-admissible and $L(S_n)$ converges in the Hartogs' sense to $L(S)$.
\end{theorem}
We say that $S$ is {\it invariant under $f^*$} or that $S$ is
{\it $f^*$-invariant} if $S$ is
$f^*$-admissible and  $L(S)=S$. 
\begin{proposition} \label{prop_pull_general}
Let $S$ be an $f^*$-admissible current in $\C_s$. Let $\U_S$,
$\U_{L(\omega^s)}$ be super-potentials of $S$ and  $L(\omega^s)$. Then 
$\lambda_{s}^{-1}\lambda_{s-1}\U_S\circ\Lambda+\U_{L(\omega^s)}$ is
equal to a super-potential of $L(S)$ on
$R\in\C_{k-s+1}$, smooth in a neighbourhood of $I^+$. 
\end{proposition}  
Similarly, one define push-forward of currents. We remark that an element in $S\in\C_s$ smooth near $I^-$ is $f^*$-admissible and that the two available definitions of $L(S)$ coincide.

\section{Additional properties}
We state now some properties of the super-potentials that we need. Recall that $f$ is a birational map of $\P^k$ satisfying Hypothesis \ref{distance} and that $s$ is such that $\text{dim}(I^+)=k-s-1$ and $\text{dim}(I^-)=s-1$. 
\begin{lemme}\label{dec_wedge}
Let $S_1 \in \C_{r_1}$ and  $S_2 \in \C_{r_2}$ be wedgeable currents with $r_1+r_2 \leq k$. There exist super-potentials $\U_{S_1\wedge S_2}$, $\U_{S_1}$ and $\U_{S_2}$ of 
$S_1\wedge S_2$, $S_1$ and $S_2$ such that:
$$ \U_{S_1\wedge S_2}(R)=\U_{S_1}(R\wedge S_2)+\U_{S_2}(R\wedge \omega^{r_1})$$
for all $R \in \C_{k-r_1-r_2+1}$ such that $R$ and $S_2$ are wedgeable.  
\end{lemme}
\emph{Proof.} Let $S_{1,n}$ and $S_{2,m}$ be sequence of smooth currents in $\C_{r_1}$ and $\C_{r_2}$ converging to $S_1$ and $S_2$ in the Hartogs' sense. If $U_{S_{1,n}}$ and  $U_{S_{2,m}}$ are smooth quasi-potentials of  $S_{1,n}$ and $S_{2,m}$, then:
$$S_{1,n}\wedge S_{2,m}=dd^c(U_{S_{1,n}}\wedge S_{2,m} +U_{S_{2,m}} \wedge \omega^{r_1}) + \omega^{r_1+r_2}.$$
So,  if $\U_{S_{1,n}\wedge S_{2,m}}$, $\U_{S_{1,n}}$ and $\U_{S_{2,m}}$ are super-potentials of mean $0$ of the currents, we have for any $R$ in $\C_{k-r_1-r_2+1}$:
$$\U_{S_{1,n}\wedge S_{2,m}}(R)= \U_{S_{1,n}}(S_{2,m}\wedge R) +  \U_{S_{2,m}}(\omega^{r_1}\wedge R) -\U_{S_{1,n}}(S_{2,m}\wedge \omega^{k-r_1-r_2+1}).  $$
Now, we take $R$ such that $R$ and $S_2$ are wedgeable. We let $n\to \infty$. By Proposition \ref{prop_cv_wedge_hartogs}, $S_{1,n}\wedge S_{2,m}$ converges in the Hartogs' sense to $S_{1} \wedge S_{2,m}$. So by Proposition \ref{prop_hartogs_sqp}, we have that:
  $$\U_{S_{1}\wedge S_{2,m}}(R)= \U_{S_{1}}(S_{2,m}\wedge R) +  \U_{S_{2,m}}(\omega^{r_1}\wedge R) -\U_{S_{1}}(S_{2,m}\wedge \omega^{k-r_1-r_2+1}),$$ 
where the super-potentials are of mean $0$. Similarly, we let $m\to \infty$. Recall that $S_{2,m}\wedge R$ converges to $S_2\wedge R$ in the Hartogs' sense, hence $\U_{S_{1}}(S_{2,m}\wedge R)=\U_{S_{2,m}\wedge R}(S_{1})$ converges to $\U_{S_{2}\wedge R}(S_{1})$ . So we have indeed:
$$\U_{S_{1}\wedge S_{2}}(R)= \U_{S_{1}}(S_{2}\wedge R) +  \U_{S_{2}}(\omega^{r_1}\wedge R) -\U_{S_{1}}(S_{2}\wedge \omega^{k-r_1-r_2+1}).$$
 Since $\U_{S_{1}}(S_{2}\wedge \omega^{k-r_1-r_2+1})$ does not depend on $R$ and is finite because $S_1$ and $S_2$ are wedgeable, we can add it to $\U_{S_{1}\wedge S_{2}}$ and we have the lemma. \hfill $\Box$ \hfill 

\begin{lemme}\label{pushpull}
Let $T_1 \in \C_{r_1}$ be an $f_*$-admissible current with $r_1 \geq k-s$.
 Let $T_2\in \C_{r_2}$ be an $f^*$-admissible current with $r_1+r_2\leq k$ such that $L(T_2)$ and $T_1$ are wedgeable and $L(T_2)\wedge T_1$ is $f_*$-admissible. Assume also that $T_2$ and $\Lambda(T_1)$ are wedgeable. Then:
 $$\Lambda(L(T_2)\wedge T_1)=T_2 \wedge \Lambda(T_1). $$ 
\end{lemme}
\emph{Proof.} Assume first that $T_2$ is smooth. Let $T_{1,n}$ and $L_{2,m}$ be sequences in $\C_{r_1}$ and $\C_{r_2}$ converging in the Hartogs' sense to $T_1$ and $L(T_2)$. Let $\Theta$ be a smooth current of bidegree $k-r_1-r_2$. When $n$ and $m$ goes to $\infty$, $\Lambda(L_{2,m}\wedge T_{1,n})$ converge to $\Lambda(L(T_2)\wedge T_1)$ in the sense of currents by Propositions \ref{prop_cv_wedge_hartogs} and \ref{continuiteL}. We want to show that:
$$ \langle T_2\wedge\Lambda(T_1),  \Theta \rangle = \langle \Lambda(L(T_2)\wedge T_1), \Theta \rangle$$
for all $\Theta$ smooth. \\

First assume that $\Theta$ is closed and (strongly) positive. Up to a multiplicative constant, we assume that $\Theta \in \C_{k-r_1-r_2}$. Since everything is smooth:
\begin{align*}
\langle \Lambda(L_{2,m}\wedge T_{1,n}), \Theta \rangle & = \langle L_{2,m}\wedge T_{1,n}, L(\Theta) \rangle \\
																											 & = \langle  T_{1,n}, L_{2,m}\wedge L(\Theta) \rangle.
\end{align*}
 Since $L_{2,m}$ converges to $L(T_2)$, we have that $L_{2,m}\wedge L(\Theta)$ converges to $L(T_2\wedge \Theta)$ in the sense of currents. Indeed, the sequence $(L_{2,m}\wedge L(\Theta))_m$ is of mass $1$. We can extract a converging subsequence (in the sense of currents). 
Observe that its limit is less than $ L(T_2)\wedge L(\omega^{k-r_1-r_2})$ which gives no mass to
$I^+$ by dimension's arguments. So its limit gives no mass to $I^+$ either. But outside $I^+$, 
  $ L_{2,m}\wedge L(\Theta)$ converges to the smooth form $ L(T_2) \wedge L(\Theta)$.
   That implies that $ L_{2,m}\wedge L(\Theta)$ converges to the trivial extension of
    $L(T_2)\wedge L(\Theta)$ which is equal to the form $L(T_2 \wedge \Theta)$ which has coefficients in $L^1$ (as in Lemma \ref{continuite} test 
    the convergence against a smooth form $\Psi$ and write it $\xi \Psi+ (1-\xi) \Psi$ 
    where $\xi$ is a cut-off function equal to $1$ in a small neighborhood of $I^+$).

  So, letting $m \to \infty$ and using the fact that $T_2\wedge \Theta$ is smooth:
 \begin{align*}
\langle  T_{1,n}, L(T_2\wedge \Theta) \rangle & = \langle \Lambda(T_{1,n}), T_2\wedge \Theta \rangle \\
																						  & = \langle T_2\wedge\Lambda(T_{1,n}),  \Theta \rangle .
\end{align*}
Now, we let $n\to \infty$, $\Lambda(T_{1,n})$ converges to $\Lambda(T_1)$ in the Hartogs' sense (Proposition \ref{continuiteL}) hence Proposition \ref{prop_cv_wedge_hartogs} gives that $T_2\wedge\Lambda(T_{1,n})$ converges to $T_2\wedge\Lambda(T_1)$ in the sense of currents. So we have indeed that: 
$$ \langle T_2\wedge\Lambda(T_1),  \Theta \rangle = \langle \Lambda(L(T_2)\wedge T_1), \Theta \rangle$$
for $\Theta$ closed. \\

Now, for $\Theta$ not necessarily closed, we can assume that $\Theta$ is positive and $\Theta \leq C \omega^{k-r_1-r_2}$ for $C$ large enough. Again, we have that
\begin{align*}
\langle \Lambda(L_{2,m}\wedge T_{1,n}), \Theta \rangle   = \langle  T_{1,n}, L_{2,m}\wedge L(\Theta) \rangle.
\end{align*}
The positive current  $ L_{2,m}\wedge L(\Theta)$ is less than $C L_{2,m}\wedge L(\omega^{k-r_1-r_2})$ 
so it is of mass less than $C$. We can extract a converging subsequence (in the sense of currents). 
Observe that its limit is less than $C L(T_2)\wedge L(\omega^{k-r_1-r_2})=CL(T_2\wedge \omega^{k-r_1-r_2})$ which gives no mass to 
$I^+$ by dimension's arguments. So its limit gives no mass to $I^+$ either. Again outside $I^+$, 
  $ L_{2,m}\wedge L(\Theta)$ converges to the smooth form $ L(T_2) \wedge L(\Theta)$.
   That implies that $ L_{2,m}\wedge L(\Theta)$ converges to the trivial extension of
    $L(T_2)\wedge L(\Theta)$ which is equal to the form $L(T_2 \wedge \Theta)$ which has coefficients in $L^1$. We have again that:
$$ \langle  T_{1,n}, L(T_2\wedge \Theta) \rangle  = \langle T_2\wedge\Lambda(T_{1,n}),  \Theta \rangle .$$ 
That gives the conclusion as before.\\

Now, for $T_2$ not necessarily smooth, we can approximate $T_2$ by a sequence of smooth currents converging in the Hartogs' sense to $T_2$.  Since both members of the equality:
$$\Lambda(L(T_2)\wedge T_1)=T_2 \wedge \Lambda(T_1),$$
depend continuously on $T_2$ for the Hartogs' convergence (wedge-product, pull-pack and push-forward are continuous for the Hartogs' convergence) we get the lemma from the smooth case. \hfill $\Box$ \hfill \\

Some of the hypothesis of the following lemma are not necessary, but the following version is enough for our purpose: 
\begin{lemme}\label{decpushpull}
Let $S_1$, $S_2$ and $S_3$ in $\C_{r_1}$,  $\C_{r_2}$ and  $\C_{r_3}$ with $r_1+r_2+r_3=k+1$. Assume that $S_2$ is smooth and that $L(S_2)$ and $S_3$ are wedgeable. Assume that $L(S_2)\wedge S_3$ is $f_*$-admissible. Assume also that the super-potential $\U_{S_1}$ of $S_1$ is finite at $\Lambda(L(S_2)\wedge S_3)$. Finally, we also assume that $S_1$ is $f^*$-admissible, that $S_3$ and $L(S_1)$ are wedgeable and that their wedge product is finite at the super-potential $\U_{L(S_2)}$ of $L(S_2)$. Then we have the formula:
\begin{align*}
\U_{S_1}(\Lambda(L(S_2)\wedge S_3))=&(\frac{\lambda_{r_1}}{\lambda_{r_1-1}})\left( \U_{L(S_2)}(S_3 \wedge L(S_1)) - \U_{L(S_2)}(S_3 \wedge L(\omega^{r_1}))\right) \\
                                        &+\U_{S_1}(\Lambda(\omega^{r_2}\wedge S_3))
\end{align*} 
\end{lemme} 
\emph{Proof.} First, observe that $\omega^{r_1}$ is more H-regular than $S_1$ hence $L(\omega^{r_1})$ is more H-regular than $L(S_1)$. So, we have that $S_3$ and $L(\omega^{r_1})$ are wedgeable and $S_3 \wedge L(\omega^{r_1})$ is more H-regular than $S_3 \wedge L(S_1)$. In particular, $\U_{L(S_2)}(S_3 \wedge L(\omega^{r_1}))$ is finite. Similarly, the expression $\U_{S_1}(\Lambda(\omega^{r_2}\wedge S_3))$ is finite and everything is well defined in it.

 Let $S_{1,{m_1}}$, $L_{2,m_2}$ and $S_{3,m_3}$ be sequences of smooth currents converging in the Hartogs' sense to  $S_1$, $L(S_2)$ and $S_3$.  Let $U_{1,m_1}$ and $U_{2,m_2}$ be smooth quasi-potential of $S_{1,m_1}$ and $L_{2,m_2}$. For smooth currents, we have the identity:
\begin{align*}
\U_{S_{1,m_1}}(\Lambda(L_{2,m_2}\wedge S_{3,m_3}))&= \langle   U_{1,m_1}, \Lambda(L_{2,m_2}\wedge S_{3,m_3})   \rangle \\
                                              &= \langle   (\lambda_{r_1-1}(f))^{-1} f^*(U_{1,m_1}), L_{2,m_2}\wedge S_{3,m_3}   \rangle.
\end{align*} 
By Stokes, we recognize:
$$\langle   (\lambda_{r_1-1}(f))^{-1} f^*(U_{1,m_1}), \omega^{r_2}\wedge S_{3,m_3} \rangle + \langle dd^c((\lambda_{r_1-1}(f))^{-1} f^*(U_{1,m_1})), U_{2,m_2} \wedge S_{3,m_3} \rangle. $$
Since $f^*$ commutes with $dd^c$, we have that 
\begin{align*}
dd^c((\lambda_{r_1-1}(f))^{-1} f^*(U_{1,m_1}))&=(\lambda_{r_1-1}(f))^{-1} f^*(S_{1,m_1}-\omega^{r_1}) \\
                                                &= \left(\frac{\lambda_{r_1}}{\lambda_{r_1-1}}\right) \left(L(S_{1,m_1})-L(\omega^{r_1}) \right).
\end{align*}  
The lemma follows then by letting $m_1$ then $m_2$ then $m_3$ go to $\infty$ and using the continuity of the wedge product, the pull-back, push-forward and value at a point for the super-potential for the Hartogs' convergence. \hfill $\Box$ \hfill \\
  
We also have the following integration by parts lemma:
\begin{lemme}\label{byparts}
Let $S_1$, $S_2$ and $S_3$ in $\C_{r_1}$,  $\C_{r_2}$ and  $\C_{r_3}$ with $r_1+r_2+r_3=k+1$. Assume that the $S_i$ are two by two wedgeable. Then if $\U_{S_1}$ and $\U_{S_2}$ are super-potentials of $S_1$ and $S_2$ finite at $S_2 \wedge S_3$ and $S_1\wedge S_3$:
\begin{align*}
\U_{S_1}(S_2\wedge S_3)-\U_{S_1}(\omega^{r_2}\wedge S_3)=\U_{S_2}(S_1\wedge S_3)-\U_{S_2}(\omega^{r_1}\wedge S_3)
\end{align*} 
\end{lemme} 
\emph{Proof.} First observe that $\omega^{r_2}\wedge S_3$ and $\omega^{r_1}\wedge S_3$ are more H-regular than $S_2\wedge S_3$ and $S_1 \wedge S_3$ so every term is finite. 

Now, if every term is smooth, we write $U_{S_1}$ and $U_{S_2}$ quasi-potentials of $S_1$ and $S_2$. By Stokes:
\begin{align*}
\U_{S_1}(S_2\wedge S_3)-\U_{S_1}(\omega^{r_2}\wedge S_3)&= \langle U_{S_1}, S_2\wedge S_3-\omega^{r_2}\wedge S_3\rangle=\langle U_{S_1}, dd^c U_2\wedge S_3\rangle \\
                                                        &=\langle dd^c U_{S_1},  U_2\wedge S_3\rangle = \langle S_1,  U_2\wedge S_3\rangle- \langle \omega^{r_1},  U_2\wedge S_3\rangle \\
                                                        &=\U_{S_2}(S_1\wedge S_3)-\U_{S_2}(\omega^{r_1}\wedge S_3).
\end{align*} 
And the result follows in the general case by Hartogs' convergence. \hfill $\Box$ \hfill \\

We also have the following refinement of Lemma \ref{lemma_wedge_regular} whose proof is similar:
\begin{lemme}\label{lemma_wedge_regular2} 
Let $R_{1,n}$ and $R_{2,m}$ be sequence of currents in $\C_{p_1}$ and $\C_{p_2}$ converging in the Hartogs' sense to
 $R_1$ and $R_2$ which are wedgeable. Then there exists a constant $A_{n,m} >0$ such that 
 $$\U_{R_{1,n}\wedge R_{2,m}} \geq \U_{R_{1}\wedge R_{2}}-A_{n,m}$$ 
where the super-potentials are of mean 0 and where $A_{n,m}$ is uniformly bounded from above in $n$ and $m$ and is arbitrarily small for $n$ and $m$ large enough. 
\end{lemme}
\emph{Proof.}
By Lemma \ref{lemma_wedge_regular}, $R_{1,n}$ and $R_{2,m}$ are
wedgeable. 

The symbols $U$ and $\U$ below denote quasi-potentials and 
super-potentials of mean 0. Assume first that all the terms are smooth.
By hypothesis, there is a constant $a$ such that
$\U_{R_{1,n}}+a\geq \U_{R_1}$ and $\U_{R_{2,m}}+a\geq \U_{R_2}$.
Write $r=k-p_1-p_2+1$. Consider a smooth form $R$ in
$\C_r$ and choose $U_R$ smooth. 
We have the computation:
\begin{eqnarray*}
\U_{R_{1,n}\wedge R_{2,m}}(R) & = &  \langle R_{1,n}\wedge R_{2,m},U_R\rangle =
 \langle R_{2,m},\omega^{p_1}\wedge U_R\rangle + \langle
 R_{1,n}-\omega^{p_1},R_{2,m}\wedge U_R\rangle\\
& = &   \langle R_{2,m},\omega^{p_1}\wedge U_R\rangle + \langle
 dd^c  U_{R_{1,n}},R_{2,m}\wedge U_R\rangle\\
& = &   \langle R_{2,m},\omega^{p_1}\wedge U_R\rangle + \langle
  U_{R_{1,n}},R_{2,m}\wedge dd^c U_R\rangle\\
&=& \langle R_{2,m},\omega^{p_1}\wedge U_R\rangle +
\U_{R_{1,n}}(R_{2,m}\wedge R)-\U_{R_{1,n}}(R_{2,m}\wedge\omega^{r}). \\
&=& \U_R( R_{2,m}\wedge \omega^{p_1}) +
\U_{R_{1,n}}(R_{2,m}\wedge R)-\U_{R_{1,n}}(R_{2,m}\wedge\omega^{r}).
\end{eqnarray*}
And that identity holds when the currents are not smooth by Hartogs' convergence. We have the same identity for $R_1\wedge R_{2,m}$ and $R_1\wedge R_2$. By difference, we have:
\begin{align*}
\U_{R_{1,n}\wedge R_{2,m}}(R)-\U_{R_{1}\wedge R_{2,m}}(R)+\U_{R_{1}\wedge R_{2,m}}(R)-\U_{R_{1}\wedge R_{2}}(R) = \\
\U_{R_{1,n}}(R_{2,m}\wedge R)-\U_{R_{1}}(R_{2,m}\wedge R)-\U_{R_{1,n}}(R_{2,m}\wedge\omega^{r})+\U_{R_{1}}(R_{2,m}\wedge\omega^{r}) \\
+\U_{R_{2,m}}(R_{1}\wedge R)-\U_{R_{2}}(R_{1}\wedge R)-\U_{R_{2,m}}(R_{1}\wedge\omega^{r})+\U_{R_{2}}(R_{1}\wedge\omega^{r}). 
\end{align*}
 So:
 \begin{align*}
\U_{R_{1,n}\wedge R_{2,m}}(R)-\U_{R_{1}\wedge R_{2}}(R) \geq \\
-2a-\U_{R_{1,n}}(R_{2,m}\wedge\omega^{r})+\U_{R_{1}}(R_{2,m}\wedge\omega^{r})-\U_{R_{2,m}}(R_{1}\wedge\omega^{r})+\U_{R_{2}}(R_{1}\wedge\omega^{r}). 
\end{align*}
The last quantity does not depend on $R$ and is uniformly bounded from below: the terms with a minus sign are greater than $-M$ since the super-potentials are of mean $0$, and since $R_{2,m}\wedge\omega^{r}$ converges to $R_{2}\wedge\omega^{r}$ in the Hartogs' sense and $\U_{R_{1}}(R_{2}\wedge\omega^{r})$ is finite, we have that $\U_{R_{1}}(R_{2,m}\wedge\omega^{r})$ and $\U_{R_{2}}(R_{1}\wedge\omega^{r})$ are uniformly bounded from below.  \\

This gives that the constant $A_{n,m}$ of the lemma is uniformly bounded from above in $n$ and $m$. Now, we can choose $A_{n,m}$ going to zero by Proposition \ref{prop_cv_wedge_hartogs}: if not, we can extract subsequences such that $A_{n_i,m_i} \geq \varepsilon >0$ and it contradicts the fact that $R_{1,n_i}\wedge R_{2,m_i}$  converges in the Hartogs' sense to $R_{1}\wedge R_{2}$. \hfill $\Box$ \hfill

\bigskip

\noindent Henry de Thélin, Mathématiques - Bât. 425, UMR 8628,\\
Université Paris-Sud, 91405 Orsay, France.  \\
\noindent Email: Henry.De-Thelin@math.u-psud.fr

\vspace{1cm} 

\noindent Gabriel Vigny, Mathématiques - Bât. 425, UMR 8628,\\ 
Université Paris-Sud, 91405 Orsay, France.  \\
\noindent Email: gabriel.vigny@math.u-psud.fr

 \end{document}